%% file: HAND.tex
\documentclass[12pt]{amsart} 
\usepackage{amssymb,amsmath,amscd,a4,graphicx,latexsym,psfrag,epsfig}

\usepackage[all]{xy}
\everymath{\displaystyle}

\newtheorem{teo}{Theorem}[section]
\newtheorem{lem}[teo]{Lemma}
\newtheorem{cor}[teo]{Corollary}
\newtheorem{exa}[teo]{Example}
\newtheorem{prop}[teo]{Proposition}
\newtheorem{defi}[teo]{Definition}

\newtheorem{remark}[teo]{Remark}

\newcommand{\Pm}{\mathbb{P}}
\newcommand{\mr}{\mathbb{R}}
\newcommand{\mc}{\mathbb{C}}
\newcommand{\mz}{\mathbb{Z}}
\newcommand{\mh}{\mathbb{H}}

\newcommand{\mm}{\mathbb{M}}

\newcommand{\mx}{\mathbb{X}}

\newcommand{\md}{\mathbb{D}}

\newcommand{\Yy}{{\mathcal Y}}
\newcommand{\Aa}{{\mathcal A}}
\newcommand{\Bb}{{\mathcal B}}
\newcommand{\Cc}{{\mathcal C}}
\newcommand{\Dd}{{\mathcal D}}
\newcommand{\Ee}{{\mathcal E}}
\newcommand{\Ff}{{\mathcal F}}
\newcommand{\Gg}{{\mathcal G}}
\newcommand{\Hh}{{\mathcal H}}

\newcommand{\Kk}{{\mathcal K}}
\newcommand{\Ll}{{\mathcal L}}
\newcommand{\Mm}{{\mathcal M}}

\newcommand{\Pp}{{\mathcal P}}
\newcommand{\Qq}{{\mathcal Q}}

\newcommand{\Ss}{{\mathcal S}}
\newcommand{\Tt}{{\mathcal T}}
\newcommand{\Uu}{{\mathcal U}}
\newcommand{\Ww}{{\mathcal W}}
\newcommand{\Vv}{{\mathcal V}}

\newcommand{\SG}{{\mathfrak S}}

\newcommand{\cG}{{\mathfrak c}}

\newcommand{\sG}{\mathfrak s}
\newcommand{\lG}{\mathfrak l}
\newcommand{\pG}{\mathfrak p}
\newcommand{\mG}{\mathfrak m}
\newcommand{\CG}{{\mathfrak C}}
\newcommand{\tG}{{\mathfrak t}}
\newcommand{\C}{{\mathbb C}}

\newcommand{\R}{{\mathbb R}}

\def\ort#1{#1^\perp}

\def\arctgh{\mathrm{arctgh\, }}

\def\Dim{\emph{Proof : }}
\def\cvd{\nopagebreak\par\rightline{$_\blacksquare$}}

\def\fut{\mathrm{I}^+}                         %futuro
\def\pass{\mathrm{I}^-}                         %passato
\def\SL#1#2{\mathrm{SL}(#1,\,\mathbb{#2})}

\def\eps{\epsilon}                        % scrive epsilon

\def\sign{\mathrm{sign}}                      % segno
\def\d{\mathrm{d}}
\def\ch{\mathrm{ch\,}}
\def\sh{\mathrm{sh\,}}
\def\tgh{\mathrm{tgh\,}}

\def\arctgh{\mathrm{arctgh\, }}

\def\Shat{\hat S}

\title{(2+1) Einstein spacetimes of finite type} 

\author {Riccardo Benedetti and Francesco Bonsante}

\begin{document}

\maketitle

\vspace{0.5cm}

\centerline {Dipartimento di Matematica, Largo B. Pontecorvo 5,
Pisa, ITALY}  

\centerline {benedett@dm.unipi.it}
\smallskip

\centerline {Scuola Normale Superiore, Pisa, ITALY} 

\centerline{f.bonsante@sns.it}

\tableofcontents
 
\input HANDINTRO.tex

\input HANDGRAV.tex
\input HANDMLS.tex
\input HANDWR.tex
\input HANDADS.tex

\input HANDPART.tex

\input HANDBIB.tex
\end{document}

%% file: HANDINTRO.tex
\section{Introduction}
A surface $S$ is said to be of {\it finite type} if: 
\smallskip

(1) it is of the form 
$$S= \hat S \setminus V$$ where $(\hat S, V)$ is a compact closed
oriented surface of genus $g\geq 0$, with a set of $r\geq 0$ marked
points $V=\{p_1,\cdots,p_r\}$;
\smallskip

(2) the fundamental group of $S$ is {\it
  non Abelian}, equivalently $2-2g -r < 0$.
\medskip

The main aim of this survey is to widely describe, for every $S$ of
finite type, and for every $\kappa = 0,\ \pm 1$, the geometry of
3-dimensional {\it maximal globally hyperbolic} Lorentzian spacetimes
of {\it constant curvature $\kappa$}, that contain a {\it complete}
Cauchy surface homeomorphic to $S$. We call them generically {\it
Einstein $MGH$ spacetimes of finite type}. The (3-dimensional) general
relativity background will be briefly recalled in Section
\ref{grav}. These spacetimes are supported by the product $S\times
\mr$; considered up to Lorentzian isometry homotopic to the identity
of $S\times \mr$, they form, for every $\kappa$, a {\it
Teichm\"uller-like space} denoted
$$\Mm\Gg\Hh_\kappa(S) \ .$$ Clearly these notions make sense also if
$S$ is not necessarily of finite type. In the monograph \cite{Be-Bo}
we have developed a {\it canonical Wick rotation-rescaling} theory on
such general $MGH$ spacetimes. It is easy to see that
$\Mm\Gg\Hh_\kappa(S) \neq \emptyset \ $ for every $\kappa$, if and only
if the open disk $D^2$ is the universal covering of $S$. In
\cite{Be-Bo} we have actually analyzed $\Mm\Hh\Gg_\kappa(D^2)$, by
developing also an {\it equivariant version} of the theory, with
respect to the action of {\it any} deck transformation
group. WR-rescaling theory includes a wide generalization of Mess
classification \cite{M} (completed by Scannell \cite{Sc} for $\kappa
=1$) of $MGH$ spacetimes with {\it compact} Cauchy surfaces ({\it
i.e.}  $V = \emptyset$). Moreover, it establishes explicit geometric
correlations between spacetimes of different curvatures, or between
spacetimes and {\it complex projective structures} on $S$. In
particular, this gives a clear geometric explanation of the occurrence
of a certain {\it ``universal'' parameter space}
$$\Mm\Ll(S)$$ shared by all $\Mm\Gg\Hh_\kappa(S)$, $\kappa = 0,\ \pm 1$,
and by $\Pp(S)$, the Teichm\"uller-like space of complex projective
structures on $S$.

A large part of this survey just reports on such a theory, by
specializing it to $S$ of finite type. This class is large enough to
display the main features of the theory; on the other hand,
spacetimes of finite type are possibly easier to figure out than
completely general ones. In fact we will spell out several specific
statements that are quite implicit in the general treatment given in
\cite{Be-Bo}. Hence the present paper should be an actual complement
to that monograph. Moreover, there are in this case direct relations
between $\Mm\Ll(S)$ and more familiar {\it Teichm\"uller spaces of
hyperbolic structures on $S$} and, to some extent, with the
corresponding {\it tangent bundles} (see Section \ref{ML(S)}). For
example, when $S$ is compact $\Mm\Ll(S)$ coincides with the
(topologically trivialized) bundle $\Tt_g\times \Mm\Ll_g$ of {\it
measured geodesic laminations} on hyperbolic structures on $S$; in
general we will deal with hyperbolic structures $F$ on $S$ whose
completions $F^\Cc$ have (non necessarily compact) geodesic boundary,
and with a kind of measured geodesic laminations $\lambda$ on
$F^\Cc$. In fact, another goal should be to convince a reader familiar
with such topics of hyperbolic geometry, that not only this provides
some important tools for studying Einstein spacetimes; in the reverse
direction, via Lorentzian geometry we get a new insight into several
fundamental hyperbolic constructions such as {\it grafting},
(3-dimensional hyperbolic) {\it bending} and {\it earthquakes} along
laminations. To support this claim we just mention here the ``AdS
proof'' of {\it Thurston Earthquakes Theorem} for hyperbolic
structures on compact surfaces $S$, that Mess obtained in \cite{M} as
a by-product of his classification of spacetimes in
$\Mm\Gg\Hh_{-1}(S)$. An AdS look at earthquakes theory {\it beyond the
compact case} will be a theme of this paper (see Section 
\ref{moreAdS}).
\smallskip

Finally, we note that spacetimes of finite type occur (via canonical
Wick rotation) as ``ending spacetimes'' of geometrically finite
hyperbolic $3$-manifolds, which furnish basic examples for a bordism
category supporting (2+1) QFT pertinent to 3-dimensional gravity (see
Section 1.11 of \cite{Be-Bo}, and \cite{BB}).

\smallskip

In Section \ref{moreAdS} we focus the AdS case that displays the
richer phenomenology, mostly referring (besides \cite{Be-Bo}) to
\cite{Ba}(1, 2) and \cite{BSK}.  In particular we will describe the
common maximal {\it causal extension} $\Omega(h)$ of the $MGH$
spacetimes of finite type that share a given AdS holonomy $h$.  We
will see that $\Omega(h)$ is still supported by the product $S\times
\R$ but it is not in general globally hyperbolic. This is a
particularly interesting case, because we can detect a 
specific one among the maximal globally hyperbolic spacetimes
contained in $\Omega(h)$ that can be truly considered as a {\it black hole}. 
The analysis of the causal extension is also important to achieve a
proof of the Earthquake Theorem.
\smallskip

Finally, in Section \ref{part} we will outline (by following \cite{BG}(2) and
mostly \cite{BS}) how the Wick rotation - rescaling theory (partially)
extends to $MGH$ spacetimes of finite type that include world lines of
``particles'' ({\it i.e} inextensible time-like lines of space-like
{\it conical singularities}).
\medskip
  
We stress that this paper is not intended to be exhaustive of the
subject. We have made a few partial and subjective choices, organized
around our favorite Wick rotation-rescaling view point. Nevertheless,
we hope that this would be enough to show that 3-dimensional gravity
is a fairly non-trivial and beautiful ``toy model''. In particular, we
have neglect a classical analytic approach to the classification of
constant curvature $MGH$ spacetimes in terms of solutions of the
Gauss-Codazzi equation at a Cauchy surface, possibly imposing some
supplementary conditions to such solutions, that translate some
geometric property of the embedding of $S$ as Cauchy surface (see also
Section \ref{grav}). A widely studied possibility requires that the
surface has {\it constant mean curvature} (see for instance \cite{Mon,
A-M-T, B-Z, Kra-Sch}). At least for compact $S$, the classical
Teichm\"uller space of {\it conformal} structures on $S$, with its
complex cotangent bundle arises in this way towards the classification.
This approach also selects a distinguished {\it CMC global time} on
$MGH$ spacetimes, that basically coincides with the mean curvature of
its level surfaces.
\smallskip

Wick rotation-rescaling theory is based on a rather different more
geometric approach, initiated by Mess in \cite{M}. It turns out that a
key ingredient is another canonical time, the so called {\it
cosmological time}.  Every $MGH$ spacetime is in a sense determined by
the ``asymptotic states'' of the corresponding level surfaces, rather
than the embedding data of some Cauchy surface. The Wick
rotation-rescaling mechanism is ultimately based on the fact that
$MGH$ spacetimes (of different curvatures) can be associated in such a
way that the intrinsic geometry of these level surfaces does not
depend on the curvature, up to some scale factor.

%% file: HANDGRAV.tex
\section{3-dimensional gravity}\label{grav}
\subsection{General background}\label{gen-rel}
For the basic notions of global Lorentzian geometry and causality we
refer for instance to \cite{BEE, H-E}.  

A $(n+1)$ {\it spacetime} consists of a $(n+1)$-manifold $M$ equipped
with a Lorentzian metric $h$ and with a {\it time orientation}, so
that the {\it causal past/future} of every {\it event} $p\in (M,h)$ is
determined. We also stipulate that $M$ is oriented.

Roughly speaking, the general problem of gravity can be stated as
follows. Given a $(n+1)$-manifold $M$, a symmetric $(0,2)$-Tensor
$T$ on $M$ and a constant $\Lambda$ (called the \emph{cosmological
constant}), find out all spacetimes $(M,h)$ such that:
\smallskip

(a) The metric $h$ satisfies the {\it Einstein equation}
\[
   Ric_h+ (\Lambda-1/2 R_h) h=T
\]
where $Ric_h$ is the Ricci tensor of $h$, $R_h$ is the scalar
curvature.
\smallskip

(b) The global {\it causal structure} of $(M,h)$ satisfies determined
conditions. 
\medskip

These spacetime structures are considered up to diffeomorphism of $M$
that preserves the tensor $T$. Both the features of the tensor $T$ and
of the causality conditions are determined by physical (even logical)
considerations. Normally they also impose some constraints on the
topology of $M$.  Requirements in (a) and (b) are basically of
independent nature.
\smallskip

The \emph{pure gravity} case is when $T=0$. In such a case the
solutions of Einstein equation coincide with the so called {\it
Einstein metrics}: $Ric_h = \frac{2}{n-1}\Lambda h$.
\smallskip

The basic causality condition is that $(M,h)$ is {\it chronological}
({\it causal}), that is it does not contain any {\it closed} timelike
(causal) curve $c$. A curve is said timelike (causal) if the velocity
field $v(t)$ along $c$ is timelike (causal): $h(v(t),v(t))<0$ (
$h(v(t),v(t)) \leq 0$ ).

The strongest causality condition is that $(M,h)$ contains a {\it
Cauchy surface} $S$; this means that $S$ is a {\it spacelike}
hypersurface of $M$ (the restriction of $h$ to $S$ is Riemannian),
such that every causal inextensible line of $(M,h)$ intersects $S$
exactly once.  In such a case we say that $(M,h)$ is \emph{globally
hyperbolic}. If $(M,h)$ is globally hyperbolic the $M$ turns to be a
product manifold $M \cong S\times \mr$ so that (up to some
diffeomorphisms of $M$) the Cauchy surface $S$ coincides with $S\times
\{0\}$, and every slice $S\times\{t\}$ is $h$-spacelike (indeed we can
also require that every such a slice eventually is a Cauchy surface of
$(M,h)$). Such a picture is coherent with the intuitive idea of {\it a
space evolving in time}. Globally hyperbolic spacetimes naturally
arise as {\it dependence domains} $(D(S),h_{|D(S)})$ of spacelike
hypersurfaces $S$ in arbitrary spacetimes $(M,h)$; $S$ turns to be a
Cauchy surface of $D(S)$.  Hence globally hyperbolic spacetimes make a
fundamental sector of gravity theory.

\subsection{$(2+1)$-spacetimes}\label{2+1ST} 
3D gravity is much simpler (but non trivial) because in three
dimensions the Riemann tensor is determined by the Ricci tensor. In
particular 3D Einstein metrics actually have {\it constant (sectional)
curvature}. The sign of the curvature coincides with the sign of the
cosmological constant.  We will be mainly concerned by $(2+1)$
globally hyperbolic Einstein spacetimes $(M,h)$ ({\it i.e.} of constant
curvature $\kappa$).

We recall two possible ways to study such spacetimes.  The first
\emph{analytic} one is based on the important fact that the germ of
the metric $h$ at a Cauchy surface $S$ determines, in a sense that we
will make precise, the whole spacetime. This leads to consider the
pairs $(g,b)$ of a Riemannian metric, $g$, on the surface $S$ and a
$g$-symmetric endomorphism, $b$, of $TS$, that verify the
Gauss-Codazzi equation
\[
\begin{array}{l}
d^\nabla b=0\\
\det b = \kappa -\kappa_g
\end{array}
\]
where $d^\nabla$ is the differential with respect to the Levi-Civita
connection of $g$, $\kappa$ is a constant and $\kappa_g$ is the Gauss
curvature of $g$.

It is possible to associate to such a pair $(g,b)$ a Lorentzian metric
$h$ on $M=S\times\mr$ of constant curvature $\kappa$, such that
$S=S\times \{0\}$ is a Cauchy surface, the {\it first fundamental
form} of $S$ in $(M,h)$ is $g$ and the {\it shape operator} is $b$
(recall that the shape operator of a spacelike surface $F$ in a
Lorentzian or Riemannian manifold $M$ is the endomorphism of $TF$ that
coincides with the covariant derivative of the normal field of $F$ in
$M$).  A priori $(g,b)$ determines only the germ of $h$ around
$S\times\{0\}$. On the other hand, it is proved in~\cite{Cho-Ge} that
there exists a unique (up to isometries) such a globally hyperbolic
spacetime $(M_{\max},h_{max})$ that is {\it maximal} in the following
sense:

{\it Given any globally hyperbolic spacetime $(M,h)$ as above, there
exists an isometric embedding of $(M,h)\rightarrow(M_{max},h_{max})$
that is the identity on $S\times\{0\}$ (and preserves the
orientations).}
\smallskip

At a first sight, this definition of ``maximality'' involves the choice
of a Cauchy surface (i.e. $S\times\{0\}$). On the other hand, one can
see that it is equivalent to the following one: 

{\it Every isometric embedding of $(M_{max},h_{max})$ into an
Einstein spacetime $(N,k)$ that sends {\rm any} Cauchy surface of
$M_{max}$ onto a Cauchy surface of $N$ actually is a global isometry.}
\smallskip

This last property gives a good definition of the class of {\it
maximal globally hyperbolic (MGH) Einstein spacetimes}, that makes
intrinsic sense, not depending on the analytic approach we are
outlining. It is reasonable to restrict to this class in order to get
a classification.

Continuing with the analytic approach, a well-defined map eventually
associates to every pair $(g,b)$ as above the (isotopy class of the)
maximal globally hyperbolic spacetime
$(M_{max},h_{max})_{(g,b)}$. Such map is surjective, but not
injective. In fact it establishes a bijective correspondence between
pairs $(g,b)$ and spacetimes with a marked Cauchy surface.  To get rid
of this excess of degrees of freedom, some additional condition on
$(g,b)$ has to be imposed, possibly translating some geometric
property of the Cauchy surface embedding.  A widely investigated
possibility consists in requiring that the trace of $b$ is
constant, that is, $S\times\{0\}$ is a surface of {\it constant
mean curvature}.
\medskip

The second \emph{geometric} approach makes use of the $(X,G)$-{\it
structure} technology.  Indeed any $(2+1)$ Einstein spacetime $M$ is a
$(\mx_\kappa, {\rm Isom}(\mx_\kappa))$-{\it manifold}, where
$\mx_\kappa$ is a suitable {\it isotropic model of constant curvature}
$\kappa$.

Denote $\tilde M$ an universal covering of $M$. A very general
``analytic continuation'' procedure allows to associate to every
$(X,G)$-manifold, $M$, a {\it compatible} couple $(d,h)$, where $d$ is a
\emph{developing map}, that is a local isomorphism $$d:\tilde
M\rightarrow X$$ and $h$ is a {\it holonomy representation}
$$h:\pi_1(M)\rightarrow G$$
such that ({\it $\pi_1$-equivariance}):
\[
     d(\gamma x)=h(\gamma)d(x)
\]
(where $\pi_1(M)$ is identified to the covering transformation group
of $\tilde M$).  The developing map is determined up to
post-composition by any element of $G$, whereas the holonomy is determined
up to conjugation by the same element.

Conversely a local diffeomorphism $d:\tilde M\rightarrow X$
equivariant with respect to a representation $h:\pi_1(M)\rightarrow G$
produces a well-defined $(X,G)$-structure on $M$.

In this paper we will mainly focus on this second geometric
approach. For this reason we will briefly recall the principal
features of the isotropic models of constant curvature $\kappa$, that
we will normalize to be $\kappa= 0,1,-1$.
\medskip

{\bf Minkowski space.}  The isotropic model of flat spacetimes,
$\mx_0$, is the Minkowski space, that is $\mr^3$ equipped with the
flat metric $-\d x_0^2+\d x_1^2+\d x_2^2$.  Isometries of $\mx_0$ are
affine transformations whose linear part preserves the Minkowski
product (that is $\mathrm{Isom}(\mx_0)=SO(2,1)\rtimes\mr^3$).  We
consider the time-orientation on $\mx_0$ such that the $x_0$-component
of future-directed timelike vectors is positive. The set of future
directed unit timelike vectors is a hypersurface of $\mx_0$ that
inherits from $\mx_0$ a Riemannian metric. This is the hyperboloid
model of the hyperbolic plane $\mh^2$.  The isometric action of
$SO^+(2,1)$ on it induces an identification between $SO^+(2,1)$ and
$PSL(2,\mr)\cong {\rm Isom^+}(\mh^2)$ (by using also the Poincar\'e
half-plane model of $\mh^2$).  The main advantage of the hyperboloid
model is that geodesics are just obtained by intersecting $\mh^2$ with
timelike planes. In particular the duality between linear planes and
linear straight lines given by the orthogonality relation induces an
identification between the set of geodesics of $\mh^2$ and the set of
un-oriented spacelike directions of $\mx_0$.  The projection of
$\mh^2$ in the projective plane $\Pm(\mr^3)$ is injective and the
image is the set of timelike directions. Notice that in this
projective (Klein) model geodesics are just projective
lines. Moreover, the set of lightlike directions is the boundary of
$\mh^2$ and the end-points of a geodesic $l$ in $\mh^2$ are the two
lightlike directions contained in the plane of $\mx_0$ containing $l$.
\smallskip

By using the $4$-dimensional Minkowski space, in a similar way we get
the different models of the hyperbolic space $\mh^3$.

\medskip

{\bf De Sitter space.} The set of unit spacelike vectors in
$4$-dimensional Minkowski space, is a Lorentzian submanifold of
constant curvature $1$. It is called the de Sitter spacetime and will
be denoted by $\hat\mx_{1}$. The isometric action of $SO(3,1)$ on
$\hat\mx_{1}$ shows that this model is isotropic and that its isometry
group coincides with $SO(3,1)$. Also in this model geodesic are
obtained by intersecting $\hat\mx_1$ with a linear plane of the Minkowski
space. In particular spacelike geodesics are closed with length equal
to $2\pi$, whereas timelike geodesics are embedded lines with infinite
length.

It is often convenient to consider the projection of $\hat\mx_1$ into
the projective space $\Pm(\mr^4)$. Notice that the image, $\mx_1$, is
the set of spacelike directions, that is, it is the exterior of
$\mh^3$ into $\Pm(\mr^4)$.  Clearly the projection
$\hat\mx_1\rightarrow\mx_1$ is a $2-$to$-1$ covering, so $\mx_1$ is not
simply connected. On the other hand, since $\mx_1=\hat\mx_1/\{\pm
Id\}$, and $\{\pm Id\}$ is the center of $SO(3,1)$, also $\mx_1$ is an
isotropic model of the de Sitter geometry. Its isometry group if
$SO(3,1)/\{\pm Id\}\cong SO^+(3,1)$.  An advantage in using this model
is that $\mx_1$ and $\mh^3$ share the same asymptotic boundary and
their isometry groups actually coincide. By means of the duality
between geodesic planes of $\mh^3$ and spacelike directions of
Minkowski space, $\mx_1$ can be regarded as the set of un-oriented
geodesic planes of $\mh^3$.
\medskip

{\bf Anti de Sitter space.}  Consider on $\mr^4$ a scalar product
$\eta$ with signature $(2,2)$, then the set of unit timelike vectors
is a Lorentzian submanifold $\hat\mx_{-1}$ of constant curvature
$-1$. Let $\mr^4$ be identified with the set of $2\times 2$ matrices,
and consider the form $\eta$ such that $\eta(X,X)=-\det X$. The
signature of $\eta$ is $(2,2)$, so an explicit model of $\hat\mx_{-1}$
is $SL(2,\mr)$ equipped with its Killing form.  The isometric action
of $SL(2,\mr)\times SL(2,\mr)$ on $SL(2,\mr)$ by left and right
multiplication shows that $\hat\mx_{-1}$ is isotropic and that its
isometry group is $SL(2,\mr)\times SL(2,\mr)/(-Id,-Id)$.

As in the previous case the projection of $\hat\mx_{-1}$ into the projective
space $\Pm(\mr^4)$, is a $2$-to-$1$ covering map on a open set $\mx_{-1}$ of
$\Pm(\mr^4)$ that is $PSL(2,\mr)$. Since the covering transformations of
$\hat\mx_{-1}\rightarrow\mx_{-1}$ are $\pm Id$ it follows that $\mx_{-1}$
inherits from $\hat\mx_{-1}$ an isotropic Lorentzian metric of constant
curvature $-1$.  The isometry group of $\mx_{-1}$ turns to be
$PSL(2,\mr)\times PSL(2,\mr)$.

Topologically $\mx_{-1}$ is a solid torus and its boundary in $\Pm(\mr^4)$
can be identified with the projective classes of rank $1$ matrices. The Segre
embedding produces a double foliation on $\partial\mx_{-1}$ by projective
lines (actually it induces a product structure
$\partial\mx_{-1}=\Pm^1\times\Pm^1$). Isometries of $\mx_{-1}$ extends on
the boundary: the left multiplication preserves the leaves of the
left foliation, whereas right multiplication preserves the leaves of the right
foliation. Notice that the product structure on the boundary can be
regarded as a conformal Lorentzian structure.

Geodesics and geodesic planes of $\mx_{-1}$ are the intersection of $\mx_{-1}$
with projective lines and projective planes of $\Pm(\mr^4)$. In particular
projective lines contained in $\mx_{-1}$ are timelike geodesic of
length $\pi$, projective lines tangent to the boundary are lightlike lines and
projective lines intersecting the boundary in two points are spacelike
geodesic of infinite length. Notice that spacelike geodesics are determined by
their end-points on the boundary. Conversely given two points on the boundary
that do not lie on the same left nor right leaf, there exists a unique
spacelike geodesic connecting them.

Projective planes intersecting $\mx_{-1}$ along compression disks are
spacelike planes and turn to be isometric to $\mh^2$.  Points of $\mx_{-1}$
bijectively corresponds to spacelike planes via the duality induced by $\eta$
between points of $\Pm(\mr^4)$ and projective planes. There is a geometric
interpretation of such a duality: given a point $x\in\mx_{-1}$, its dual plane
$P(x)$ is the set of points at distance $\pi/2$ from $p$ along some timelike
geodesic. Conversely given a spacelike plane, its normal geodesics intersect
at the dual point of the plane. Given a spacelike geodesic line $l$, the points
$x$ such that $l\subset P(x)$ form another spacelike line $l^*$, that is the
{\it dual} geodesic of $l$.

%%% Local Variables: 
%%% mode: latex
%%% TeX-master: "HAND"
%%% End: 

%% file: HANDMLS.tex
\section{The space $\Mm\Ll(S)$}\label{ML(S)}  
This section is entirely settled in the framework of (2-dimensional)
{\it hyperbolic geometry}, and several facts that we are going to
recall are well-known. However, we will give later a new insight (if
not an outline of foundation) to many constructions and concepts in
terms of {\it Lorentzian geometry}.
\medskip

Let us fix once for ever some {\it base} surfaces that will support
several geometric structures:
$$(\hat S,V)$$ is a compact closed oriented surface of genus $g\geq
0$, with a set of $r\geq 0$ {\it marked points}
$V=\{p_1,\cdots,p_r\}$;
$$S= \hat S \setminus V \ .$$ $ \overline {\Sigma}$ is obtained by
removing from $\hat S$ a small open disk around each point
$p_j$. Hence $\overline \Sigma$ is compact with $r$ boundary
components $C_1,\cdots,C_r$. We denote by $\Sigma$ the interior of
$\overline \Sigma$.  We fix also a continuous map
$$\overline \phi: \overline \Sigma \to \Shat$$ such that for every
$j$, $\overline \phi(C_j) = p_j$, and the restriction $\phi: \Sigma
\to S$ is an oriented diffeomorphism that is the identity outside a
regular neighbourhood of the boundary of $\overline \Sigma$. In this
way, we will often tacitly confuse $S$ and $ \Sigma$.  We will also
assume that $S$ is {\it not elementary}, that is its fundamental group
is {\it non-Abelian}, equivalently $2-2g -r <0$. Such an $S$ is said
to be {\it of finite type}.

\subsection{The Teichm\"uller space $\widetilde \Tt(S)$}\label{teich} 
\smallskip

We denote by $$\widetilde \Hh(S)$$ the space of non-necessarily
complete hyperbolic structures $F$ on $S$ such that its completion
$F^\Cc$ is a {\it complete hyperbolic surface with geodesic
boundary}. Note that we do not require that the boundary components of
$F^\Cc$ are closed geodesics.  Denote by Diff$^0$ the group of
diffeomorphisms of $S$ homotopic to the identity. Set
$$\widetilde \Tt(S) = \widetilde \Hh(S)/{\rm Diff}^0$$ in other words,
two hyperbolic structures in $\widetilde \Hh(S)$ are identified up to
isometries homotopic to the identity. This is the ``full''
Teichm\"uller space we will deal with.

\subsection{The convex-core map}\label{convex-core}
\smallskip

Let us point out some distinguished subspaces of $\widetilde \Tt(S)$.
$$\Hh(S)\subset \widetilde \Hh(S)$$ denotes the space of {\it
complete} hyperbolic structures on $S$ ({\it i.e.} $F=F^\Cc$). Hence
every $F\in \Hh(S)$ can be realized as the quotient $\mh^2/\Gamma$ by
a discrete, torsion free subgroup $\Gamma \subset {\rm Isom}^+(\mh^2)\cong
PSL(2,\mr)$, isomorphic to $\pi_1(S)$. The corresponding quotient
space
$$\Tt(S)\subset \widetilde \Tt(S)$$ can be identified with the space
of conjugacy classes of such subgroups of $PSL(2,\mr)$.
$$\CG(S)\subset \widetilde \Hh(S)$$ denotes the set of $F$ of {\it
finite area and such that all boundary components of $F^\Cc$ are 
closed geodesics}.
$$\Tt_\cG(S) \subset \widetilde \Tt(S)$$
is the corresponding quotient space.

Clearly if $S$ is compact ($V=\emptyset$)
$$\Tt_g := \Tt_\cG(S) = \Tt(S) = \widetilde \Tt(S)$$
is the classical Teichm\"uller space.

In general, set $$\Tt_{g,r} := \Tt(S)\cap \Tt_\cG(S) \ .$$ Via 
 Uniformization Theorem, $\Tt_{g,r}$ is isomorphic to the
 Teichm\"uller space of {\it conformal structures} on $\hat S$ ({\it
 i.e.} on $S$ that extend to $\hat S$) mod Diff$^0(\hat S,\  {\rm rel}
 \ V)$. $\Tt(S)$ is isomorphic to the Teichm\"uller space of arbitrary
 conformal structures on $S$.
\medskip

\begin{prop}\label{convcore-map} There is an natural isomorphism
$$ \Kk: \Tt(S) \to \Tt_\cG(S) \ .$$
\end{prop}
Basically $\Kk[F]$ coincides with $[\Kk(F)]$, where $\Kk(F)$ denotes
the interior of the {\it convex core} $\ \overline \Kk(F)$ of $F$. Note
that $\Kk(F)^\Cc = \overline \Kk(F)$. This is a bijection because the
convex core determines the whole complete surface. 

\begin{prop}
There is a natural projection $$\beta: \widetilde \Tt(S)\to \Tt(S)$$
such that $\beta_{|\Tt(S)}={\rm Id}$.
\end{prop}
In fact the holonomy of any $[F]\in \widetilde \Tt(S)$ is the
conjugacy class of a faithful representation of $\pi_1(S)$ onto a
discrete, torsion free subgroup $\Gamma$ of $PSL(2,\R)$, hence
$\beta([F])=[\hat F]$, $\hat F=\mh^2/\Gamma$. Finally we can lift the 
map of Proposition \ref{convcore-map} to define the {\it convex-core map}
$$\Kk: \widetilde \Tt(S) \to \Tt_\cG(S), \ \ \Kk([F])= \Kk([\hat F]) \
. $$ In fact we can realize the representatives of the involved
classes in such a way that
$$ \overline\Kk(\hat F) \subset F^\Cc \subset \hat F$$ for $F^\Cc$ is
a closed convex  set in $\hat F$ homotopically equivalent to
$S$, and $\overline \Kk(\hat F)$ is the minimal one with these
properties. In what follows we will often made the abuse of confusing
the classes with their representatives.
\medskip

{\bf Partition by types.}
\begin{prop}\label{typeF} For every complete surface $F\in \Tt(S)$ 
there is a partition
$$V = V_\Pp \cup V_\Hh$$
such that $p$ belongs to $V_\Pp$ iff the following equivalent
properties are satisfied:
\smallskip

(1) $F$ is of finite area at $p$ (that is $F$ has a
{\rm cusp} at $p$);

(2) the holonomy of a circle in $S$ surrounding $p$ is of {\rm parabolic
type};
\medskip

\noindent $p$ belongs to $V_\Hh$ iff the following equivalent
properties are satisfied:
\smallskip

(i) $p$ corresponds to a boundary component of the convex core
$\overline \Kk(F)$;

(ii)  the holonomy of a circle in $S$ surrounding $p$ is of 
{\rm hyperbolic type}.
\end{prop}

The partition $V=V_\Pp \cup V_\Hh$, so that $r=r_\Pp + r_\Hh$, is said
the {\it type $\theta (F)$} of $F$.  More generally, for every $F\in
\widetilde \Tt(S)$, set $\theta (F)= \theta (\hat F)$.  Any fixed type
$\theta$ determines the subspace $\widetilde \Tt^\theta (S)$ of
hyperbolic structures that share that type; varying $\theta$ we get
the {\it partition by types} of $\widetilde \Tt(S)$.
\medskip

{\bf The fibers of the convex-core map.} 

We want to describe the fibers of the convex-core map
$$\Kk: \widetilde \Tt (S) \to \Tt_\cG (S) \ .$$
\smallskip

Let $h \in {\rm Isom}^+(\mh^2)$ be of hyperbolic type. Denote by
$\gamma=\gamma_h$ its invariant geodesic. Let $P$ be the closed
hyperbolic half-plane determined by $\gamma$ such that the orientation
of $\gamma$ as boundary of $P$ is {\it opposite} to the sense of the
translation $h_{|\gamma}$. 

\begin{defi}\label{fring-end} 
{\rm An {\it $h$-crown} is of the form
$$\Ee = H/h$$ where $H$ is the convex hull in $P$ of a $h$-invariant
closed subset, say $\ \Ee_\infty \subset S^1_\infty$, contained in the
frontier at infinity of $P$.
}
\end{defi}

An $h$-crown $\Ee$ is complete and has geodesic boundary made by
the union of the closed geodesic $\gamma/h$ and complete open geodesics.
$\Ee\setminus \partial \Ee$ is  homeomorphic to $S^1\times
(0,+\infty)$.

Now, let $F\in \widetilde \Tt^\theta (S)$ and
$ \overline\Kk(\hat F) \subset F^\Cc \subset \hat F$ be as above.
Then $F^\Cc$ is obtained by gluing a (possibly empty) crown at
each boundary component $C$ of $\overline\Kk(\hat F)$, associated to
some point $p\in V_\Hh$. This is possible iff, for every $C$ we take
an $h$-crown $\Ee$ such that $h$ is in the same conjugacy class of the
$\hat F$-holonomy of the loop $C$, endowed with the boundary
orientation of $\overline\Kk(\hat F)$ (in other words, ${\rm
length}(\gamma/h)= l(C)$ and both orientations of $\overline\Kk(\hat
F)$ and $\Ee$ are induced by the one of $\hat F$).
\begin{lem}\label{fin-area-fend} $F$ is of finite area iff
all crowns are. A crown $\Ee$ is of finite area iff one of the
following equivalent conditions are satisfied:

(1) $\Ee_\infty/h$ is a finite set. 

(2) $\Ee$ has finitely many boundary components. For every boundary
component $l$, the distance between each end of $l$ and $\partial \Ee
\setminus l$ is $0$.

Every $h$-crown $\Ee$ (every
$F\in \widetilde \Tt(S)$) is the union of exhaustive sequences of
increasing sub-crowns $\Ee_n\subset \Ee$ (sub-surfaces $F_n\subset F$)
of finite area such that $\Ee_{n,\infty}\subset \Ee_\infty$.
\end{lem}
In fact if $\Ee_\infty$ is finite, then the area of $\Ee$ can be
bounded by the sum of the area of a finite set of ideal triangles.  If
$\Ee_\infty$ is not a finite set, then $\Ee$ contains an infinite
family of disjoint ideal triangles.
\medskip

Finally, for every $F\in \Tt_\cG (S)$, the fiber $\Kk^{-1}(F)$
can be identified with the set of all possible patterns of $r_\Hh$
gluable crowns.
\medskip

{\bf Parameters for $\Tt_\cG(S)$.}  The fibers of the convex-core map
are in any sense ``infinite dimensional''. On the other hand, the base
space $\Tt_\cG(S)$ is tame and admits nice parameter spaces, that we
are going to recall.
\medskip

\begin{figure}[ht]
\begin{center}
\includegraphics[width= 4cm]{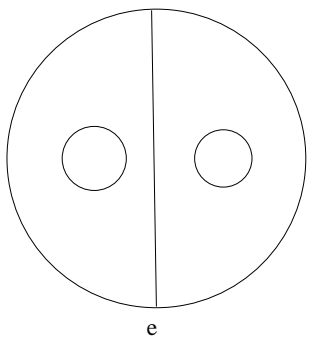}
\caption{\label{pant} A pant and an arc $\rho$.}
\end{center}
\end{figure}

{\bf Length/twist parameters.} This is based on a
fixed {\it pant decomposition} $\Dd$ of $\overline \Sigma$. 
$\Dd$ contains $2g+r-2$ pants obtained by cut/opening $\overline
\Sigma$ at $3g-3+r$ (ordered) disjoint essential simple closed curves
$z_1,\cdots,z_{3g-3+r}$ in $\Sigma$, not isotopic to any boundary
component. Everyone of the $r$ boundary components $C_1,\cdots, C_r$
of $\overline \Sigma$ is in the boundary of some pant. For every
boundary component of a pant $P_k$, corresponding to some $z_j$, we
fix also the unique ``essential'' arc $\rho$ in $P_k$ (shown in
Fig. \ref{pant}) that has the end-points on that component, and we select
furthermore one among these end-points, say $e$.

Denote by $$\mr_+ = \{l\in \mr; \ l>0\}, \ \ \overline \mr_+ = \{l\in
\mr; \ l\geq0\} \ .$$ Consider first the simplest case of $S$ having
$(g,r)=(0,3)$. In this case, set $$\Tt_\cG(S)=\Tt_\cG(0,3) \ .$$  We have
just one pant.  Let us vary the types. If $r_\Hh = 3$, every
hyperbolic structure is determined by the $3$ lengths $(l_1,l_2,l_3)$
of the geodesic boundary components.  If $r_\Hh =2$, by the
corresponding $2$ lengths, and it is natural to associate the value
$0$ to the boundary component that corresponds to the cusp, and so
on. Eventually the octant
$$\overline \mr_+^3=\{(l_1,l_2,l_3); \ l_j\geq 0\}$$ is a natural
parameter space for the whole $\Tt_\cG(0,3)$. The canonical
stratification by open cells of this closed octant corresponds to the
partition by types.

In the general case, let $F\in \Tt_\cG(S)$; then every pant of the
topological decomposition $\Dd$ is associated to a suitable hyperbolic
pant $P_i=P_i(F)$ belonging to $\Tt_\cG(0,3)$. Pant geodesic boundary
components corresponding to some curve $z_j$ have the same length, so
that $F^\Cc$ is obtained by isometrically gluing the hyperbolic pants
at the curves $z_j$. Summing up, $F$ is of the form $$F=F(l,t)$$
$$(l,t)=(l_{C_1},\cdots,l_{C_r},l_{z_1},\cdots
l_{z_{3g-3+r}},t_{z_1},\cdots , t_{z_{3g-3+r}})$$ where $l_{C_i}$
($l_{z_j}$) is the length of the geodesic boundary component (the
simple closed geodesic) of $F^\Cc$ corresponding to $C_i$ ($z_j$). The
{\it twist} parameter $t_{z_j}\in \mr$ specifies the isometric gluing
at $z_j$ as follows.  For every hyperbolic pant, an arc $\rho$ is
uniquely realized by a geodesic arc orthogonal to the boundary.  Then
$F(l,0)$ is the unique hyperbolic structure such that the selected
end-points $e$ of such geometric $\rho$-arcs match by gluing. A
generic $F(l,t)$ is obtained from $F(l,0)$ by modifying the gluing as
follows: if $t_{z_j}>0$, the two sides at any geodesic line $\widetilde
z_j$ in $\mh^2$ over the closed geodesic $z_j$ of $F(l,0)$ translate
by $t_{z_j}$ along $\widetilde z_j$ on the {\it left} to each other.
If $t_{z_j}<0$, they translate on the {\it right} by $|t_{z_j}|$ (``left''
and ``right'' are well defined and only depend on the orientation of
$S$).

We eventually realize in this way that
$$\overline \mr_+^r\times \mr_+^{3g-3+r}\times \mr^{3g-3+r}$$ is a
parameter space (depending on the choice of $\Dd$) for the whole
$\Tt_\cG(S)$. The product by $\mr_+^{3g-3+r}\times \mr^{3g-3+r}$ of
the natural stratification by open cells of $\overline \mr_+^r$,
corresponds to the partition by types. Every cell has
dimension
$$ 6g-6+2r+r_\Hh $$ according to the type. The top-dimensional cell
($r_\Hh = r$) corresponds to the hyperbolic surfaces $F$ without
cusps. $\Tt_{g,r}$ is the lowest dimensional one. Cells that share the
same $r_\Hh$ are isomorphic as well as the corresponding $\Tt^\theta
_\Cc(S)$.  By varying $\Dd$ we actually get an atlas for $\Tt_\cG(S)$
that gives it a {\it real analytic manifold with corner} structure.
\medskip

{\bf Marked length spectrum.}  Length and twist parameters are of
somewhat different nature; in fact we can deal with {\it length
parameters only}. For every $j$, consider: the ``double pant'' obtained
by gluing the two pants of $\Dd$ at $z_j$; the simple closed
curve $z_j'$ obtained by gluing the respective two $\rho$ arcs, and 
$z''_j$ the curve obtained from $z'_j$ via a Dehn twist along
$z_j$. Thus we have further $6g-6 +2r$ simple closed curves on $S$,
and for every $F$ we take the length of the corresponding simple
closed geodesics. In this way we get an {\it embedding}
$$ \Tt_\cG(S) \subset \overline \mr_+^r\times \mr_+^{9g-9+3r} \ .$$
This is the projection onto this finite set of factors of the 
{\it marked length spectrum} injection
$$ {\rm L}: \Tt_\cG(S) \to \overline \mr_+^r \times \mr_+^{\SG'}$$
where $\SG'$ denotes the set of isotopy classes of essential simple
closed curves in $S$, not isotopic to any boundary component.

For more details about the  length/twist parameters and the length spectrum
see for instance \cite{F-L-P, Be-Pe}.
\medskip

{\bf Shear parameters.} This is based on a fixed {\it topological
ideal triangulation} $T$ of $(\hat S,V)$, and works only if $V\neq
\emptyset$. By definition $T$ is a (possibly singular - multi and self
adjacency of triangles are allowed) triangulation of $\hat S$ such
that $V$ coincides with the set of vertices of $T$. There are
$6g-6+3r$ edges $E_1,\cdots, E_{6g-6+3r}$. The idea is to consider
every triangle of $T$ as a hyperbolic ideal triangle and realize
hyperbolic structures $F$ on $S$ by isometrically gluing them at the
geodesic edges, according to the pattern of edge-identifications given
by $T$. By the way, $T$ will be converted in a {\it geometric} ideal
triangulation $T_F$ of $F$. Let us decorate every edge $E$ of $T$ by a
real number $s_E$ and get
$$s=(s(E_1),\cdots,s(E_{6g-6+3r}))\in \mr^{6g-6+3r} \ .$$ These {\it
shear} parameters encode the isometric gluing at each $E_j$, and are
of the same nature of the above twist parameters. Every edge of an
ideal triangle has a distinguished point, say $e$, that is the
intersection of the edge with the unique geodesic line emanating from
the opposite ideal vertex and which is orthogonal to it.  Then set
$F=F(0)$ to be the unique hyperbolic structure such that the
distinguished points match by gluing. A generic $F=F(s)$ is obtained
from $F(0)$ by modifying the gluing according to the left/right moving
rule as before. It turns out that all so obtained hyperbolic
structures $F$ belong to $\Tt_\cG(S)$, and all elements of
$\Tt_\cG(S)$ arise in this way. For every $s$ and every $p_i\in V$,
set
$$s(p_i)=\sum_{E_j\in {\rm Star}(p_i)}s(E_j) \ .$$ We realize that
$$l_{C_i}(F(s))=|s(p_i)|$$ so that, in particular, $p_i\in V_\Pp$ iff
$s(p_i)=0$ and this determines the type $\theta = \theta(F(s))$. This
also shows that the map
$$\Ss:\mr^{6g-6+3r} \to \Tt_\cG(S), \ \ \ F=F(s)$$
is {\it not} injective. For every $p_i\in V_\Hh$,  define the sign
$\epsilon_s(p_i)$ by 
$$ |s(p_i)|= \epsilon_s(p_i)s(p_i) \ .$$ Then, the generic fiber
$\Ss^{-1}(F)$ consists of $2^{r_\Hh}$ points, that is $\Ss$ realizes
all the possible {\it signature} $V_\Hh \to \{\pm 1\}$.  For
the geometric meaning of these signs, see below. For more details about
shear parameters see for instance \cite{Bon}(4).

\medskip

{\bf The enhanced $\Tt_\cG(S)^\#$}. Let us reflect a length/twist
parameter space
$$\overline \mr_+^r\times \mr_+^{3g-3+r}\times \mr^{3g-3+r}$$ of
$\Tt_\cG(S)$ along its boundary components to get
$$ \mr^r\times \mr_+^{3g-3+r}\times \mr^{3g-3+r} \ .$$

This can be considered as a parameter space of the {\it enhanced
Teichm\"uller space} $\Tt_\cG(S)^\#$, obtained by decorating each $F$
with a signature $$\epsilon: V_\Hh \to \{\pm 1\} \ .$$ Moreover, we
stipulate that the sign $\epsilon_i$ associated to $i$ has the meaning
of selecting an orientation of the corresponding $C_i$, by the rule:
{\it $\epsilon_i= +1$ ~iff $C_i$ is equipped with the boundary
orientation.}
\smallskip

To make the notation simpler, it is convenient to extend the signature $\eps$
to the whole of $V$ by stating that $\eps_i = 1$ on $V_\Pp$.  In this way an
enhanced surface can be written as $(F,\eps_1,\ldots,\eps_r)$ with
$\eps_i\in\{\pm 1\}$ and $\eps_i=1$ for $i$ corresponding to a cusp of $F$.

In the same way one can show that the shearing parameters are global
coordinates on $\Tt_\cG(S)^\#$, namely the map
\[
   \Ss^\#:\mr^{6g-6+3r}\rightarrow\Tt_C^\#(S)
\]
defined by $\Ss^\#(s)=(F(s), \sign(s(p_1)),\ldots, \sign(s(p_n)))$, is a
homeomorphism (see~\cite{F-G} for details).

\smallskip

There is a natural {\it forgetting projection}
$$ \phi^\#: \Tt_\cG(S)^\# \to \Tt_\cG(S) \ .$$ We can also define in a
coherent way the {\it enhanced length spectrum}
$${\rm L}^\#: \Tt_\cG(S)^\# \to \mr^r\times \mr_+^{\SG'}$$
by setting
$$l^\#_{C_i}(F,\epsilon)=\eps_i l_{C_i}(F)$$ on the peripheral loops,
and $l^\#_{\gamma}(F,\epsilon)= l_{\gamma}(F)$ elsewhere.  This is an
injection of $\Tt_\cG(S)^\#$, and already the projection onto the
usual finite set of factors as above is an embedding.

\begin{remark}\emph{
    For each $C_i$, the enhanced length is a continuous function of
    $\Tt_\cG^\#(S)$.  On the other hand notice that $\eps_i$ coincides with
    the $\sign$ of $l^\#_{C_i}$, with the rule that the sign of $0$ is $1$.
  }\end{remark}

\subsection{The space of measured geodesic laminations}
\label{lamination}

\begin{defi}\label{lam}{\rm A {\it simple} (complete) geodesic
in $F\in \tilde \Tt(S)$ is a geodesic which admits an arc length
parametrization defined on the whole real line $\R$ that either is
injective (and we call its image a {\it geodesic line} of $F$), or its
image is a simple closed geodesic.  A {\it geodesic lamination} $\Ll$
on $F$ consists of:
\smallskip

(1) A {\it closed} subset $L$ of $F$ (the {\it support}) ;
\smallskip

(2) A partition of $L$ by simple geodesics (the \emph{leaves}).
\smallskip

The leaves together with the connected components of $F\setminus L$
make a {\it stratification} of $S$.}
\end{defi}

\begin{defi}\label{measure}{\rm
Given a geodesic lamination $\Ll$ on $F\in \tilde \Tt(S)$, a
rectifiable arc $k$ in $F$ is {\it transverse} to the lamination if
for every point $p\in k$ there exists a neighbourhood $k'$ of $p$ in
$k$ that intersects each leaf in at most a point and each $2$-stratum
in a connected set.  A {\it transverse measure $\mu$} on $\Ll$ is the
assignment of a positive measure $\mu_k$ on each rectifiable arc $k$
transverse to $\Ll$ (this means that $\mu_k$ assigns a non-negative
{\it mass} $\mu_k(A)$ to every Borel subset of the arc, in a countably
additive way) in such a way that:

(1) The support of $\mu_k$ is $k\cap L$;
\smallskip

(2)  If $k' \subset k$, then $\mu_{k'} = \mu_k|_{k'}$;
\smallskip

(3) If $k$ and $k'$ are homotopic through a family of arcs transverse to
$\Ll$, then the homotopy sends the measure $\mu_k$ to $\mu_{k'}$;
\smallskip
}

Notice that we allow an arc $k$ hitting the boundary of $F^\Cc$ to have
infinite mass, that is $\mu_k(k)=+\infty$.

\end{defi}
\begin{defi}\label{MLS}{\rm
A {\it measured geodesic lamination on $F$} is a pair
$\lambda=(\Ll,\mu)$, where $\Ll$ is a geodesic lamination and $\mu$ is
a transverse measure on $\Ll$. For every $F\in \widetilde \Tt(S)$,
denote by $\Mm\Ll(F)$ the set of measured geodesic laminations on
$F$. Finally, let us define $\Mm\Ll(S)$ to be the set of couples
$(F,\lambda)$, such that $F\in \widetilde \Tt(S)$, and $\lambda \in
\Mm\Ll(F)$. We have the natural projection
$$\pG : \Mm\Ll(S)\to \widetilde \Tt(S) \ .$$ 
}
\end{defi}

\begin{defi}\label{W-S-parts}{\rm
Given $(F,\lambda)\in \Mm\Ll(S)$, the {\it simplicial part} $\Ll_S$ of
$\Ll$ consists of the union of the isolated leaves of $\Ll$. Hence
$\Ll_S$ does not depend on the measure $\mu$.  A leaf, $l$, is called
\emph{weighted} if there exists a transverse arc $k$ such that $k\cap
l$ is an atom of $\mu_k$.  The {\it weighted part} of $\lambda$ is the
union of all weighted leaves.  It depends on the measure and it is
denoted by $\Ll_W=\Ll_W(\mu)$.  }
\end{defi}

\begin{remark}\label{more rem}{\rm The word ``simplicial''
mostly refers to the ``dual'' geometry of the initial singularity of
the spacetimes that we will associate to every $(F,\lambda)$, see
Section~\ref{WR}.  
\smallskip

By property (3) of the definition of transverse measure, if $l$ is
weighted then for every transverse arc $k$ the intersection of $k$
with $l$ consists of atoms of $\mu_k$ whose masses are equal to a
positive number $A$ independent of $k$. We call this number the weight
of $l$.  Since every compact set $K\subset F$ intersects finitely many
weighted leaves with weight bigger than $1/n$, it follows that $\Ll_W$
is a countable set.  As $L$ is the support of $\mu$, then we have the
inclusion $\Ll_S\subset\Ll_W(\mu)$.
}
\end{remark}

\begin{remark}\label{equi-def}{\rm 
There is a slightly different but equivalent definition of $\Mm\Ll(S)$
that runs as follows. We can consider measured geodesic laminations
$\lambda = (\Ll,\mu)$ of $F^\Cc$ requiring furthermore that:
\smallskip

(1) The boundary components of $F^\Cc$ are leaves of $\Ll$;

\smallskip

(2) Every arc $k$ hitting the boundary of $F^\Cc$ {\it
  necessarily} has infinite mass ($\mu_k(k)=+\infty$).
\smallskip

If a boundary component of $F^\Cc$ is isolated in $\Ll$ we stipulate
that it has {\it weight $+\infty$}.  Notice that while a geodesic
lamination on $F^\Cc$ can be regarded also as a particular lamination
on the associated complete surface $\hat F$, condition (2) ensures
that such a {\it measured} lamination cannot be extended beyond
$F^\Cc$. On the other hand, a lamination on $F$ in not in general a
lamination on $\hat F$.

Given any $\lambda$ of $F$ we get a corresponding $\hat \lambda$ of
$F^\Cc$ by adding the (possibly $+\infty$-weighted) boundary
components to the lamination and keeping the same measure. Given $\hat
\lambda$ in $F^\Cc$ we get $\lambda$ in $F$ by just forgetting the
boundary leaves. In particular the empty lamination on $F$ corresponds
to the lamination on $F^\Cc$ reduced to its boundary
components. Clearly this establishes a canonical bijection, hence an
equivalent definition of $\Mm\Ll(S)$. This second definition would
sound at present somewhat unmotivated, so in this section we prefer to
deal with $F$ instead of $F^\Cc$. However, we will see in Section
\ref{WR} that it is the suitable one when dealing with the Lorentzian
``materializations'' of $\Mm\Ll(S)$.}
\end{remark}

{\bf Marked measure spectrum.}\label{mark-spect} Similarly to the
above length spectrum ${\rm L}$, for every $F\in \widetilde \Tt(S)$,
it is defined the {\it marked measure spectrum}
$${\rm I}: \Mm\Ll(F)\to \overline \mr_+^r\times \overline
\mr_+^{\SG'}$$ where for every $\lambda \in \Mm\Ll(F)$ and for every
isotopy class $\gamma$ of essential simple closed curves on $S$, ${\rm
I}_\gamma(\lambda)$ is the minimum of the {\it total variation}
$\mu(c)$ of the ``$\lambda$-transverse component'' of $c$, $c$ varying
among the representatives of $s$. The first $r$ factors correspond as
usual to the curves parallel to the boundary components.
\medskip

{\bf Ray structure.} Every $\lambda = (\Ll,\mu)\in \Mm\Ll(F)$
determines the ray
$$ R_\lambda= \{t\lambda = (\Ll, t\mu); \ t\in [0,+\infty)\}\subset
\Mm\Ll(F)$$ where we stipulate that for $t=0$ we take the empty lamination
 of $F$.  If ${\rm I}_\lambda \neq 0$, then
${\rm I}(R_\lambda)= R_{{\rm I}_\lambda}$, that is the corresponding
ray in $\overline \mr_+^\SG$.

\subsection {The sub-space $\Mm\Ll_\cG(S)$}\label{distinguish}
$$\Mm\Ll_\cG(S)=\{(F,\lambda)\in \Mm\Ll(S);\ F\in \Tt_\cG(S)\}$$
$$ \pG_\cG:\Mm\Ll_\cG(S)\to \Tt_\cG(S)$$ being the natural restriction
of $\pG$ with fibers $\Mm\Ll_\cG(F)$.  

For any $F\in\Tt_\cG(S)$ denotes by $\Mm\Ll_\cG(F)^0$ the set of
laminations on $F$ that do not enter any cusp (namely the closure in
$F^\Cc$ of the lamination support is compact). For a fixed type
$\theta$ we denotes by
$$
  \Mm\Ll_\cG(S)^\theta =\{(F,\lambda)| F\in\Tt^\theta_\cG(S)\,,\
  \lambda\in\Mm\Ll_\cG(S)^0\}
$$ 
and we still denote by $\pG_\cG$ the restriction of the projection on every
$\Mm\Ll_\cG(S)^\theta$.

The spectrum $\rm I$ and the ray structure naturally restrict
themselves. In particular, if $\lambda \in \Mm\Ll_\cG(F)^0$, and $s$
surrounds a cusp of $F$, then ${\rm I}_\lambda(s)=0$. On the other
hand, if $s$ is parallel to a boundary component of $F^\Cc$, then
${\rm I}_\lambda(s)=0$ iff the closure in $F^\Cc$ of the lamination
support $L$ does not intersect that boundary component.

The following Proposition summarizes some basic properties
of the fibers of $\pG_\cG$.

\begin{prop}\label{property} Let $\lambda \in \Mm\Ll_\cG(F)$. Then:
\smallskip

(1) $F\setminus L$ has a finite number of connected components, and
each component belongs to some $\widetilde \Tt(S')$, providing that
we drop out the requirement that $S'$ is non-elementary.
\smallskip

(2) $\lambda$ is disjoint union of a finite set of {\it minimal} {\rm
  [with respect to the inclusion]} measured sublaminations {\rm
  [recall that a lamination $\Ll$ is minimal iff every half-leaf is
  dense in $\Ll$]}. Every minimal sublamination either is compact or
  consists of a geodesic line such that each sub half-line either
  enters a cusp or spirals towards a boundary component of $F^\Cc$.

(3)  $\Ll_W=\Ll_S$.
\smallskip

(4) Either any cusp or any boundary component has a neighbourhood
$U$ such that $\Ll \cap U = \Ll_S \cap U$.
\smallskip

(5) For every arc $c$ in $F$ transverse to $\lambda$, $c\cap L$
is union of isolated points and of a finite union of Cantor sets.
\end{prop}   

For a proof when $F\in \Tt_{g,r}$ we refer for instance to the body
and the references of \cite{Bon}(1). The details for the extension to
the whole of $\Mm\Ll_\cG(S)$ are given for instance in \cite{BSK}.

\begin{remark}\label{equi-def2} {\rm 
If the lamination $\hat \lambda$ of $F^\Cc$ corresponds to $\lambda$
of $F$ as in Remark \ref{equi-def}, then a leaf spiraling towards a
boundary component of $F^\Cc$ as in (2) is no longer a minimal
sublamination of $\hat \lambda$.}
\end{remark}
 
\begin{exa}\label{exa-via-twist-shear} {\rm
We refer to the above
length/twist or shear parameters for $\Tt_\cG(S)$. 
\smallskip

(a) Let $F=F(l,t)$. The union of simple closed geodesics of $F$
corresponding the the curves $z_j$ is a geodesic
lamination $\Ll = \Ll_S$ of $F$. By giving 
each $z_j$ an arbitrary real weight $w_j>0$, we get 
$\lambda(w) \in \Mm\Ll_\cG(F(l,t))^0$.
\smallskip

(b) Let $F=F(s)$. The 1-skeleton of the geometric ideal triangulation
$T_F$ (which is made by geodesic lines) makes a geodesic lamination of
$F$. Every geodesic line is a minimal sublamination. By
giving each geodesic lines an arbitrary weight $w_j>0$, we get
$\lambda(w) \in \Mm\Ll_\cG(F(s))$.  For such a $\lambda =
\lambda(w)$
$$ {\rm I}_{C_i}(\lambda) = \sum_{E_j\in {\rm Star}(p_i)}w(E_j) \ .$$
}
\end{exa}

{\bf Lamination signatures.}  Let $\lambda \in \Mm\Ll_\cG(F)$.  Leaves of
$\lambda$ can spiral around a boundary component $C_i$ in two different ways.
On the other hand two leaves that spiral around $C_i$ must spiral in the same
way (otherwise they would meet each other).

This determines a {\it signature}
$$\sigma(\lambda) : V_\Hh \to \{\pm 1 \}$$ such that
$\sigma_i(\lambda)=-1$ if and only if there are leaves of $\lambda$
spiraling around the corresponding geodesic boundary $C_i$ with a
negative sense with respect to the boundary orientation. In other
words, $\sigma_i(\lambda)$ is possibly equal to $-1$ only if $p_i \in
V_\Hh$ and ${\rm I}_{C_i}(\lambda)\neq 0$, $\sigma_i(\lambda)=1$
otherwise. The signature depends indeed only on the lamination $\Ll$,
not on the measure.

\begin{remark} {\rm If $\lambda = \lambda(w)$ as in Example
\ref{exa-via-twist-shear}(b), then $\sigma_\lambda$ recovers the signs
$\epsilon_s(p_i)$ already defined at the end of Section
\ref{convex-core}.}
\end{remark}

\subsection{ Enhanced bundle $\Mm\Ll_\cG(S)^\#$ and measure spectrum}  
Here we address the question to which extent the (restricted) marked
measure spectrum determines $\Mm\Ll_\cG(S)$. For example, this is
known to be the case if we restrict to $\Mm\Ll_{g,r}^0$ {\it i.e.} to
laminations over $\Tt_{g,r}$ that do not enter the cusps (see for
instance \cite{Bon}(1)).  We want to extend this known result.
\smallskip

We have seen in Proposition \ref{property} that a measured geodesic
lamination $\lambda$ on $F\in\Tt_\cG(S)$ is the disjoint union of a
compact part, say $\lambda_c$ (that is far away from the geodesic
boundary of $F^\Cc$ and does not enter any cusps), with a part, say
$\lambda_b$, made by a finite set of weighted geodesic lines
$l_1,\ldots,l_n$ whose ends leave every compact subset of $F$. Notice
that $\sigma(\lambda)=\sigma(\lambda_b)$.

Let us take such a geodesic line $l$ on $F\in\Tt_\cG(S)$.  We can select a
compact closed interval $J$ in $l$ such that both components of $l
\setminus J$ definitely stay either within a small $\eps$-neigbourhood
of some boundary component of $F^\Cc$, or within some cusp.  $J$ can
be completed to a simple arc $c$ in $\hat S$ with end-points in $V$,
just by going straightly from each end-point of $J$ to the
corresponding puncture. It is easy to see that the homotopy class with
fixed end-points of the so obtained arc $c$ does not depend on the
choice of $J$. For simplicity we refer to it as the ``homotopy class''
of $l$. We can also give the end-points of $c$ a sign $\pm 1$ in the
very same way we have defined the signature of a lamination on $F$
(recall that the sign is always equal to $1$ at cusps). We can prove
\begin{lem} Given any $F\in  \Tt_\cG(S)$, every homotopy class $\alpha$ 
of simple arcs on $\hat S$ with end-points on $V$, and every signature
of the end-points (compatible with the type of $F$) can be realized by
a unique geodesic line $l$ of $F$ whose ends leave every compact set
of $F$.  Moreover, the members of a finite family of such geodesic
lines are pairwise disjoint iff the signs agree on every common
end-point and there are disjoint representatives with end-points on
$V$ of the respective homotopy classes. Analogously they do not
intersect a compact lamination $\lambda_c$ iff so do suitable
representatives.
\end{lem}
By using the lemma, we can prove (see ~\cite{BSK})
\begin{prop}\label{homotopy-arc} Let $\lambda \in \Mm\Ll_\cG(F)$. Then 
the support of $\lambda_b$ is determined by the homotopy classes of
its geodesic lines $l_i$ and the signature of $\lambda$. More
precisely, given any $\lambda_c$, every finite set of homotopy classes
of simple weighted arcs on $\hat S$, with signed end-points in $V$
(providing the signature being compatible with the type of $F$), admitting
representatives that are pairwise disjoint and do not intersect
$\lambda_c$, is uniquely realized by a lamination  $\lambda_b$ such that
$\lambda = \lambda_b \cup \lambda_c \in \Mm\Ll_\cG(F)$. 
\end{prop}

\begin{prop}\label{map-iota}
Let $F, F' \in \Tt_\cG(S)$. Assume that $F$ is without cusps (that is
$F$ belongs to the top dimensional cell of $ \Tt_\cG(S)$). Then there
is a natural map
\[
   \iota:\Mm\Ll_\cG(F)\rightarrow\Mm\Ll_\cG(F')
\]
such that for every (isotopy class of) simple closed curve $\gamma$ on
$S$, we have
\[
   {\rm I}_\gamma(\lambda)={\rm I}_{\gamma}(\iota(\lambda))\,.
\]
\end{prop}
\Dim Assume first that $\lambda = \lambda_c \in \Mm\Ll_\cG(F)$. Then
there is a unique $\lambda' = \lambda'_c \in \Mm\Ll_\cG(F')$ with the
same spectrum.  For we can embed $F'$ in the double surfaces of
$(F')^\Cc$, say $DF'$ which is complete and of finite area. The
measure spectrum of $\lambda_c$ induces a measure spectrum of a unique
lamination $\lambda''_c$ on $DF'$ (by applying the result on the
spectrum in the special case recalled at the beginning of this
Section). Finally we realize that the compact support of $\lambda''_c$
is contained in $F'$ giving us the required $\lambda'_c$.
So the map $\iota$ can be defined for laminations with compact support.

Given a general lamination $\lambda=\lambda_c\cup\lambda_b$, we can
define $\lambda'_c$ as before, while $\lambda'_b$ is the unique
lamination of $F'$ (accordingly with Proposition \ref{homotopy-arc})
that share with $\lambda_b$ the same homotopy classes, weights and
signs at $V_\Hh(F)\cap V_\Hh(F')$. Notice that $\lambda'_b$ is
disjoint from $\lambda'_c$: in fact one can construct an isotopy of
$S$ sending the supports of $\lambda_b$ and $\lambda_c$ to the
supports of $\lambda'_b)$ and $\lambda'_c$. Finally set
$\iota(\lambda)=\iota(\lambda_b)\cup\iota(\lambda_c)$.  \cvd
\begin{cor}\label{orbit} 
If both $F$ and $F'$ are without cusps, then the map $\iota$ is
bijective. More generally, for every $\lambda'\in\Mm\Ll_\cG(F')$,
$\iota^{-1}(\lambda')$ consists of $2^k$ points, where $k$ is the
number of cusps of $F'$ entered by $\lambda'$.
\end{cor}
 In fact, for every $F\in\Tt_\cG(S)$ (not necessarily in the top
 dimensional cell), there is a natural action of $(\mz/2\mz)^r$ on
 $\Mm\Ll_\cG(F)$ determined as follows. Let $\rho_i = (0,\dots, 1,
 \dots,0)$, $i=1,\dots r$, be the $i$th element of the standard basis
 of $(\mz/2\mz)^r$.  Let $\lambda \in \Mm\Ll_\cG(F)$. First define the new
signature $\rho_i \sigma(\lambda)$ by setting: 
\smallskip

$\rho_i\sigma(\lambda)(p_j)= \sigma(\lambda)(p_j)$ if $i\neq j$; 
\smallskip

$\rho_i\sigma(\lambda)(p_i)= \sigma(\lambda))(p_i)$ if
either $p_i\in V_\Pp (F)$ or $p_i \in V_\Hh(F)$ and ${\rm I}_{C_i}(\lambda)=0$;
\smallskip

$\rho_i\sigma(\lambda)(p_i)= -\sigma(\lambda)(p_i)$, otherwise.

This naturally extends to every $\rho \in (\mz/2\mz)^r$, giving the
signature $\rho\sigma(\lambda)$. Finally set $\rho(\lambda) =
\rho(\lambda_b) \cup \lambda_c$ where (accordingly again with
Proposition \ref{homotopy-arc}) $\rho(\lambda_b)$ is the unique
lamination that shares with $\lambda_b$ the homotopy classes and the
weights, while its signature is $\rho\sigma(\lambda)$. Clearly the
orbit of $\lambda$ consists of $2^k$ points, where $k$ is the number
of $p_i$ in $V_\Hh(F)$ such that ${\rm I}_{C_i}(\lambda)\neq
0$. Finally $\iota^{-1}(\lambda')$ in Corollary \ref{orbit} is just an
orbit of such an action.  We call the action on $\Mm\Ll_\cG(F)$ of the
generator $\rho_i$, the {\it reflection along $C_i$} (even if it
could be somewhat misleading, as in some case it is just the identity).
\smallskip

If we restrict over the top-dimensional cell of $\Tt_\cG(S)$,
$\pG_\cG$ is a bundle and we can use the first statement of the
Corollary in order to fix a trivialization. The same fact holds for
every restriction $\pG_\cG : \Mm\Ll_\cG(S)^\theta \to
\Tt_\cG(S)^\theta$, type by type. On the other hand, because of the
last statement of the Corollary, this is no longer true for the whole
$\pG_\cG$.  In order to overcome such phenomenon, one can introduce
the notion of {\it enhanced lamination}.  An enhanced lamination on
$F\in\Tt_\cG(S)$, is a couple $(\lambda, \eta)$ where $\lambda \in
\Mm\Ll_\cG(F)$, and $\eta: V \to \{\pm\}$ is a {\it relaxed signature}
such that:
\smallskip

$\eta_i=\sigma_i(\lambda)$ if either $p_i\in V_\Hh(F)$ or $p_i\in V_\Pp(F)$
and ${\rm I}_{C_i}(\lambda) = 0$;
\smallskip

$\eta_i$ is arbitrary otherwise. 
\smallskip

Notice that there are exactly $2^k$ $(\lambda,\eta)$ enhancing a given
$\lambda \in \Mm\Ll_\cG(F)$, where $k$ is the number of cusps entered
by $\lambda$. Clearly the above action of $(\mz/2\mz)^r$ extends on
enhanced laminations: $\rho(\lambda,\eta)= (\rho(\lambda), \rho(\eta))$,
where $\rho(\eta)$ is uniquely determined by the above requirements and by
the fact that $\rho\sigma(\lambda)$ possibly modifies $\sigma(\lambda)$ only
on $V_\Hh$. In particular this holds for the generating reflections $\rho_i$.

We denote by $\Mm\Ll^\#_\cG(F)$ the set of such $(\lambda,\eta))$ on $F$.
Finally we can define the {\it enhanced measure spectrum}
$$ {\rm I}^\#: \Mm\Ll^\#_\cG(F) \to  \mr^r \times\mr_+^{\SG'}$$
such that:
\smallskip

$${\rm I}^\#_{\gamma}(\lambda,\eta)={\rm I}_\gamma(\lambda)$$ for every
$\gamma\in\SG'$;
\smallskip

$$ {\rm I}^\#_{C_i}(\lambda,\eta)=\eta_i {\rm I}_{C_i}(\lambda)$$
for every peripheral loop $C_i$.

Here is the enhanced version of Proposition \ref{map-iota}.
\begin{cor}\label{map-iota-en}
Let $F, F' \in \Tt_\cG(S)$. Then there is a natural {\rm bijection}
\[
   \iota^\#:\Mm\Ll_\cG(F)^\# \rightarrow\Mm\Ll_\cG(F')^\#
\]
such that for every (isotopy class of) simple closed curve $\gamma$ on
$S$, we have
\[
   {\rm I}^\#_\gamma((\lambda,\eta))={\rm
   I}^\#_{\gamma}(\iota^\#(\lambda,\eta))\,.
\]
\end{cor}

\begin{prop}\label{m-s}
(i) The enhanced spectrum ${\rm I}^\#$ realizes an embedding of every
$\Mm\Ll_\cG(F)^\#$ into $\mr^r\times \overline \mr_+^{\SG'}$. Only the
empty lamination goes to $0$.  The image is homeomorphic to
$\mr^{6g-6+3r}$.  The image of $\Mm\Ll_\cG(F)^{\#,0}$ (that is the set
of enhanced laminations that do not enter any cusp) is homeomorphic to
$\mr^{6g-6+2r+r_{\Hh}}$
\smallskip

(ii) For every pant decomposition $\Dd$ of $\overline \Sigma$,
consider the subset {\rm [already considered to deal with the length
spectrum]}
$$\SG_\Dd =
\{C_1,\cdots,C_r,z_1,z_1',z_1'',\cdots,,z_{3g-3+r},z_{3g-3+r}',z_{3g-3+r}''\}
\subset \SG \ .$$ The projection onto this {\rm finite} set of factors
is already an embedding of $\Mm\Ll_\cG(F)^\#$.  By varying $\Dd$ we
get an atlas of a {\it PL structure} on $\Mm\Ll_\cG(F)^\#$ ({\it i.e.}
on $\mr^{6g-6+3r}$). Similar facts hold for the restriction to
$\Mm\Ll_\cG(F)^{\#,0}$.
\smallskip

(iii) Finite laminations are dense in $\Mm\Ll_\cG(F)^\#$
( $\Mm\Ll_\cG(F)^{\#,0}$).
\smallskip

(iv) For every $F,F'\in \Tt_\cG(S)$, there is a canonical
identification between the respective sets of finite enhanced measured
geodesic laminations, and this extends to a canonical PL isomorphism
between $\Mm\Ll_\cG(F)^\#$ and $\Mm\Ll_\cG(F')^\#$, which respects the
ray structures. Similarly for $\Mm\Ll_\cG(\cdot)^{\#,0}$.
 \end{prop}
 
\Dim We will sketch the proof of this proposition. We assume that the
result is known when $S$ is compact (see \cite{Bon}(1), \cite{F-L-P}).
Thanks to Proposition \ref{map-iota-en} it is enough to deal with 
$F$ without cusps. 
Then the double $DF$ of $F^\Cc$ is compact, and we
consider on $DF$ the involution $\tau$ that exchange the two copies of
$F$.  Let us denote by $ML(F)$ the set of $\tau$-invariant measured
geodesic laminations on $DF$ that do not contain any component of
$\partial F^\Cc$.  The idea is to construct a map
 \[
     T: \Mm\Ll_\cG(F)\rightarrow ML(F)
 \]
 that is surjective and such that
 
 (1) the fiber over a lamination $\lambda' \in ML(F)$ consists of $2^k$
laminations of $\Mm\Ll_\cG(F)$, where $k$ is the number of boundary
components of $F^\Cc$ that intersect the support of  $\lambda'$.
 
 (2) For every $\lambda \in \Mm\Ll_\cG(F)$, the restrictions to
 $\SG$ of both the spectrum of $T(\lambda)$ and of $\lambda$ coincide. 
\smallskip

The existence of the map $T$ and the known results in the special cases
recalled above will imply the Proposition.
 
The construction of the map $T$ runs as follows. Let $\lambda =
\lambda_b \cup \lambda_c \in \Mm\Ll_\cG(F)$ be decomposed as above.
We define $T(\lambda_c)$ to be the double of $\lambda_c$ in $DF$.  For
each leaf $l_i$ of $\lambda_b$, take a ``big'' segment $J_i\subset
l_i$, and complete it to a simple arc $l'_i$ properly embedded in
$(F^\Cc,\partial F^\Cc)$, obtained by going straightly from each
end-point of $J_i$ to the corresponding boundary component along an
orthogonal segment. Clearly the double of $l'_i$ is a simple non-trivial 
curve in $DF$, so there is a geodesic representative, say
$c_i$, that is $\tau$-invariant and simple.  Since $l_i\cap
l_j=\varnothing$ the same holds for the $c_i$'s. Moreover, since
$l_i\cap\lambda_c=\varnothing$, the intersection of $c_i$ with
$T(\lambda_c)$ is also empty. So we can define
\[
    T(\lambda)=
T(\lambda_b)\cup (c_1,a_1)\cup (c_2,a_2)\cup\ldots\cup (c_n,a_n)\,.
\] where $a_i$ is the initial weight of $l_i$.
This map satisfies (2) by construction; moreover, it follows from
Corollary \ref{orbit} that (1) holds for every $\lambda'$ belonging to
the image of $T$.  The only point to check is that the map is
surjective.  The key remark is that for every $\lambda' \in ML(F)$,
every leaf $l$ hitting $\partial F^\Cc$ is necessarily closed.  As it
is $\tau$-invariant, then $l$ is orthogonal to $\partial F^\Cc$, and
if $l$ intersects $\partial F^\Cc$ twice, then it is closed.  Suppose
that $l$ is a geodesic line, so that $l$ meets $\partial F$ exactly
once. On the other hand, we know that the closure of $l$ is a minimal
sublamination $\lambda''$, such that every leaf is dense in it. Thus
if $l''\neq l$ is another leaf in $\lambda''$, then it intersects
$\partial F^\Cc$ in a point $p$.  Since $l$ is dense in
$\lambda''$, there is a sequence of points in $l\cap\partial F^\Cc$
converging to $p$ and this contradicts the assumption that $l$
intersects $\partial F^\Cc$ once.

Thus a lamination in $ML(F)$ is given by the double of a compact
lamination $\lambda_c$ in $F$ and of a finite number of weighted
simple geodesics arcs in $F$ hitting orthogonally $\partial F^\Cc$.
These arcs can be completed to give a family of simple arcs on $\hat S$
with end-points on $V$. Fix a signature on the end-points of such arcs.
Finally we can apply Proposition \ref{homotopy-arc} to these data
and we get a suitable $\lambda = \lambda_b \cup \lambda_c \in \Mm\Ll_\cG(F)$
such that $T(\lambda)=\lambda'$.  \cvd
 
Finally we can define the map
$$\pG_\cG^\#: \Mm\Ll_\cG(S)^\# \rightarrow \Tt_\cG^\#(S) \ .$$
The total space is defined  
as the set of pairs 
$$((F,\eps), (\lambda,\eta))$$ such that
\begin{enumerate}
\item $(F,\eps)=(F,\eps_1,\ldots,\eps_r)\in\Tt_\cG(S)^\#$;
 \item $(\lambda,\eta)=(\lambda,\eta_1,\ldots,\eta_r)\in\Mm\Ll_\cG(F)^\#$
\end{enumerate}
Clearly $$ \phi^\#\circ \pG^\# = \pG \circ \phi_{\Mm\Ll}^\# $$ where
$\phi_{\Mm\Ll}^\#$ denotes the {\it forgetting projection} of  
$\Mm\Ll_\cG(S)^\#$ onto $\Mm\Ll_\cG(S)$. We are going to see that
in fact $\pG_\cG$ determines a {\it bundle} of enhanced lamination, that
admits furthermore natural {\it trivializations} $\tG$.
It follows from the previous discussion that the image of  
${\rm I}^\#$ does not depend on the choice of $F$, hence ${\rm I}^\#(S)$
is well defined. We want to define a natural bijection 
$$ \tG: \Tt_\cG^\#(S)\times {\rm I}^\#(S) \rightarrow
\Mm\Ll_\cG(S)^\# \ .$$ 
For every $\xi \in {\rm
I}^\#(S)$ and $F\in\Tt_\cG(S)$ there is a unique 
$(\lambda(\xi), \eta(\xi)) \in
\Mm\Ll_\cG(F)^\#$ that realizes $\xi$.  
So, let us put $$\tG(F,\eps,\xi)=(F,\eps,\rho_\eps(\lambda(xi),\eta(\xi))\ .$$

It follows from the previous discussion that $\tG$ is a bijection. We
stipulate that it is a homeomorphism, determining by the way a
topology on $\Mm\Ll_\cG(S)^\#$. Summing up, the map
$$\pG^\#: \Mm\Ll_\cG(S)^\# \to \Tt_\cG(S)^\# $$ can be considered as a
{\it canonically trivialized} fiber bundle having both the base space
and the fiber (analytically or PL) isomorphic to
$\mr^{6g-6+3r}$. Different choices of the base surface $F_0$ lead to
isomorphic trivializations, via isomorphisms that preserve all the
structures. These trivializations respect the ray
structures.  When $S$ is compact this specializes to the
trivialized bundle $ \Tt_g \times \Mm\Ll_g \to \Tt_g $ mentioned in
the Introduction.
\begin{remark}{\rm 
The definition of $\tG$ could appear a bit distressing at a first
sight.  However the geometric meaning is simple. Given a spectrum
of positive numbers, this determines the lamination up to choosing the
way of spiraling towards the boundary components.  If we give a
sign to the elements of the spectrum corresponding to the boundary
components, this allows to reconstruct the lamination by the rule: if
the sign is positive the lamination spiral in the positive way, if the
sign is negative the lamination spirals in the negative way {\it with
respect to a fixed orientation of the boundary component}. In the
non-enhanced set up, we have stipulated to use the boundary
orientation induced by the surface one.  Since the elements of an
enhanced Teichm\"uller space can be regarded as hyperbolic surfaces
equipped with an (arbitrary) orientation on each boundary component,
it seems natural to reconstruct the lamination from the spetrum ${\rm I}^\#$
by means of such boundary component orientations.

This choice is suitable in view of the  {\it earthquake flow} that we are
going to define on $\Tt_\cG(S)^\#$. }
\end{remark}

\subsection{Grafting, bending, earthquakes}\label{bend-quake}
Let $(F,\lambda)\in \Mm\Ll(S)$. {\it Grafting} $(F,\lambda)$ produces
a deformation $Gr_\lambda(F)$ of $F$ in $\Pp(S)$, the
Teichm\"uller-like space of {\it complex projective structures} ({\it
i.e.}  $(S^2,PSL(2,\mc)$-{\it structures}) on $S$.
\smallskip

$3$-dimensional {\it hyperbolic bending} produces the $H$-{\it hull}
of $Gr_\lambda(F)$, that is, in a sense, its
``holographic image'' in $\mh^3$.
\smallskip

The {\it left (right) earthquake} produces (in particular) a new
element $\beta^L_\lambda(F)$ ($\beta^R_\lambda(F)$) in $\widetilde
\Tt(S)$.
\smallskip

We will see in Section \ref{WR} how these constructions are {\it
materialized} within the {\it canonical Wick rotation-rescaling
theory} for $MGH$ Einstein spacetimes. For example, the grafting is
eventually realized by the {\it level surfaces of the cosmological
times}; earthquakes are strictly related to the {\it Anti de Sitter
bending} procedure.
\smallskip

Here we limit to recall a few details about earthquakes, purely in terms
of hyperbolic geometry.
\medskip

{\bf Features of arbitrary $(F,\lambda)$.}  In such a
general case, the leaves of $\lambda$ possibly enter the crowns
of $F$.  If $F$ is of finite area (see Lemma
\ref{fin-area-fend}), basically the conclusions of Proposition
\ref{property} still hold. The only new fact is that possibly there is
a finite number of isolated geodesic lines of $\lambda$ having at
least one end converging to a point of some $\Ee_\infty$.
\smallskip

The situation is quite different if $F$ is of infinite area.  The
set of isolated geodesic lines of $\lambda$ that are not entirely
contained in one crown $\Ee$ is always finite. On the other hand, (1),
(2), (3) and (5) of Proposition \ref{property} definitely fails. For
example, the support of a lamination $\lambda$ could contain bands
homeomorphic to $[0,1]\times \R$, such that every $\{t\}\times \R$
maps onto a geodesic line of $\lambda$. Both ends of every such a line
converge to some $\Ee_\infty$. We can also construct transverse
measures such that $L_W$ is dense in such bands. This also shows that
in general $\Ll_S$ is strictly contained in $\Ll_W$.

In general the fibers of ${\rm I}$ are, in any sense, infinite
dimensional.  For example we have:
\begin{lem}\label{I=0} ${\rm I}^{-1}(0) \subset \Mm\Ll(F)$ consists
of laminations such that the support is entirely contained in the
union of crowns.
\end{lem}

On the other hand, the image of ${\rm I}$ is tame, in fact:
\begin{prop}\label{I=I} ${\rm I}(\Mm\Ll(F))= {\rm I}(\Mm\Ll_\cG(\Kk(F))$.
\end{prop}

{\bf Earthquakes along finite laminations of $\Mm\Ll_\cG(F)$.} As
finite laminations are dense, and arbitrary laminations $\lambda \in
\Mm\Ll_\cG(F)$ look like finite ones at cusps and boundary components
of $F^\Cc$, it is important (and easy) to understand earthquakes in
the finite case.

\begin{exa}\label{more-exa-via-twist-shear}
{\rm Let us consider again the Examples \ref{exa-via-twist-shear}.
Let $(F(l,t)$ be such that all twist parameters are strictly positive.
Then, {\it by definition} $(F(l,t),\lambda(t))$ is obtained from
$(F(l,0),\lambda(t))$ via a {\it left earthquake (along the
measured geodesic lamination $\lambda(t)$ on $F(l,0)$)}.
$(F(l,-t),\lambda(t))$ is obtained from $(F(l,0),\lambda(t))$
via a {\it right earthquake (along the measured geodesic lamination
$\lambda(t)$ on $F(l,0)$)}.  In the reverse direction,
$(F(l,0),\lambda(t))$ is obtained from $(F(l,t),\lambda(t))$
via a {\it right earthquake}, and so on.  This pattern of earthquakes
does preserve the types.
\smallskip

Similarly, let $F(s)$ be such that all shear parameters are strictly
positive.  Then, by definition $(F(s),\lambda(s))$ is obtained from
$(F(0),\lambda(s))$ via a {\it left earthquake (along the measured
geodesic lamination $\lambda(s)$ on $F(0)$)}.  $(F(-s),\lambda(s))$ is
obtained from $(F(0),\lambda(s))$ via a {\it right earthquake (along
the measured geodesic lamination $\lambda(s)$ on $F(0)$)}.  In the
reverse direction, $(F(0),\lambda(s))$ is obtained from
$(F(s),\lambda(s))$ via a {\it right earthquake}, and so on. This
pattern does not preserve the types, for $F(0)\in \Tt_{g,r}$, while
$F(s)$ is without cusps. Moreover, $\lambda(s)$ has the following
special property:
\smallskip

{\it For every boundary component $C_i$ of $F(s)^\Cc$
$$ l_{C_i}(F(s)) = {\rm I}_{C_i}(\lambda(s)) \ .$$}
}
\end{exa}

For every $(F,\lambda) \in \Mm\Ll_\cG(S)$, $\lambda$ being finite, the
definition of {\it $(F',\lambda')$ obtained from $(F,\lambda)$ via a
left (right) earthquake} extends {\it verbatim} the one of the above
examples, so that $(F',\lambda')\in \Mm\Ll_\cG(S)$, $\lambda'$ is also
a finite lamination, and $(F,\lambda)$ is obtained from
$(F',\lambda')$ via the {\it inverse} right (left) earthquake. 
\smallskip
 
{\bf Quake cocycles and general earthquakes.}  It is convenient to
describe earthquakes by lifting everything to the universal
covering. Let us set as usual
$$\overline \Kk(\hat F) \subset F^\Cc \subset \hat F = \mh^2/\Gamma \
.$$ Then $F^\Cc$ lifts to a $\Gamma$-invariant {\it straight convex
set} $H$ of $\mh^2$ ({\it i.e} $H$ is the closed convex hull of an
ideal subset of $S^1_\infty$), and $\lambda$ lifts to a
$\Gamma$-invariant measured geodesic lamination on $\mathring{H}$, that, for
simplicity, we still denote by $\lambda$. If $F\in \Tt_\cG$, then
$\overline \Kk(\hat F) = F^\Cc$.

\begin{lem}\label{cocycle} 
Let $(F,\lambda)\in \Mm\Ll_\cG(S)$ such that $\lambda$ is finite. Then
there exists a {\rm left-quake cocycle}
\[B^L_\lambda :\mathring H \times \mathring H \rightarrow
  PSL(2,\mr)
\]
such that  
\begin{enumerate}
\item
$B^L_\lambda(x,y)\circ B^L_\lambda(y,z)=B^L_\lambda(x,z)$ for every 
$x,y,z\in \mathring H$.
\item
$B^L_\lambda(x,x)=Id$ for every $x\in\mathring H$.
\item
$B^L_\lambda$ is constant on the strata of the stratification of
$\mathring H$ determined by $\lambda$.
\item
$B_\lambda(\gamma x,\gamma y)= \gamma B_\lambda(x,y)\gamma^{-1}$, for
every $\gamma \in \Gamma$.
\item
For every $x_0$ belonging to a $2$-stratum of $\mathring H$,
\[
   \mathring H\ni x\mapsto B^L_\lambda (x_0,x)x\in\mh^2 
\]
lifts the left earthquake $\beta^L_\lambda(F)$ to $\mathring H$.  This
cocycle is essentially unique.  There exists a similar {\rm
right-quake cocycle} $B^R_\lambda$.
\end{enumerate}
\end{lem}
The proof is easy and the earthquake is equivalently encoded by its
cocycle. For a general $(F,\lambda)$ we look for (essentially unique)
{\it quake-cocycles} that satisfy all the properties of the previous
Lemma, with the exception of the last one, and requiring furthermore
that
\medskip

{(*) \it If $\lambda_n\rightarrow\lambda$ on a $\eps$-neighbourhood of
the segment $[x,y]$ and $x,y \notin L_W$, then
$B_{\lambda_n}(x,y)\rightarrow B_{\lambda}(x,y)$.}
\medskip

Given such cocycles we can use the map of (5) in the
previous Lemma as the {\it general definition of earthquakes}.
\smallskip

For example, if $(F,\lambda)\in \Mm\Ll_\cG(S)$ the cocycle can be
derived by using Lemma \ref{cocycle}, the density of finite
laminations and the fact that we require $(*)$. If $(F',\lambda')$
results from the left earthquake starting at $(F,\lambda)$, then
this last belongs to $\Mm\Ll_\cG(S)$ and $(F,\lambda)$ is obtained
from it via the inverse right earthquake.
\smallskip

In fact, in \cite{Ep-M} Epstein-Marden defined these quake-cocycles in
general (extending the construction via finite approximations).
Strictly speaking they consider only the case of (arbitrary) measured
geodesic laminations on $\mh^2$, but the same arguments holds for
laminations on arbitrary straight convex sets $H$ - see also
\cite{Be-Bo} for more details. Hence general left (right) earthquakes
$$ (F',\lambda') = \beta^L(F,\lambda)$$
so that
$$(F,\lambda) = \beta^R(F',\lambda')$$ are eventually defined for
arbitrary $(F,\lambda)\in \Mm\Ll(S)$.  We will also write $F'=
\beta^L_\lambda(F)$, $\lambda'= \beta^L_\lambda(\lambda)$.
\medskip

{\bf Earthquake flows on $\Mm\Ll_\cG(S)$.}  Let $\lambda \in
\Mm\Ll_\cG(F)$. Consider the ray $(F,t\lambda)$, $t\in
[0,+\infty)$. Then, for every $t>0$, set 
$$(F_t,\lambda_t) = (\beta^L_{t\lambda}(F),
\frac{1}{t}\beta^L_{t\lambda}(t\lambda)), \ \ t\geq 0 \ . $$
This continuously extends at $t=0$ by
$$(F_0,\lambda_0)=(F,\lambda) \ .$$
We have  
$$((F_t)_s,(\lambda_t)_s)=(F_{t+s},\lambda_{t+s})$$
hence this defines the so called {\it
left-quake flow} on $\Mm\Ll_\cG(S)$. In particular this allows to
define a sort of ``exponential''  map
$$\psi^L: \Mm\Ll_\cG(F) \to \Mm\Ll_\cG(S)$$ by evaluating the flow at
$t=1$.  We do similarly for the {\it right-quake flow}. 

Let $p_i\in V$ and $C_i$ be the curve surrounding it; as ${\rm
I}_{C_i}(t\lambda) = t {\rm I}_{C_i}(\lambda)$, there is a unique ``critical
value'' $t_i$ (see below) such that ${\rm I}_{C_i}(t\lambda)=
l_{C_i}(F)$.

For every $t$, we denote by $l(t)$ the marked length spectrum of $F_t$,
by $\theta(t)$ its type, by ${\rm I}(t)$ the marked measure spectrum of
$\lambda_t$, by $\sigma_t: V \to \{\pm\}$ its signature, and so
on.  The following Lemma describes the behaviour of these objects
along the flow.

\begin{lem}\label{ray-quake} 
The marked measure spectrum is constant for every $t$, that is
\[
    {\rm I}_\gamma(t)={\rm I}_\gamma(0)\ \ \textrm{for every }\gamma\in\SG\,.
\]

Let $p_i\in V$ and $C_i$ be
the curve surrounding it.
\smallskip

If $p_i\in V_{\Hh}(0)$, then:
 $$ l_{C_i}(t)= |l_{C_i}(0)-t\sigma_i(\lambda) {\rm I}_{C_i}(0)| $$
 and 
$$ \sigma_i(t)=\sign [l_{C_i}-t\sigma_i(\lambda) 
{\rm I}_{C_i}(0)]\sigma_i(0) \ .$$
 \medskip

 If $p_i\in V_\Pp(0)$ then:
 $$ l_{C_i}(t)=t{\rm I}_{C_i}(0) $$
 and
 $$  \sigma_i(t)=-1$$

\end{lem}

As every $\lambda \in \Mm\Ll_\cG(F)$ looks finite at cusps and
boundary components of $F^\Cc$, it is enough (and fairly easy) to
check the Lemma in the finite case, by using also Examples
\ref{more-exa-via-twist-shear}.

\begin{remark}{\rm
    If $p_i\in V_\Pp$ and the lamination enters the corresponding
    cusp, then for $t>0$ the cusp opens on a geodesic boundary
    component whose length linearly depends on $t$ with slope equal to
    ${\rm I}_{C_i}(0)$. The way of spiraling of $\lambda_t$ around $p_i$ is
    always negative (positive for right earthquakes).
    
    Let us consider more carefully the case $p_i\in V_\Hh$.  Notice
    that if $\lambda$ does not spiral around $C_i$ then the length of
    $C_i$ is constant.  In the other cases let us distinguish two
    possibilities according to the sense of spiraling of $\lambda$.

(1) Case $\sigma_i(0)=-1$. Then for every $t>0$,
$$\sigma_i(t)=-1 \ , \ \ l_{C_i}= l_{C_i}(0)+t{\rm I}_{C_i}(0) \ . $$
Thus the length of $C_i$ increases linearly of slope $ {\rm
I}_{C_i}(0)$ and the laminations continues to spiral in the negative
direction.

\smallskip

(2) Case $\sigma_i(0)=1 \ $.  There is a critical time
$t_i=l_{C_i}(0)/{\rm I}_{C_i}(0)$. Before $t_i$ the length of $C_i$
decreases linearly and the lamination spiral in the positive
direction.  At $t_i$, $C_i$ is become a cusp. After $t_i$, $C_i$ is
again a boundary component but the way of spiraling is now negative.
}\end{remark}

\begin{remark}\label{I=l}{\rm The above Proposition points out
in every $\Mm\Ll_\cG(F)$ the set :
$$ \Vv_\cG(F) = \{ \lambda; \ {\rm I}_{C_i}(\lambda)<  l_{C_i}(F);\  
i \in V_\Hh \} \ .$$

Note that this set is {\it not} preserved by the canonical
bijections stated in Proposition
\ref{m-s}(iv).}
\end{remark}
\begin{cor}\label{preserve-type} The restriction of the 
exponential-like map $\psi^L$ to $\Vv_\cG(F)\cap \Mm\Ll_\cG(F)^0$
preserves the type and the signatures.  The restriction of this map to
the whole of $ \Vv_\cG(F)$ has generic image over the top-dimensional cell
of $\Tt_\cG(S)$. 
\end{cor}

{\bf The quake-flow on $\Mm\Ll_\cG(S)^\#$.}  We will define an
earthquake flow on $\Mm\Ll_\cG(S)^\#$ that will satisfy the following
properties

(1) $\beta^\#_t\circ\beta^\#_s=\beta^\#_{t+s}$.
\smallskip

(2) Every flow line $\{\beta^\#_t(F,\eps,\lambda,\eta)|t>0\}$ is 
{\it horizontal}
with respect to the trivialization of $\Mm\Ll_\cG^\#(S)$.  This means that the
enhanced lamination is constant along the flow.  \smallskip

(3) If we include $\Mm\Ll_\cG(S)$ into $\Mm\Ll_\cG(S)^\#$ by sending
$(F,\lambda)$ to $(F,\eps,\lambda,\eta)$ with $\eps_i=1$ for every $i$
and $\eta_i=1$ for every $i\in V_\Pp$ then
$\beta=\phi_{\Mm\Ll}^\#\circ\beta^\#$ (where
$\phi_{\Mm\Ll}^\#$ is the usual forgetting map).
 \begin{remark}{\rm Before giving the actual definition, we describe
the qualitative idea. Earthquakes paths on $\Tt_\cG(S)$ rebounce when
     reaches a cusp. Since $\Tt_\cG(S)^\#$ is obtained by reflecting
     $\Tt_\cG(S)$ along its faces, it is natural to lift such a paths
     to horizontal paths on $\Tt_\cG(S)^\#$.  Instead of rebouncing
     the enhanced lamination after a cusp is obtained by a reflection
     along a boundary component of the initial lamination.  This
     liftings are unique (up to the choice of a initial signature
     $\eps$) when $F$ does not contain cusp. When $F$ contains a cups
     then there are many possible liftings due to the possible choices
     of the signature of the cusp after the earthquake.  Thus data
     $(F,\eps,\lambda)$ are not sufficient to determines the
     lifting. On the other hand the information of a signature of
     $\lambda$ around the cusp solves this ambiguity.  }\end{remark}
  
Let us come to the actual definition:
\[
\beta_t^\#(F,\eps,\lambda,\eta)=
(\overline F,\overline \eps,\overline\lambda, \overline\eta)
\]
where
\smallskip

(a) Similarly to the definition of the map $\tG$, $(\overline
F,\overline\lambda)= \beta(F,
\rho_\eps(\lambda))$;
\smallskip

(b) $\overline\eps_i=\eps_i\sign (l_{C_i}(F)+t\eta_i {\rm I}_{C_i}(\lambda))$.
\smallskip

(c) $\overline\eta_i=\eta_i\sign (l_{C_i}(F)+t\eta_i {\rm I}_{C_i}(\lambda))$.
\smallskip

Property (1) follows from the fact that $\beta$ is a flow. Point (2)
depends on the fact the spectrum of $\lambda_t$ is constant and the
products $\eps_i(t)\eta_i(t)$ are constant. Point (3) is
straightforward. The only point to check is that $\beta^\#$ is
continuous, as a map
$\mr_\geq0\times\Mm\Ll_\cG^\#(S)\rightarrow\Mm\Ll_\cG^\#(S)$.  By the
definition of the topology of $\Mm\Ll_\cG^\#(S)$ it is enough to show
that for every $\gamma\in\SG$ the functions
\[
(t,(F,\eps,\lambda,\eta))\mapsto l^\#_\gamma(\beta_t^\#(
F,\eps,\lambda,\eta))\qquad (t,(F,\eps,\lambda,\eta))\mapsto
{\rm I}^\#_\gamma(\beta_t^\#(F,\eps,\lambda,\eta))
\]
are continuous.
If $\gamma$ is not peripheral, then $
  l^\#_\gamma(\beta_t^\#(F,\eps,\lambda,\eta))$ and
  ${\rm I}^\#_\gamma(\beta^\#(t, F,\eps,\lambda,\eta)$ depend only on $F$
  and $\lambda$ so the continuity is a consequence of the continuity
  of $\beta$.

If $\gamma$ is peripheral, then by Lemma~\ref{ray-quake} we have
\[
\begin{array}{l}
 l^\#_\gamma(\beta_t^\#(F,\eps,\lambda,\eta))= 
l^\#_\gamma(F,\eps)-t {\rm I}^\#_\gamma(F,\eps,\lambda,\eta)\\
 {\rm I}^\#_\gamma(\beta_t^\#(F,\eps,\lambda,\eta))={\rm I}^\#_\gamma(F,\eps)\,.
 \end{array}
 \]
 
 For every $\xi\in {\rm I}^\#(S)$ let us consider the map of
   $\mr_{\geq 0}\times\Tt_\cG^\#(S)\rightarrow\Tt_\cG^\#(S)$ that
   associates to $t,(F,\eps)$ the projection on $\Tt_\cG^\#(S)$ of
   $\beta_t(F,\eps,\xi(F))$ (where $\xi(F)$ is the realization of
   $\xi$ with respect to the structure given by $F$). By (2) it is a
   flow on $\Tt_\cG(S)^\#$. We will denote by $\Ee^\#_\xi$ the
   homeomorphism of $\Tt_\cG(S)^\#$ corresponding to such a flow at
   time $1$ (notice that $\Ee_\xi\circ\Ee_\xi=\Ee_{2\xi}$), it will be
   called the {\it enhanced earthquake along $\xi$}.
 
{\bf Earthquake Theorem.}
\begin{teo}\label{quake-teo}{\rm [Earthquake Theorem on $\Tt_\cG(S)$]}
For every $F_0,\ F_1 \in \Tt_\cG(S)$, denote by $m$ the number of
points in $V$ that do not correspond to cusp of $F_1$ nor of
$F_2$. Then there exist exactly $2^m$ left earthquakes such that $F_1=
\beta^L_\lambda(F_0)$.  The similar statement holds with respect to
right-quakes.
\end{teo}
This is a consequence of the somewhat more precise
\begin{teo}\label{quake-teo-bis}{\rm [Earthquake Theorem on $\Tt_\cG(S)^\#$]}
For every $(F_0,\epsilon_0),\ (F_1,\epsilon_1) \in \Tt_\cG(S)^\#$,
there is a unique $\xi\in  {\rm I}^\#(S)$ such that
$\Ee_\xi^\#(F_0,\eps_0)=(F_1,\eps_1)$ Similarly for the right quakes.
\end{teo}

Given two ``signed'' surfaces $(F_0,\sigma_0)$ and $(F_1,\sigma_1)$ in
$\Tt_\cG(S)$, where the respective signatures are arbitrary maps
$\sigma_j: V \to \{ \pm 1\}$), we say that they are {\it left-quake
compatible} if there exists a left earthquake
$(F_1,\lambda_1)=\beta^L(F_0,\lambda_0)$ such that $\sigma_j=
\sigma_{\lambda_j}$. The following is an easy Corollary of Lemma
\ref{ray-quake} and of Theorem \ref{quake-teo}.
\begin{cor}\label{nec-comp} The signed surfaces 
$(F_0,\sigma_0)$ and $(F_1,\sigma_1)$ are left-quake compatible if and
only if for every $i=1,\ldots,r$ the following condition is satisfied:
\smallskip

If $l_{C_i}(F_1)< l_{C_i}(F_0)$, then $\sigma_0(i)= 1 \ $. 
If $l_{C_i}(F_1)> l_{C_i}(F_0) \ $, then  $\sigma_1(i)= 1 \ $.
\smallskip

Symmetric statements hold w.r.t. the right-quake compatibility.
\end{cor}

In Section \ref{moreAdS} we will outline an {\it AdS proof} of the Earthquake
Theorem (by following \cite{BSK}) that generalizes Mess's proof in the 
special case of compact $S$.
\medskip
 
{\bf $\Mm\Ll_\cG(S)$ as tangent bundle of $\Tt_\cG(S)$.}
We have seen above that the bundle
$$\pG^\#: \Mm\Ll_\Cc(S)^\# \to \Tt_\Cc(S)^\# $$ shares some properties
with the {\it tangent bundles} $T\Tt_\cG^\#$ of its base space. We are
going to substantiate this fact by means of quake-flows.  In fact we
have associated to every $\xi\in {\rm I}^\#(S) $ a flow of
$\Tt_\cG(S)^\#$, So we can consider the infinitesimal generator of
such a flow, that is a vector field on $\Tt_\cG(S)^\#$, say $X_\xi$.
\begin{prop}
The map
\[
   \Pi:\Tt_\cG(S)^\#\times  {\rm I}^\#(S)\rightarrow T\Tt_\cG(S)^\#
 \]
defined by $\Pi(\xi, F)=X_\xi(F)$ is a trivialization of $T\Tt_\cG(S)$.
\end{prop}

As in the case of compact $S$, it is a consequence of the convexity of the
length function along earthquakes paths.

\begin{remark}{\rm
  The map $\Pi$ is only a {\it topological trivialization}. This means
  that the identifications between tangent spaces arising from $\Pi$
  are not linear.}
\end{remark}

For a fixed type $\theta$, denotes by $ {\rm I}^\#(S)^\theta$ the
points corresponding to laminations that do not enter any cusp. It is
clear that for a point $F\in\Tt^\theta_\cG(S)^\#$ we have that
$X_\xi(F)$ is tangent to $T^\theta_\cG(S)^\#$. So we get that the
restriction of $\Pi$ to $\Tt^\theta_\Cc(S)^\#\times {\rm
I}^\#(S)^\theta$ is a trivialization of $T\Tt^\theta_\cG(S)^\#$.

%%% Local Variables: 
%%% mode: latex
%%% TeX-master: "HAND"
%%% End: 

%% file: HANDWR.tex
\section{Wick rotation-rescaling theory}\label{WR}

We refer to \cite{Be-Bo}. Let $S$ be a base surface of finite
type. Recall from the Introduction and Section \ref{grav},
that $\Mm\Gg\Hh_\kappa(S)$ denotes the Teichm\"uller-like
space of Einstein maximal globally hyperbolic spacetimes of constant
curvature $\kappa = 0, \pm 1$, that contain a {\it complete} Cauchy
surface homeomorphic to $S$.
\smallskip

Denote by $\Pp(S)$ the Teichm\"uller-like space of {\it complex
projective} (that is $(S^2,PSL(2,\mc))$-{\it manifold}) structures on
$S$. Here $S^2$ is the Riemann sphere, identified with $S^2_\infty =
\partial \mh^3$, and $PSL(2,\mc)\cong {\rm Isom^+}(\mh^3)$.
\medskip

The aim of this section is to illustrate the following pattern of
statements (given here in somewhat informal way):
\smallskip

{\bf (Classifications)} {\it For every surface $S$ of finite type, and
every $\kappa = 0,\pm 1$, there are geometrically defined {\rm
``materialization'' maps}
$$ \mG_\Pp: \Mm\Ll(S)\to \Pp(S)$$
$$ \mG_\kappa: \Mm\Ll(S) \to \Mm\Gg\Hh_\kappa(S)$$
that actually make $\Mm\Ll(S)$ an {\rm universal parameter
space}.}

\medskip

{\bf (Canonical correlations)} {\it For every $(F,\lambda)\in
\Mm\Ll(S)$, there are geometrical correlations between the spacetimes
$\mG_\kappa(F,\lambda)$ or between them and the projective surface
$\mG_\Pp(F,\lambda)$. Such correlations are either realized by means
of canonical  {\rm rescalings} or {\rm Wick rotations} directed by 
the respective {\rm cosmological times}, with {\rm universal rescaling
functions}.}
\medskip

Let us explain first some terms entering the last statement.
\smallskip

\begin{defi}\label{WR-resc} {\rm Let $(M,h)$ be any spacetime and
$X$ be a nowhere vanishing $h$-timelike and future directed vector
field on $M$. Let $\alpha, \beta: M\to \mr_+$ be positive functions.
\smallskip

We say that the Riemannian manifold $(M,g)$ is obtained from $(M,h)$
via the {\it Wick rotation} directed by $X$, with {\it vertical}
(resp. {\it horizontal}) rescaling function $\beta$ (resp. $\alpha$),
if for every $y\in M$, the $g$- and $h$-orthogonal spaces to $X(y)$
coincide (denoted it by $\ort{<X(y)>}$), and
$$||X(y)||_g = - \beta(y)||X(y)||_h$$
$$g|_{\ort{<X(y)>}} \ =\ \alpha(y)h|_{\ort{<X(y)>}} \ . $$
\smallskip

Similarly, the spacetime $(M,h')$ is obtained from $(M,h)$ via the
{\it rescaling} directed by $X$, with {\it vertical} (resp. {\it
horizontal}) rescaling function $\beta$ (resp. $\alpha$), if for every
$y\in M$, the $h'$- and $h$-orthogonal spaces to $X(y)$ coincide, and
$$||X(y)||_{h'} = \beta(y)||X(y)||_h$$
$$h'|_{\ort{<X(y)>}} \ =\ \alpha(y)h|_{\ort{<X(y)>}} \ . $$
}
\end{defi}

\subsection {Cosmological time}\label{CT}
We refer to \cite{A} for a general treatment of this
matter.  Here we limit ourselves to recalling the main features of
this notion. Let $(M,h)$ be any spacetime. The {\it cosmological function}
$$\tau: M \to (0,+\infty]$$ is defined as follows: let $C^-(q)$ be the
set of past-directed causal curves in $M$ that start at $q\in M$. Then
$$\tau (q) = \sup\{L(c)|\ c\in C^-(q)\}\ ,$$ where $L(c)$ denotes the
Lorentzian length of $c$. Roughly speaking, this gives the (possibly
infinite) {\it proper time that every event $q \in M$ has been in
existence in $M$}. The function $\tau$ is said {\it regular} if it is
finite valued for every $q\in M$, and $\tau \to 0$ along every
past-directed inextensible causal curve. In such a case it turns that
$\tau$ is a continuous {\it global time} on $M$, called its {\it
cosmological time}.  This cosmological time (if it exists) represents
an intrinsic feature of the spacetime.  Having cosmological time has
strong consequences for the structure of $M$, and $\tau$ itself has
stronger properties (it is locally Lipschitz and twice differentiable
almost everywhere).  In particular: $M$ is globally hyperbolic; for
every $q\in M$, there exists a future-directed time-like unit speed
geodesic ray whose length equals $\tau(q)$. Up to a suitable
past-asymptotic equivalence, these rays form the {\it initial
singularity} of $M$. In a sense $\tau$ gives the Lorentzian distance
of every event from the initial singularity.  \smallskip

\subsection{Grafting and Lorentzian grafting}

Before describing in some formal way how to get a parameterizations of
$\Mm\Gg\Hh_\kappa(S)$ and  $\Pp(S)$ in terms of $\Mm\Ll(S)$, we will explain
how to associate to a pair $(F,\lambda)\in\Mm\Ll(S)$ a projective structure on
$S$ and a spacetime of constant curvature $\kappa$, in some simple cases.

First consider the case $S$ compact closed and $\lambda$ empty. Given
a hyperbolic structure $F=(S,h)$ on $S$, the projective structure
associated to it, that, with a little abuse, we will denote simply by
$F$, is the structure whose developing map coincides (up to
post-composition with $g\in PSL(2,\mc)$) with the isometric developing
map of $F$. Structures obtained in this way are called Fuchsian and
are characterized by the following requirements:

(1) the developing map is injective,

(2) the holonomy is conjugated in $PSL(2,\mr)$.
\smallskip

For the Lorentzian side, define $\mG_\kappa(F)$ to be the spacetime $(S\times
I, g_\kappa)$ where $I$ is the interval $(0,+\infty)$ for $\kappa\geq 0$ and
$I=(0,\pi/2)$ for $\kappa=-1$ and $g_\kappa=g_\kappa(F)$ is so defined
\begin{equation}\label{cc:eq}
g_\kappa=\left\{
\begin{array}{ll}
-dt^2+t^2 h & \textrm{ if }\kappa=0\\
-dt^2+\sh^2(t) h &\textrm{  if }\kappa=1\\
-dt^2+\sin^2(t) h & \textrm{ if }\kappa=-1
\end{array}\right.
\end{equation}

The fact that $g_\kappa$ has constant curvature $\kappa$ is just a local
computation independent of the compactness of $F$. Thus one does the
computation assuming $F=\mh^2$.  For instance, for $\kappa=0$, one embeds
$\mh^2$ in the Minkowski space $\mx_0$ and take the normal evolution of
$\mh^2$ (that is a map $\mh^2\times\mr_\geq 0\rightarrow\mx_0$ sending $(x,t)$
to $tx$) the pull-back of the Minkowski metric takes the form (\ref{cc:eq}).

\begin{remark}{\rm
    Strictly speaking $\mG_{-1}(F)$ is not maximal. In fact the metric
    $g_\kappa$ can be defined as well on the interval $(0,\pi)$.  On the other
    hand, for some reason that will appear clear it is better to define
    $\mG_{-1}(F)$ in this way and then to take its maximal extension.
  }\end{remark}

Now suppose $S$ to be closed and $\lambda$ to be a weighted curve $(c,a)$.
The projective surface $\mG_\Pp(F,\lambda)$ is the grafting of $F$ along
$\lambda$, that we sometimes denote by $Gr_\lambda(F)$.  We cut $F$ along $c$
and grafts a projective annulus $A=c\times[0,a]$ whose developing map can be
explicitly described in the following way. We can choose a developing map
$dev:\tilde F\rightarrow\mh^2$ such that $c$ lifts to a geodesic $\tilde c$
with end-points at $0$ and $\infty$. The developing map of $A$ is given by
\[
     \tilde c\times [0,a]\ni (x,t)\rightarrow dev(x)e^{it}\in \C
     \subset S^2\,.
\]     
The fact that $A$ can be grafted on $F$ descends on the fact that the
developing map of each component of $\partial A$ is conjugated in $PSL(2,\mc)$
to the developing map of $c$.  Notice that $A$ carries a natural Euclidean
metric. The length of each boundary component of $A$ is equal to the length of
$c$ whereas the width of $A$ is equal to $a$.  Thus we can consider on
$Gr_\lambda(F)$ the metric that is hyperbolic on $F\setminus c$ and Euclidean
on $A$. Such a metric is $C^1$ and compatible with the conformal structure
underlying the projective structure of $Gr_\lambda(F)$. We call it the the
Thurston metric of $Gr_\lambda(F)$, in what follows we often indicate with
$Gr_\lambda(F)$ both the projective structure and the metric structure on $S$.

\begin{remark}{\rm
    Thurston distance is defined on every projective structure on $S$ and is a
    metric compatible with the conformal class of the projective surface. The
    interesting point showed by Thurston is that Thurston metric determines
    the projective structure. This means that a map between projective
    surfaces is a projective equivalence iff it is an isometry with respect
    to the corresponding Thusrton distances.}
\end{remark} 

\begin{remark}{\rm
    If $a$ is \emph{small}, then the holonomy group of $Gr_\lambda(F)$, say
    $\Gamma$, is quasi-Fuchsian and the developing map is injective with image
    a component of the discontinuity domain. Thus, $Gr_\lambda(F)$ can be
    regarded as an asymptotic end of the quasi-Fucshian manifold
    $\mh^3/\Gamma$. In fact, the boundary component of the convex core facing
    $Gr_\lambda(F)$ is isometric to $F$ bent along $c$ with bending angle $a$,
    and the annulus $A$ coincides with the set of points in $Gr_\lambda(F)$
    that are sent by the retraction on the convex core to the bending line.
    
    Moreover, let us consider the component of the complement of the
    convex core in $\mh^3/\Gamma$, facing $Gr_\lambda(F)$. Then the
    distance $d$ from the convex core is a $\mathrm C^1$ function on
    it whose level surfaces are isometric to $\ch d\cdot Gr_{\tgh
    d\lambda}(F)$ (if $X$ is a metric space $\lambda\cdot X$ denotes
    the metric space obtained by multiplying the distance by
    $\lambda$).
    
    Thurston generalized this idea and showed how to associate to each
    projective structure on $S$ a non-complete hyperbolic structure on
    $S\times (0,1)$, called the $H$-hull such that \smallskip\par\noindent (1)
    its completion is $S\times[0,1)$ and $S\times\{0\}$ is a locally convex
    bent surface $F$ along a lamination $\lambda$ \smallskip\par\noindent (2)
    the asymptotic end $S\times\{1\}$ carries the original projective surface
    that in tunrs coincides with $Gr_\lambda(F)$.  \smallskip\par Moreover the
    distance $d$ from $S\times\{0\}$ is a $C^1$ function and level surfaces are
    isometric to 
\begin{equation}\label{dist:eq}
 \ch d\cdot Gr_{\tgh d\lambda}(F)\,.
\end{equation}
 Clearly in the
    quasi-Fuchsian case the $H$-hull is simply the end of the corresponding
    quasi-Fuchsian manifold facing the projective surface.  }\end{remark}

Consider now the Lorentzian case.

To construct $\mG_\kappa(F,\lambda)$ we will deform the structure on
$\mG_\kappa(F)$ by means of a construction that is reminescent of the grafting
procedure, so we call it the \emph{Lorentzian grafting}.

With a little abuse let us denote by $c$ the geodesic representative of $c$
with respect to the hyperbolic structure $F$.  Then one shows that the
timelike surface $c\times I$ is totally geodesic in $\mG_\kappa(F)$ (it is
still a local computation -  for instance, in the flat case it is a direct
consequence of the fact that geodesics of $\mh^2$ are intersection of $\mh^2$
with linear time-like planes of Minkowski space).  Then one cuts
$\mG_\kappa(F)$ along $c\times I$ and \emph{grafts} a piece, say
$\mG_\kappa(A)$ such that \medskip

(1) topologically $\mG_\kappa(A)=(c\times[0,a])\times I$ that is a the product
of the annulus $A=c\times[0,a]$ by the time interval $I$.  \smallskip

(2) the restriction of the metric on each slice $A\times\{t\}$ is a Euclidean
annulus, whose width depends only on $a$ and on $t$ and whose boundary length
is equal to the length of $c\times\{t\}\subset\mG_\kappa(F)$.  \smallskip

(3) the boundary of $\mG_\kappa (A)$ (that is $\partial A\times
I=c\times I\times\{0,a\}$) is totally geodesic and each component is
isometric to $c\times I$
\medskip

For instance in the flat case $\mG_0 (A)$is just $(c\times I)\times [0,a]$
with the product metric (that is flat since it is the product of two flat
metrics).  For the other curvatures, the expression of the metric on
$\mG_0(A)$ takes the more complicated form given by
\begin{equation}\label{bo:eq}
\left\{\begin{array}{ll}
-dt^2+(\ch^2(t) dr^2 +\sh^2(t)  d\theta^2) & \textrm{ for }\kappa=1\\
-dt^2+(\cos^2(t) dr^2 +\sin^2(t) d\theta^2) & \textrm{ for }\kappa=-1
\end{array}\right.
\end{equation}
 where $\theta$ is an arc parameter on $c$ and $r$ is the variable on
$[0,a]$.  Notice that the width of $A\times\{t\}$ is independent of
$t$ only in the flat case.

\begin{remark}{\rm
    The piece $\mG_{-1}(A)$ is well-defined only for $t\in(0,\pi/2)$ and this
    explains the definition of $\mG_{-1}(F)$.  In general the spacetime
    obtained for $\kappa=-1$ is never maximal, so more correctly $\mG_{-1}(F)$
    will denote the maximal extension of the spacetime we have defined.  In
    the next sections we will explain the reason of this asymmetry and also
    how the spacetime we have defined is uniquely determined by its maximal
    extension.  }\end{remark}

\begin{remark}{\rm
    A way to define $\mG_\kappa(F,\lambda)$ for a generic $\lambda$ is by
    means of an approximation argument.  We take a sequence of simple weighted
    curves $\lambda_n=(c_n,a_n)$ converging to $\lambda$ and define
    $\mG_\kappa(F,\lambda)=\lim \mG_\kappa(F,\lambda_n)$. Clearly the
    existence of this limit has to be checked: to this aim it is better to
    work in the framework of $(X,G)$-structures and study the behaviour of the
    developing maps of $\mG_\kappa(F,\lambda)$.  This will be the theme of the
    next sections.  }
\end{remark}

Notice that the construction of $\mG_\kappa(F,\lambda)$ gives, as a
by-product, a natural foliation of the spacetime in spacelike surfaces
homeomorphic to $S$.  In fact in both $\mG_\kappa(F)$ and $\mG_\kappa(A)$ we
have pointed out a time-function $t$ to express the metric in some explicit
way, these functions glue to a time-function on $\mG_\kappa(F,\lambda)$.
Notice however that the function $t$ in $\mG_\kappa(F,\lambda)$ is not smooth:
its level surfaces are made by hyperbolic pieces and Euclidean annuli. In fact
they are reminescent of the usual grafted surfaces.

Let us consider the flat case. In such a case the $t$ level surface
corresponding to some value $t_0$ is obtained by multiplying the hyperbolic
metric on $F$ by the factor $t_0^2$ (that is by multiplying the hyperbolic
distance by the factor $t_0$), by cutting along $c$ and gluing a Euclidean
annulus of width $a$.  This is the same as grafting an annulus of width $a/t_0$
on $F$ and then multiplying the grafted distance by the factor $t_0$.

More generally one can check explicitly that for a weighted multi-curve
$\lambda=(c,a)$ the surface $t^{-1}(t_0)\subset\mG_\kappa(F,\lambda)$ is
metrically equal to
 \begin{equation}\label{CT:eq}
 \begin{array}{ll}
   t_0 \cdot Gr_{\lambda/t_0}(F) & \textrm{ if } \kappa=0\\
   \sh t_0 \cdot Gr_{\lambda/\tgh t_0}(F) & \textrm{ if } \kappa=1\\
   \sin t_0 \cdot Gr_{\lambda/\tan t_0}(F) & \textrm{ if }\kappa=-1
\end{array}
\end{equation}

The point that makes this remark interesting is that the function $t$ is the
cosmological time of $\mG_\kappa(F,\lambda)$, so it is somehow independent of
the parameterization and the same formulae to express the level surface work
for every $(F,\lambda)$.  This remark motivates the idea to find a canonical
rescaling directed by the gradient of the cosmological time transforming
$\mG_0(F,\lambda)$ into $\mG_{\pm 1}(F,\lambda)$ and a Wick Rotation
transforming $\mG_0(F,\lambda)$ into the $H$-hull of $Gr_\lambda(F)$.

\begin{remark}{\rm
    Consider the case $S$ of finite type. For $F\in\Tt_\Cc(S)$ we could try to
    define $\mG_\kappa(F)$ as in the closed case.  Notice however that the
    slice $S\times\{t\}$ is in general not complete. In fact such a spacetime
    have a natural totally geodesic timelike boundary that is homeomorphic to
    $\partial F^\Cc\times I$. A way to get a complete level surface is then
    for each boundary component $c$ of $F$ to glue a piece
    $\mG_\kappa(\Delta)$ where $\Delta=c\times [0,+\infty)$ is a annulus with
    infinite width and $\mG_\kappa(\Delta)=\Delta\times I$ with a metric given
    in~(\ref{bo:eq}). Notice that the definition of $\mG_\kappa(F)$ is then
    consistent with the previous case provided to allow boundary component of
    $F$ to carry an infinite weight.
    
    In fact one can show that to define $\mG_\kappa(F,\lambda)$ it is
    necessary to glue this cylindrical ends for every boundary component of
    $F$ that is not close to the lamination. On the other hand, if $\lambda$
    contains a leaf $l$ spiraling around a boundary curve, it is clear that it
    is possible to define the analogous of $\mG_\kappa(A)$ for this leaf (that
    now is the product of a infinite band of width equal to the weight of $l$
    and the time-interval $I$) and apply the grafting procedure. Notice that
    if $l$ spirals around a boundary component $c$, the corresponding end on
    the slice $S\times\{t\}$ in $\mG_\kappa(F,\lambda)$ appears complete (in
    fact a path entering the ends meets the band infinite times so its length
    cannot be bounded).
    
    From this discussion it appears clear that in this context it is
    more convenient to use the notion of geodesic lamination on a
    surface given in Remark ~\ref{equi-def}.  That is we require that the
    boundary components of $F$ are contained in the lamination and
    that paths arriving on the boundary have infinite total mass.  In
    particular for each boundary component either a leaf spirals
    around it or it carries an infinite weight.  With this definition
    the $0$ lamination on $F$ is obtained by putting the weight
    $+\infty$ on each boundary component.  }
\end{remark}

\subsection{Wick rotation - rescaling set up}
Let us go back to the statement concerning the canonical
correlations. We will see that every spacetime $\mG_\kappa(F,\lambda)$
has (rather tame) cosmological time, so that the geometry of the
initial singularity will quite naturally arise. The above mentioned
Wick rotation-rescaling (possibly only defined on suitable ``slabs''
of the spacetimes) will be directed by the {\it gradient of the
cosmological times}. The rescaling functions will be {\it universal}
in the sense that their values only depend on the cosmological time
values: for every $y$ in the domain of definition, $\beta(y) =
\beta(\tau(y))$, $\alpha(y) = \alpha(\tau(y))$. We stress that they do
not depend on $(F,\lambda)$.
\medskip

We are going to outline the {\it linked} auguments establishing both 
the constructions of the maps $\mG_*$, the geometric correlations 
and the fact that the materialization maps induce bijections.
\smallskip

Let $(F,\lambda)\in \Mm\Ll(S)$. With the notations of Section
\ref{ML(S)}, we have
$$F\subset F^\Cc \subset \hat {F} = \mh^2/\Gamma \ .$$
$(F^\Cc,\lambda)$ lifts to a $\Gamma$-invariant couple
$(H,\tilde{\lambda})$ where $H$ is a straight convex set in $\mh^2$
equipped with the measured geodesic lamination $\tilde{\lambda}$. The
universal covering map $$\mh^2 \to \hat {F}$$ restricts to the
universal covering maps $H \to F^\Cc$, $\mathring {H}\to F$, where
$\mathring {H}$ is the interior of $H$. To simplify the notations we
make the abuse of always writing $\lambda$ instead of either
$\tilde{\lambda}$, $(F,\lambda)$ or $(\mathring {H},\tilde{\lambda})$,
that is we will understand $F$ or $\mathring H$.
\smallskip

The projective surface 
$$S^\lambda_\Pp = \mG_\Pp(\lambda)$$ will be given in
terms of a couple $(d^\lambda_\Pp, h^\lambda_\Pp)$ of compatible
developing map
$$d^\lambda_\Pp: \widetilde{S} \to S^2$$
and holonomy representation
$$h^{\lambda}_\Pp: \pi_1(S) \to PSL(2,\mc) \ .$$ We denote
$$p^\lambda_\Pp:\tilde S^\lambda_\Pp \to S^\lambda_\Pp$$ the
corresponding local isomorphic universal covering.
\smallskip

Similarly, every spacetime $$Y^\lambda_\kappa=\mG_\kappa(\lambda)$$ will
be specified by a compatible couple
$(d_\kappa^\lambda,h^\lambda_\kappa)$,
$$d_\kappa^\lambda: \widetilde S\times \mr \to \mx_\kappa$$
$$h^\lambda_\kappa: \pi_1(S) \to {\rm Isom}^+(\mx_\kappa) \ . $$
We denote
$$p^\lambda_\kappa:\Uu_\kappa^\lambda \to Y^\lambda_\kappa $$ 
the corresponding local isomorphic universal covering.
\smallskip

For simplicity, we will often identify $S$ with $F$, $\widetilde S$
with $\mathring{H}$, $\pi_1(S)$ with $\Gamma$, and so on.
\smallskip

For every $F$ as above, denote $\lambda_0$ the measured geodesic
lamination just consisting of the $+\infty$ weighted boundary
components of $F^\Cc$. Recall that $\lambda_0$ is the initial
end-point of any ray in $\Mm\Ll(F)$. We will describe explicitely the
corresponding surface $S^0_\Pp$ and spacetimes $Y^0_\kappa$,
$\Uu^0_\kappa$. Every $\lambda \in \Mm\Ll(F)$ somehow encodes the
instructions in order to {\it deform} $(d^0_*,h^0_*)$ towards
$(d^\lambda_*,h^\lambda_*)$ as it has been make  explicit in the case
of finite laminations.
\medskip

\subsection{Flat spacetimes classification}\label{flat} 
Take the hyperboloid model $\mh^2 \subset \mx_0$ of the hyperbolic
plane. The {\it chronological future} of $0$ in $\mx_0$ is the cone 
$I^+(0)= \{-x_0^2 + x_1^2 + x_2^2 < 0, \ x_2>0\}$ from $0$ over
$\mh^2$.  $I^+(0)$ has cosmological time $\tau = (x_0^2 -
(x_1^2+x_2^2))^{1/2}$, so that $\mh^2 = \{\tau = 1\}$ and $0$ is the
initial singularity.  The future $I^+(r)=\{ -x_0^2 - x_2^2 <0, \ x_2>0 \}$ of
the spacelike geodesic $r= \{x_0=x_2=0\}$ has cosmological time $\tau
= (x_0^2 - x_2^2)^{1/2}$; $r$ is the initial singularity.
\smallskip

{\bf Construction of $\Uu^0_0$.}  The cone $C_0H$ from $0$ over
$H\subset \mh^2$ is contained in $I^+(0)$.  The boundary of $C_0H$ is
made by the cone over the boundary of $H$. Each component of $\partial
C_0H$, corresponding to a geodesic line $\gamma \subset \partial H$,
is the intersection with $I^+(0)$ of a hyperplane $P_\gamma$,
orthogonal to a determined unitary spacelike vector $v_\gamma$, that
points out of $C_0H$. The developing map $d^0_0$ is an {\it embedding} onto
the convex domain $\Uu^0_0$ of $\mx_0$ made by the union of $C_0H$
with the future of all the rays of the form $\{tv_\gamma + x|t\geq 0\}$.
\smallskip

A convenient description of the domain $\Uu^0_0$ is as the
intersection of half-planes.  In fact we have
\[
   \Uu^0_0=\bigcap_{x\in H}\fut(x^\perp)
\]   
This shows that $\Uu^0_0$ is convex and future complete.

Up to isometry of $\mx_0$, the local model for $\Uu^0_0$ at each
component of $\partial C_0(H)$, is the future $I^+(r_+)$ of the ray
$\{x_1\geq 0\}\subset r$, that is
$$I^+(r_+)= (I^+(0)\cap \{ x_1\leq 0\}) \cup (I^+(r)\cap \{x_1 \geq
0\} \ .$$ The above cosmological times match at the intersection,
producing the cosmological time of the union, that turns to be a C$^1$
function.  The ray $r_+$ is the initial singularity. Then $\Uu^0_0$
has cosmological time that coincides with the one of $I^+(0)$ on
$C_0H$; the initial singularity is the spacelike tree made by one
vertex at $0$ and the rays $tv_\gamma$, $t\geq 0$, emanating from the
origin. The action of $\pi_1(S)$ on $H$ naturally extends to the whole
of $\Uu_0^0$, giving the holonomy $h^0_0$.
\smallskip

{\bf Construction of $\Uu^\lambda_0$.}  Let us consider now an
arbitrary lamination $\lambda=(L,\mu) \in \Mm\Ll(F)$.  The developing
map $d^\lambda_0$ will be always an {\it embedding} onto a convex
domain $\Uu^\lambda_0$ in $\mx_0$, obtained as follows. Fix a
base-point $x_0\in\mathring{H}$ not belonging to the weighted part
$L_W$ of $\lambda$. For every $x\in\mathring H \setminus L_W$ choose
an arc $c$ transverse to $\lambda$ with end-points $x_0$ and $x$. For
$t\in c\cap L$, let $v(t)\in \R^3$ denote the unitary spacelike vector
tangent to $\mh^2$ at $t$, orthogonal to the leaf through $t$ and
pointing towards $x$. For $t\in c\setminus L$, let us set $v(t)=0$.
In this way we define a function
\[
    v:c\rightarrow\mr^3
\]
that is continuous on the support of $\mu\,$.
We can define
\[
s(x)=\int_{c}v(t)\d\mu(t).
\]
It is not hard to see that $s$ does not depend on the path $c$.
Moreover, it is constant on every stratum of the stratification
determined by $\lambda$, and it is a continuous function on
$\mathrm{H}\setminus L_W$.

The domain $\Uu_0^\lambda$
can be defined in the following way
\[
  \Uu_0^\lambda=\bigcap_{x\in\mathrm{ H}\setminus
  L_W}\fut(s(x)+\ort{x})\ 
\]
where $\ort{x}$ denote the orthogonal 2-plane to $x$ in $\mx_0$.
Note that this definition is compatible with the one already given for
$\Uu^0_0$.

The holonomy of $Y_0^\lambda $ can be defined in this way:
\[
    h_0^\lambda(\gamma)=h_0^0(\gamma)+\tau(\gamma)
\]
where $h_0^0:\pi_1(S)\rightarrow SO(2,1)$ is the hyperbolic holonomy
of $F$ and $\tau(\gamma)$ is the translation by the vector $s(\gamma
x_0)$.  Since the lamination $\lambda$ is $h$-invariant (being the
pull-back of a lamination on $F$) the domain $\Uu_0^\lambda$ turns to
be $h_0^\lambda$-invariant and $Y_0^\lambda$ is the quotient of
$\Uu_0^\lambda$ by this action.

Let us summarize the main properties of this constructions (see
\cite{Be-Bo, Ba} for all details).

\begin{teo}\label{flat2} 
(1) $\Uu^\lambda_0$   coincides with the
intersection of the future of its null support planes. In particular it
is future complete.
\smallskip

(2) $\Uu^\lambda_0$ has C$^1$ cosmological time $T^\lambda_0$ with
range $(0,+\infty)$. Every level surface
$\Uu^\lambda_0(a)=(T^\lambda_0)^{-1}(a)$ is a {\rm complete} Cauchy
surface. For every $y\in \Uu^\lambda_0$, there is a unique past
timelike geodesic segment emanating from $y$ that realizes
$T^\lambda_0(y)$. The union of the past end-points of such segments
makes the {\rm initial singularity} $\Sigma^\lambda_0$. This is a
spacelike {\rm $\R$-tree} injectively immersed in $\mx_0$.
\smallskip

(3) The action of $\pi_1(S)$ on $H$ induces a natural flat spacetime
holonomy action on $\Uu^\lambda_0 \cup \Sigma^\lambda_0$. The
cosmological time descends to the quotient spacetime $Y^\lambda_0$.
\end{teo}

It is convenient to give a general definition of convex subset of
$\mx_0$ satisfying statement (1) in this theorem.

\begin{defi}\emph{
    A {\it regular domain} is an open convex subset of $\mx_0$ that
    coincides with the intersection of the future of its null support
    planes and admits at least two non-parallel null support planes.
    }\end{defi}

Hence we have a well defined map
$$\mG_0: \Mm\Ll(S)\to \Mm\Gg\Hh_0(S) \ .$$

The spacetimes $\Uu^\lambda_0$ (and $Y^\lambda_0$) are particularly
simple to figure out when $\lambda$ is a {\it finite} lamination. In
such a case, the local model consists of the future, say $\Uu_0$, of a
segment $I=[0,\alpha_0v_0]$, where $v_0$ is a unitary spacelike vector
and $0<\alpha_0<\pi$.  Here local model means that there is a
neighbourhood of each point $p\in Y^\lambda_0$ that embeds in $\Uu_0$
via an isometry that preserves the cosmological time.

The cosmological time on $\Uu_0$ is realized by geodesics with starting point
on $[0,\alpha_0v]$, so there is a natural projection say
$r:\Uu_0\rightarrow[0,\alpha_0v]$ sending $p$ to the point on the segment that
relaizes the cosmological time.

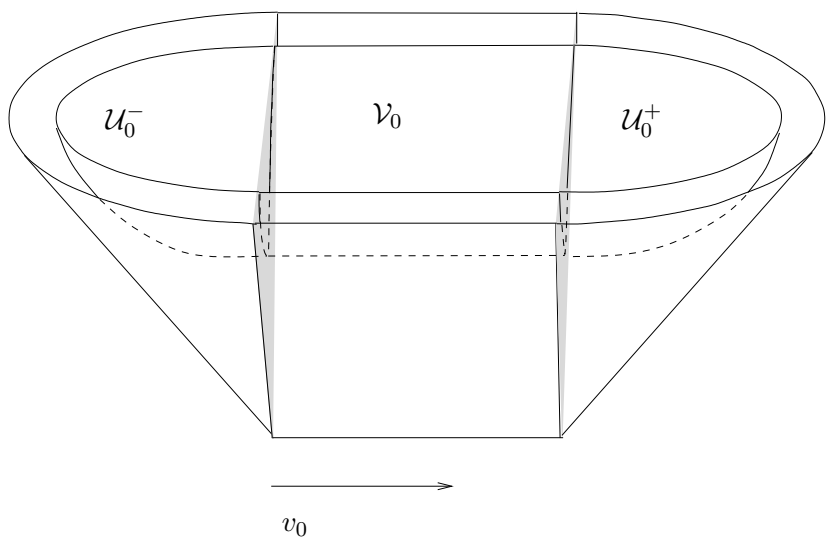
\begin{figure}[h!]
\begin{center}
\input{1geo-flat.pstex_t}
\caption{The domain $\Uu_0$, its decomposition, and a level
    surface.}
\end{center}
\end{figure}

We have a decomposition of $\Uu_0$ in three pieces
$\Uu_0^-,\Uu_0^+,\Vv$ defined in the following way:
\[
\begin{array}{l}
\Uu_0^-=r^{-1}(0)\,;\\
\Vv=r^{-1}(0,\alpha_0 v_0)\,;\\
\Uu_0^+=r^{-1}(\alpha_0 v_0)\, .
\end{array}
\]

We denote by $\Uu_0^+(a),\Uu_0^-(a),\Vv(a)$ the intersections of
corresponding domains with the surface $\Uu_0(a)$. Surfaces
$\Uu_0^+(a)$ are hyperbolic of constant curvature $-1/a^2$.  On the
other hand, the parametrization of $\Vv$ given by
\[
  (0,\alpha_0)\times l_0\ni(t, y)\mapsto a y+t v_0\in \Vv(a)
\]
produces two orthogonal geodesic foliations on $\Vv$. The
parametrization restricted to horizontal leaves is an isometry,
whereas on the on vertical leaves it acts as a rescaling of factor
$a$. Thus $\Vv(a)$ is a Euclidean band of width $\alpha_0$.
Note that by formally setting $\alpha_0= +\infty$, and removing 
$\Uu_0^+$, we recover the above local model for $\Uu^0_0$ at each
component of $\partial C_0(H)$.\\

The initial singularity of a flat spacetime corresponding to a finite
lamination is a {\it simplicial} metric tree.  On the other hand, in
\cite{Be-Bo} we prove also a suitable {\it continuous dependence} of
$\Uu^\lambda_0$ on $\lambda$. By using the density of finite
laminations, this implies that spacetimes corresponding to finite
laminations provide us with good approximations of arbitrary ones.
\smallskip

%\begin{remark}{\rm
%The spacetime $\Uu_0$  provides a local model of $\mG

{\bf Asymptotic states.} In general, the cosmological time level
surface $\Uu^\lambda_0(1)$ ($Y^\lambda_0(1)$) is a C$^1$ spacelike
surface; with the induced Riemannian metric it realizes the {\it
grafting} of $\mathring {H}$ (the hyperbolic surface $F$) at the
measured geodesic lamination $\lambda$. By taking the rescaled level
surface $(1/a)Y^\lambda_0(a)$, we get a 1-parameter family of grafting
of $F$. More precisely we get that $(1/a)Y^\lambda_0(a)$ is obtained
by grafting $F$ along $\lambda/a$.

When $a\to +\infty$ the geometries of $(1/a)Y^\lambda_0(a)$
converge to $F$. The geometry of the initial singularity
$\Sigma^\lambda_0$ of $\Uu^\lambda_0$, together with the isometric
action of $\Gamma$ on it, is ``dual'' to the geometry of
the measured lamination $\lambda$, and can be recovered by means of
the asymptotic behaviour of the level surfaces $\Uu^\lambda_0(a)$
(equipped with the respective isometric actions of $\Gamma$ on them),
when $a\to 0$.
 
\smallskip

{\bf The inverse map of $\mG_0$.}  The image of $\mG_0$ consists of
spacetimes whose universal covering is a regular domain that is, in
particular, future complete. On the other hand, general results due to
Barbot \cite{Ba} on flat spacetimes, applied in our finite type
situation, imply that, {\it possibly reversing the time orientation},
every spacetime $Y$ in $\Mm\Gg\Hh_0(S)$ is future complete, and its
universal covering is a regular domain $\Uu \neq I^+(r)$. So it is
natural to consider the quotient $\Mm\Gg\Hh_0(S)/\pm$, up to time
orientation reversing. We are going to outline the steps leading to
the inverse map of $\mG_0$, defined on it. First one shows that every
regular domain $\Uu$ has cosmological time $T$ that satisfies point
(2) of Theorem \ref{flat}.  We consider the level surface $\Uu(1)$. We
have a natural continuous {\it retraction}
$$r: \Uu(1)\to \Sigma_\Uu $$ onto the initial singularity. Moreover,
the gradient of $T$ is a unitary vector field, hence it induces the
{\it Gauss map}
$$N: \Uu(1)\to \mh^2 \ .$$ The closure $H_\Uu$ of the
image of $N$ in $\mh^2$ is a straight convex set.  If $\Uu \to
Y$ is a universal covering of $Y\in \Mm\Gg\Hh_0(S)$, the action of
$\pi_1(S)$ extends to $H_\Uu$, and makes it a universal covering of
$F_\Uu^\Cc$, for some $F_\Uu \in \widetilde{\Tt}(S)$.  We take the
partition of $\Uu(1)$ given by the closed sets $r^{-1}(y)$, $y\in
\Sigma_\Uu$. Via the retraction, we can pullback to this partition the
metric structure of $\Sigma_\Uu$, and (in a suitable sense) we can
project everything onto $H_\Uu$, by means of the Gauss map.  More
precisely, if $r^{-1}(y)$ is 1-dimensional, then it is a geodesic
line, so that the union of such lines makes a lamination in
$\Uu(1)$. We can define on it a transverse measure such that the mass
of any transverse path is given by the integral of the Lorentzian norm
of the derivative of $r$. A measured geodesic lamination $\lambda_\Uu$
on $H_\Uu$ is obtained via the push-forward by $N$ of this lamination
on $\Uu(1)$.  This descends to a lamination $\lambda_\Uu$ on $F$. So
we eventually get $\mG_0^{-1}(Y)=(F_\Uu,\lambda_\Uu)$.  This achieves
our classification of {\it flat} $MGH$ spacetimes of finite type.

\subsection {Wick rotation: flat Lorentzian vs 
hyperbolic geometry}\label{theWR} Although we adopt a slightly different
definition of the involved measured geodesic laminations, the
bijective map
$$ \mG_\Pp: \Mm\Ll(S)\to \Pp(S)$$ is due to Kulkarni-Pinkall \cite{Ku}
and extends one due to Thurston for compact $S$. This is unfolded in
terms of a 3-dimensional hyperbolic construction. We are going to describe
it, by performing at the same time the canonical Wick rotation establishing
the correlation between the flat spacetimes $Y \in \Mm\Gg\Hh_0(S)$
and suitable hyperbolic 3-manifolds.
\smallskip

For every  $Y \in \Mm\Gg\Hh_0(S)$ ($Y=Y^\lambda_0$), with universal
covering $\Uu \to Y$, and cosmological time $T$, we construct
a local $\mathrm C^1$-diffeomorphism
\[
   d_\mh: \Uu(> 1)\rightarrow\mh^3
\]

and a compatible holonomy

\[
   h_\mh: \pi_1(S) \to PSL(2,\C)
\] 

realizing a (non complete) hyperbolic structure $M=M_Y$ on $Y(>1)$.
This verifies the following properties:
\smallskip

(1) The hyperbolic metric of $M$ is obtained by the Wick rotation
of the flat Lorentzian metric on $Y(>1)$, directed by the gradient of
$T$, with universal rescaling functions:

\[
       \alpha =  \frac{1}{T^2-1} \ , \qquad\qquad \beta=\frac{1}{(T^2-1)^2} \ .
\]

(2) Recall that the closure $H$ of the Gauss map image is the straight
    convex set realizing the future asymptotic geometry of $\Uu$. Then
    the map $d_\mh$ extends (in $h_\mh$ equivariant way) to
\[
   d_\mh: \Uu(\geq 1)\cup H \rightarrow\mh^3
\]
such that:
\smallskip

(a) The restriction of $d_\mh$ to $\Uu(>1)\cup H$ corresponds to the
completion of the manifold $M$. The restriction to $\mathring {H}$ is
a locally isometric {\it pleated immersion} in $\mh^3$, having the
measured geodesic lamination $\lambda$ as {\it bending} locus. This
gives the so called {\it hyperbolic boundary} of $M$. The level
surfaces $\Uu(a)$, $a>1$, correspond via $d_\mh$ to level surfaces of
the distance function $\Delta$ on $M$ from its hyperbolic boundary, so
that the {\it inverse} Wick rotation is directed by the gradient of
$\Delta$. More precisely the following formula holds:
\[
   \Delta=\arctgh (1/T) \ .
\]
\smallskip

(b) The restriction $d_\mh|_{\Uu(1)}$ actually coincides with
the developing map of the complex projective structure $S_\Pp^\lambda$
\[
   d_\Pp: \Uu(1)\rightarrow S^2
\]
so that

\[
  h_\Pp =  h_\mh \ .
\]
The spacelike metric of $\Uu(1)$ ($Y(1)$) coincides with the {\it
Thurston metric} of this projective surface, as well as its {\it
canonical stratification} coincides with the stratification induced by
the retraction $r$ of $\Uu(1)$ onto the initial singularity. This
gives the so called {\it asymptotic complex projective boundary} of
$M$. In fact $M$ turns to be the $H$-hull of $Y(1)$.

\begin{remark}\label{thurston_metric}{\rm We recall here the definition
of the above mentioned ``Thurston metric'' and ``canonical stratification''. 
Let us take a complex projective structure on our  surface $S$ and consider a
developing map
\[
    D:\tilde S\rightarrow S^2 \ .
\]
Pulling back the standard unitary-sphere metric of $S^2$ on $\tilde S$
is {\it not} a well-defined operation, as it depends on the choice of
the developing map.  Nevertheless, by the compactness of $S^2$, the
completion $\overline S$ of $\tilde S$ with respect to such a metric
is well-defined. It turns out that in our finite-type situation,
$\overline S\setminus\tilde S$ contains at least $2$ points (we say
that it is of \emph{hyperbolic} type). A \emph{round disk} in $\tilde
S$ is a set $\Delta$ such that $D|_\Delta$ is injective and the image
of $\Delta$ is a round disk in $S^2$ (this notion is well defined
because $PSL(2,\mc)$ sends round disks onto round disks).  Given a
\emph{maximal} disk $\Delta$ (with respect to the inclusion), we can
consider its closure $\overline\Delta$ in $\overline S$.\par
$\overline\Delta$ is sent by $D$ to the closed disk
$\overline{D(\Delta)}$.  In particular, if $g_\Delta$ denotes the
pull-back on $\Delta$ of the standard {\it hyperbolic} metric on
$D(\Delta)$, we can consider the boundary of $\Delta$ in $\hat S$ as
its ideal boundary. Since $\Delta$ is maximal, $\overline\Delta$ is
not contained in $\tilde S$. So, if $\Lambda_\Delta$ denotes the set
of points in $\overline\Delta\setminus\tilde S$, let $\hat \Delta$ be
the convex hull in $(\Delta, g_\Delta)$ of $\Lambda_\Delta$ (by
maximality $\Lambda_\Delta$ contains at least two points). In
\cite{Ku} it is proved that for every point $p\in\tilde S$, there
exists a unique maximal disk $\Delta$ containing $p$ such that
$p\in\hat\Delta$. So, $\{\hat\Delta|\Delta\textrm{ is a maximal
disk}\}$ is a partition of $\tilde S$. We call it the \emph{canonical
stratification} of $\tilde S$. Clearly the stratification is invariant
under the action of $\pi_1(S)$.

Let $g$ be the Riemannian metric on $\tilde S$ that coincides at $p$
with the metric $g_\Delta$, where $\Delta$ is the maximal disk such
that $p\in\hat\Delta$.  It is a conformal metric, in the sense that it
makes $D$ a conformal map.  It is $\mathrm C^{1,1}$ and is invariant
under the action of $\pi_1(S)$. So, it induces a metric on $\tilde S$. We
call it the \emph{Thurston metric} on $\tilde S$.

Finally let us recall the construction of the $H$-hull of $S$. For 
$p\in\tilde S$, let $\Delta(p)$ be the maximal disk such that
$p\in\hat\Delta$. The image of $\Delta$ via $dev$ is a round disk in $S^2$, so
its boundary is the trace of a hyperbolic plane $P(p)$ in $\mh^3$. Let $c_p$
the geodesic half-line with an end-point at $dev(p)$ an end-point on $P(p)$
and orthogonal to $P(p)$. Then the developing map of the $H$-hull of $S$ is
the map
\[
    \tilde S\times(0,+\infty)\ni (p,t)\mapsto c_p(t)\in\mh^3\,.
\]
Notice that if $S$ is quasi-Fuchsian, the $H$-hull is simply the end of the
corresponding quasi-Fuchsian manifold facing $S$.}
\end{remark}

{\bf About the rescaling functions.}  Before proving the Theorem we
want to give some euristhic motivation for formulae of Wick Rotation.
The point is that we want to construct a Wick Rotation transforming
$Y_0^\lambda$ (or some slab) into the $H$-hull, say $H$, of
$Gr_\lambda(F)$, in such a way that the CT level surfaces are sent to
level surfaces of the distance from the hyperbolic boundary and
rescaling functions are constant on level surfaces.  Now suppose that
such a Wick Rotation exists.  Let $\Delta(T)$ be such a way that the
Wick Rotation transforms $Y(T)$ into $M(\Delta(T))$, and let
$\alpha(T)$ and $\beta(T)$ be the horizontal and vertical rescaling
functions.

By formulae~(\ref{dist:eq}) and (\ref{CT:eq}) we should have
\[
   (\alpha(T))^{1/2} T Gr_{\lambda/T}(F)=\ch\Delta(T) Gr_{\tgh\Delta(T)\lambda}
\]
Since $Gr_{t\lambda}(F)$ is conformally equivalent to $Gr_{s\lambda}(F)$ iff
$s=t$ we deduce that 
\[
   T=1/\tgh(\Delta(T))
\]
that is $\Delta(T)=\arctgh (1/T)$. Moreover we have
\[
  \alpha(T)= \ch^2(\Delta(T))/T^2=1/(T^2-1)\,.
\]

Finally, let $X$ denote the gradient of $T$ with respect to the flat
metric and $Y$ denote the gradient of $\Delta$ with respect to the hyperbolic
metric. We have $X=-\beta^{1/2} Y$. On the other hand 
$\langle X, Y\rangle_{Hyp}= d\Delta(X)=\Delta'(T) d T(X)=-1/(T^2-1)$. 
Thus $\beta(T)=1/(T^2-1)^2$.

Summarizing, if some Wick Rotation exists satysfying properties 
we have required, then necessarily $\alpha=1/(T^2-1)$ and $\beta=1/(T^2-1)^2$.

{\bf Bending cocycle}
A key step in the construction is the {\it bending} of $\mathring H$ in $\mh^3$
along a measured geodesic lamination $\lambda$. We mostly refer to the
Epstein-Marden paper \cite{Ep-M} where this hyperbolic bending has
been carefully studied (in the case of $\mathring H = \mh^2$; however
the constructions extend straightforwardly to the general case). In
fact in \cite{Ep-M} one considers {\it quake-bend} maps, more
generally associated to {\it complex-valued} transverse measures on a
lamination $\Ll$. Bending maps correspond to imaginary valued
measures. So, given a measured geodesic lamination $\lambda =
(\Ll,\mu)\in \Mm\Ll(F)$, we take $i\mu$ in order to get the
corresponding bending map.
\smallskip

{\bf The bending cocycles.} We fix once and for all an embedding of
$\mh^2$ into $\mh^3$ as a totally geodesic hyperbolic plane.

Given $\lambda$ on $H$ as usual, we define first the associated {\it
bending cocycle} (recall a similar notion already introduced in
Section \ref{ML(S)} relatively to the  earthquakes).  This is a map
\[
    B_\lambda:\mathring H \times\mathring H \rightarrow PSL(2,\mc)
\]
which satisfies the following properties:
\begin{enumerate}
\item
$B_\lambda(x,y)\circ B_\lambda(y,z)=B_\lambda(x,z)$ for every
$x,y,z\in\mathring H$.
\item
$B_\lambda(x,x)=Id$ for every $x\in\mathring H$.
\item
$B_\lambda$ is constant on the strata of the stratification of $\mathring H$
determined by $\lambda$.
\item
If $\lambda_n\rightarrow\lambda$ on a $\eps$-neighbourhood of the
segment $[x,y]$ and $x,y \notin L_W$, then
$B_{\lambda_n}(x,y)\rightarrow B_{\lambda}(x,y)$ .
\end{enumerate}

If $\lambda$ is finite, then there is an easy
description of $B_{\lambda}$. If $l$ is an oriented geodesic of
$\mh^3$, let $X_l\in\sG\lG(2,\mc)$ denote the infinitesimal generator
of the positive rotation around $l$ such that $\exp(2\pi X_l)=Id$
(since $l$ is oriented the notion of \emph{positive} rotation is well
defined). Now take $x,y\in\mathring H$. If they lie in the same leaf of
$\lambda$ then put $B_\lambda(x,y)=Id$. If both $x$ and $y$ do not lie
on the support of $\lambda$, then let $l_1,\ldots,l_s$ be the
geodesics of $\lambda$ meeting the segment $[x,y]$ and $a_1,\ldots,
a_s$ be the respective weights.  Let us consider the orientation on
$l_i$ induced by the half plane bounded by $l_i$ containing $x$ and
non-containing $y$. Then put
\[
   B_\lambda(x,y)=\exp(a_1 X_1)\circ\exp(a_2 X_2)\circ\cdots\circ\exp(a_s X_s)
\, .
\]
If $x$ lies in $l_1$ use the same construction, but replace $a_1$
by $a_1/2$; if $y$ lies in $l_s$ replace $a_s$ by $a_s/2$.

\smallskip

The bending cocycle is not continuous on the whole definition set.  However,
there is a natural continuous ``pull-back'' of it to a cocycle defined on the
flat spacetime $\Uu = \Uu^\lambda_0$
\[
    \hat B_\lambda:\Uu \times\Uu \rightarrow PSL(2,\mc)
\]  
such that 
\[
   \hat B_\lambda (p,q)=B_\lambda(N(p), N(q))
\]
for $p,q$ such that $N(p)$ and $N(q)$ do not lie on $L_W$. 

This map is locally Lipschitz (with respect to the Euclidean distance
on $\Uu$).  Moreover, for every compact set $K$ of $\Uu$, the Lipschitz
constant on $K\times K$ depends only on $N(K)$, on the diameter of
$r(1,\cdot)(K)$ and on the maximum $M$ and minimum $m$ of $T$ on $K$.
\smallskip

{\bf The bending map.} Fix a base point $x_0$ of $\mathring{H}$ ($x_0$
is supposed not to be in $L_W$). The {\it bending map} of $\mathring{H}$
along $\lambda$ is
\[
    F=F_\lambda:\mathring{H}\ni x\mapsto B(x_0,x)x\in\mh^3\ .
\] 
$F$ satisfies the following properties:
\begin{enumerate}
\item
It does not depend on $x_0$ up to post-composition of elements of
$PSL(2,\mc)$.
\item
It is a $1$-Lipschitz map.
\item
If $\lambda_n\rightarrow\lambda$ then $F_{\lambda_n}\rightarrow F_\lambda$
with respect to the compact open topology.
\end{enumerate}
\smallskip

{\bf The Wick rotation.}
We are ready to construct the local $\mathrm C^1$-diffeomorphism
\[
   d_\mh:\Uu(>1)\rightarrow\mh^3
\]
with the properties outlined at the beginning of this Section.

Recall the continuous cocycle $\hat B = \hat
B_\lambda$ defined above on the whole of $\Uu \times
\Uu$. Since both $\mh^3$ and $\mh^2 \subset \mh^3$ 
are oriented, the normal bundle is oriented too. Let
$v$ denote the normal vector field on $\mh^2$ that is positive
oriented with respect to the orientation of the normal bundle.
Let us take $p_0 \in N^{-1}(x_0)$ and for
$p\in\Uu(>1)$ consider the geodesic ray $c_p$ of $\mh^3$
starting from $F(N(p))$ with speed vector equal to $w(p)=\hat
B(p_0,p)_*(v(N(p)))$.  Thus $d_\mh$ is defined in the
following way:
\[
    d_\mh(p)=c_p(\arctgh(1/T(p)))=
\exp_{F(N(p))}\left(\arctgh\left(\frac{1}{T(p)}\right)w(p)\right) \ .
\]

As usual, we make everything explicit on the local models of
$\Uu^0_0$ and of flat spacetimes associated to  finite laminations.
\smallskip

{\it Local model of the Wick rotation for finite laminations.}
Consider as above the future $\Uu_0$ of a spacelike segment in $\mx_0$
(adopting the same notations). We introduce suitable $\mathrm C^{1,1}$
coordinates on $\Uu_0$.  Denote by $l_a$ the boundary of
$\Uu_0^{-}(a)$ and by $d_a$ the intrinsic distance of $\Uu_0(a)$.  Fix
a point $z_0$ on $l_0$ and denote by $\hat z_a\in l_a$ the point such
that $N(\hat z_a)=z_0$.

For every $x\in\Uu_0(a)$ there is a unique point $\pi(x)\in l_a$ such
that $d_a(x,l_a)=d_a(x,\pi(x))$.  Then we consider coordinates $T,\zeta,
u$, where $T$ is again the cosmological time, and $\zeta, u$ are
defined in the following way
\[
\begin{array}{l}
\zeta(x)=\eps(x) d_{T(x)}(x,l_{T(x)})/T(x)\\
u(x)=\eps'(x)d_{T(x)}(\pi(x), \hat z_{T(x)})/T(x) 
\end{array}
\]
where $\eps(x)$ (resp. $\eps'(x)$ ) is $-1$ if $x\in\Uu_0^{-}$
(resp. $\pi(x)$ is on the left of $\hat z_{T(x)}$) and is $1$ otherwise.

Choose coordinates $(y_0,y_1,y_2)$ of the Minkowski space  such that
$v_0=(0,0,1)$ and $z_0=(1,0,0)$. Thus the parametrization induced by
$T,\zeta,u$ is
\[
 (T,u,\zeta)\mapsto\left\{\begin{array}{ll} T(\ch u\ch \zeta,\ \sh
                     u\ch\zeta,\ \sh\zeta) & \textrm{ if }\zeta<0\\
                     T(\ch u,\ \sh u,\ \zeta) & 
                     \textrm{ if }\zeta\in[0,\alpha_0/T]\\ 
                     T(\ch u\ch\zeta',\ \sh
                     u\ch\zeta',\ \sh\zeta'+\alpha_0/T) &\textrm{otherwise}
                     \end{array}\right.
\]
where we have put $\zeta'=\zeta-\alpha_0/T$.\par
With respect to these coordinates the metric take the following form: 
\[
 h_0(T,\zeta, u)=\left\{\begin{array}{ll}
                           -\d T^2+T^2(\d \zeta^2+\ch^2\zeta \d u^2) &
                           \textrm{ if }\zeta<0\,,\\
                           -\d T^2 + T^2(\d \zeta^2+\d u^2) & 
                           \textrm{ if }\zeta\in[0,\alpha_0/T]\,,\\
                           -d T^2+ T^2(\d \zeta^2+\ch^2(\zeta')\d u^2) &
                           \textrm{ otherwise.}
                            \end{array} \right.
\]
Notice that the gradient of $T$ is just the coordinate field
$\frac{\partial\,}{\partial T}$.

The Gauss map takes the following form
\[
   N(T,\zeta,u)=\left\{\begin{array}{ll} (\ch u\ch \zeta,\ \sh
                     u\ch\zeta,\ \sh\zeta) & \textrm{ if }\zeta<0\\
       (\ch u,\ \sh u,\ 0) &\textrm{if }\zeta\in[0,\alpha_0/T]\\ 
       (\ch u\ch\zeta',\ \sh u\ch\zeta',\ \sh\zeta') & \textrm{otherwise}
                     \end{array}\right.
\]
and the bending cocycle $\hat B_0(p_0, (T,\zeta,u))$ is the rotation
around $l_0$ of angle equal to $0$ if $\zeta<0$, $\zeta$ if
$\zeta\in[0,\alpha_0/T]$, $\alpha_0/T$ otherwise.

Let $\mh^3$ be identified with the set of timelike unit vectors in the
$3+1$-Minkowski space $\mm^4$. We can choose affine coordinates on
$\mm^4$ in such a way the inclusion $\mh^3\subset\mh^4$ is induced by
the inclusion $\mx_0\rightarrow\mm^4$ given by
$(x_0,x_1,x_2)\mapsto(x_0,x_1,x_2,0)$.  Thus the general rotation
around $l_0$ of angle $\alpha$ is represented by the linear
transformation $T_\alpha$, such that
\[
T_\alpha(e_0)=e_0,\ T_\alpha(e_1)= e_1, \ T(e_2)=\cos\alpha\
e_2+\sin\alpha\ e_3,\ T_\alpha(e_3)=-\sin\alpha\ e_2+\cos\alpha\
e_3\,
\]
where $(e_0, e_1, e_2, e_3)$ is the canonical basis of $\mr^{4}$.
Thus,  we can write in local coordinates $d_\mh= D_0$
\[
D_0(T,u, \zeta)\mapsto\left\{\begin{array}{ll}
                          \ch\delta\left(\ch\zeta\ch u,\ \ch\zeta\sh u,\
                          \sh\zeta,\ 0\right)\ +\ \sh\delta(0,0,0,1) \\
                          \textrm{if }\zeta\leq 0\ ;\\ \ch\delta\left(\ch
                          u,\ \sh u,\ 0,0\right)\ + \
                          \sh\delta\left(0,\ 0,\
                          -\sin\frac{\zeta}{\tgh\delta},\,
                          \cos\frac{\zeta}{\tgh\delta}\right) \\
                          \textrm{if }\zeta\in[0,\alpha_0/T]\ ;\\

                          \ch\delta\left(\ch\zeta'\ch u,\ \ch\zeta'\sh
                          u,\ \sh\zeta'\cos\alpha_0,\
                          \sh\zeta'\sin\alpha_0\right)\ + & \\
                          \sh\delta(0,0,-\sin\alpha_0,\
                          \cos\alpha_0)\\ \textrm{otherwise}
                           \end{array}\right.
\] 
where $\delta=\arctgh(1/T)$ and $\zeta'=\eta-\alpha_0/T$. This map is
clearly smooth for $\zeta\neq 0,\alpha_0/T$.  Since the derivatives of
$D_0$ with respect the coordinates fields glue along $\zeta=0$ and
$\zeta=\alpha_0 T$ the map $D_0$ is $\mathrm C^1$.  It is not hard to
see that the derivatives are locally Lipschitz.  One can check by
direct computation that $D_0^*(g)$ is obtained by the canonical Wick
rotation. The same formulae hold on $\Uu^0_0$, providing that we replace
$\Uu_0^+\cup \Vv$ by $r^{-1}(0,+\infty v_0)$, the inverse image
of the open ray.

\begin{figure}
\begin{center}
\input{1geo-hyp.pstex_t}
\caption{The image $\Ee_0$ of $D_0$ and its decomposition.}
\end{center}
\end{figure}
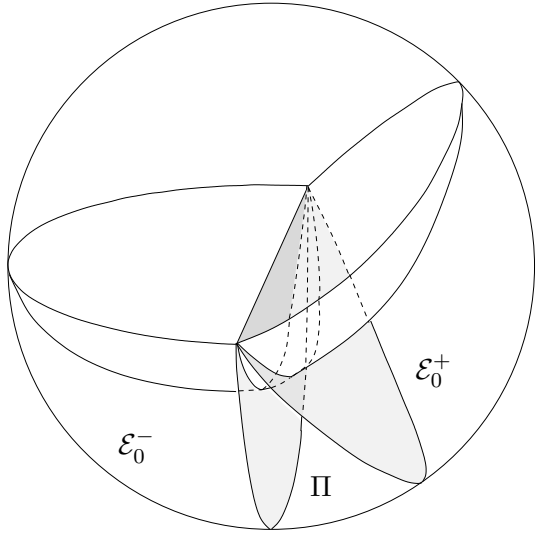

{\bf The holonomy $h_\mh$. } Recall that $(F,\lambda)\in \Mm\Ll(S)$,
$F\subset F^\Cc \subset \hat {F} = \mh^2/\Gamma $, $F= \mathring H/\Gamma$.
Then (see \cite{Ep-M}) the bending cocycle satisfies:

\[
   B_\lambda(\gamma x,\gamma y)=\gamma\circ B(x,y)\circ\gamma^{-1}
\]
for every $\gamma\in\Gamma$.

Consider a bending map 
\[
   F_\lambda:\mathring H\rightarrow\mh^3\ .
\]

For $\gamma\in\Gamma$ let us define
\[
  h_\mh(\gamma)= B_\lambda(x_0,\gamma x_0)\circ\gamma\in PSL(2,\mc)
\]

clearly $F_\lambda$ is $h_\mh$-equivariant. We eventually get that the
Wick rotation descends on the quotient spacetime $Y=Y^\lambda_0$, this
gives the required hyperbolic structure $M$ on $Y(>1)$, having as
asymptotic boundary the projective surface $S^\lambda_\Pp$.

\subsection {Flat vs de Sitter Lorentzian geometry}\label{dS}
In order to classify $MGH$ de Sitter spacetimes of finite type in
terms of complex projective structures, we widely refer to \cite{Sc}
where the case of {\it compact} Cauchy surfaces was treated. In fact we
can check that all constructions work as well by simply letting the
Cauchy surface be complete of finite type. Let us summarize the main steps
of this classification:
\smallskip

(1) We associate to every complex projective structure on a surface of
finite type $S$ a so called \emph{standard} spacetime belonging to
$\Mm\Gg\Hh_1(S)$. It turns that it is future complete.  By composing
with the parametrization $\mG_\Pp:\Mm\Ll(S)\to \Pp(S)$, we eventually
construct the {\it injective} map $\mG_1:\Mm\Ll(S) \to
\Mm\Gg\Hh_{1}(S)$.
\smallskip

(2) We show that, {\it possibly inverting the time orientation}, every
    spacetime in $\Mm\Gg\Hh_{-1}(S)$ is standard, that is 
$\mG_1:\Mm\Ll(S) \to \Mm\Gg\Hh_{1}(S)/\pm$ is a {\it bijection} (with the
same meaning of $\pm$ as for $\mG_0$).
\medskip

We recall the construction of these standard spacetimes.  Given a
projective structure on $S$, with developing
map
\[
     d:\tilde S\rightarrow S^2_\infty
\]
we perform a construction which is dual to the one made for the
$H$-hulls. Recall the canonical stratification of $\tilde S$.  
For every $p\in\tilde S$ let $U(p)$ denote the stratum
passing through $p$ and $U^*(p)$ be the maximal ball containing
$U(p)$.  Now $d(U^*(p))$ is a ball in $S^2_\infty$ which determines a
hyperbolic plane in $\mh^3$. Let $\rho(p)$ denote the point in $\mx_1$
corresponding to this plane: the map $\rho:\tilde S\rightarrow\mx_1$
turns out to be continuous. There exists a unique timelike geodesic
$c_p$ in $\mx_1$ joining $\rho(p)$ to $d(p)$ so we can define the map
\[
    \hat d:\Delta\times(0,+\infty)\ni (p,t)\mapsto c_p(t)\in\mx_1
\]
This map is a developing map for the required standard de Sitter
spacetime. A compatible holonomy follows by a natural equivariant version
of the construction. 
\smallskip

Assume now that the the projective structure is encoded by
$(F,\lambda)\in \Mm\Ll(S)$, via $\mG_\Pp$.  We eventually realize that
the construction of $\mG_1$ can be obtained via a canonical rescaling
performed on $Y^0_\lambda(<1)$.  More precisely we realize $\hat d$ as
a C$^1$ developing map
\[
   d^\lambda_1:\Uu_\lambda(<1)\rightarrow\mx_1
\]
obtained as a sort of semi-analytic continuation of the hyperbolic
developing map $d_\mh$ constructed in the previous Section, and we have:

\begin{teo}\label{memgen:dsmain:teo}
The spacetime $\Uu^1_\lambda$ ($Y^1_\lambda$), obtained from
$\Uu^0_\lambda(<1)$ ($Y^0_\lambda(<1)$) via the rescaling directed by
the gradient of its cosmological time $T$ and with rescaling functions
\[
   \alpha =\frac{1}{1-T^2} \qquad \beta= \frac{1}{(1-T^2)^2}\ .
\]
is the standard \emph{de Sitter} spacetime
corresponding to the projective structure on $\Uu^0_\lambda(1)$
($Y^0_\lambda(1)$) produced by the Wick rotation.
\end{teo}

The construction of $d^\lambda_1$ is very simple. We regard both $\mh^3$ and
$\mx_1$ as open sets of the real projective space (Klein models),
separated by the quadric $S^2_\infty$. If $s$ is a geodesic
integral line of the gradient of the cosmological time,
$s_{>1}=s\cap\Uu^0_\lambda(>1)$ is sent by $d_\mh$  onto a geodesic ray of
$\mh^3$. We define $d^\lambda_1$ on $s_{<1}$ in such a way that it
parameterizes the timelike geodesic ray in $\mx_1$ contained in the
projective line (in the Klein model) determined by $d_\mh(s_{>1})$.
\\

\begin{figure}
\begin{center}
\input{1geo-ds2.pstex_t}
\caption{A standard de Sitter spacetime - local model}
\end{center}
\end{figure}
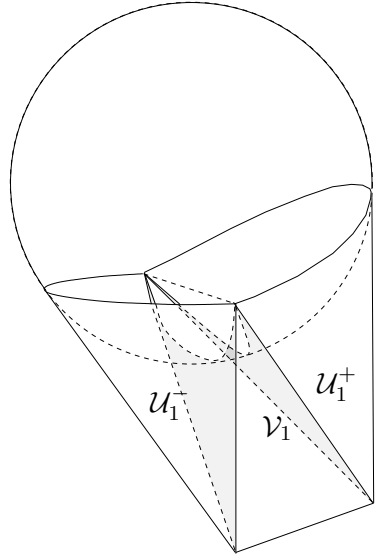

The proof, as well as the explicit computation for our favourite local
models (evoked in the Figure) are similar to the ones made for the Wick
rotation, so we omit them. 

A main point in proving that $\mG_1$ is a bijection consists in

\begin{prop}\label{CT-1}
(1) Every $Y\in \Mm\Gg\Hh_1(S)$ has C$^1$ cosmological time, and every
level surface is a complete Cauchy surface. 

(2) If $Y=Y^1_\lambda$ with universal covering $\Uu^1_\lambda$, then
    the cosmological time of $\Uu^1_\lambda$ is
\[
   \tau=\arctgh(T).
\]
$T$ being the cosmological time of $\Uu^0_\lambda(<1)$. Hence the
inverse rescaling is directed by the gradient of $\tau$ and has
universal rescaling functions.

(3) Let $\Sigma^0$ be the initial singularity of $\Uu^0_\lambda$. Then
the map $d^\lambda_1$ extends to a continuous map
\[
   \Uu^0_\lambda(\leq 1)\cup\Sigma^0\rightarrow \mx_1\cup S^2_\infty \ .
\]
Moreover, its restriction to $\Uu^0_\lambda(1)$ coincides with $d_\mh$;
the restriction to $\Sigma^0$ is an (equivariant) isometry onto
the initial singularity $\Sigma^1$ of $\Uu^1_\lambda$.
\end{prop}
Note that, in contrast with the flat Lorentzian case, these de Sitter
developing maps, as well as the dual hyperbolic ones are in general
{\it not} injective.

\subsection {Flat vs Anti de Sitter Lorentzian geometry}
\label{AdS}
We are going to outline first a few features of the spacetimes in
$\Mm\Gg\Hh_{-1}(S)$. Recall the content of Section \ref{2+1ST}; in
particular the duality between points $x$ of $\mx_{-1}$ and spacelike
planes $P(x)$, or between spacelike lines, $l\to l^*$. Recall also
that the boundary $\partial \mx_{-1}$ has a natural causal structure,
so that the notion of a {\it nowhere timelike simple closed curve}
embedded in $\partial\mx_{-1}$ makes sense.

{\bf Standard AdS spacetimes.}
Given such a curve $C\subset \mx_{-1}$, assume furthermore that 
 $C$ is different from a (left or right) leaf of the natural
 double foliation of $\partial\mx_{-1}$ (that it is a so called {\it
 admissible achronal curve}). Then its \emph{Cauchy development}
is defined as
\[
   \Yy(C)=\{p\in \mx_{-1}| \partial P(p)\cap C=\varnothing\}
\]
and the so obtained spacetime is called a (simply
 connected) {\it standard AdS spacetime}. $C$ is said its {\it curve at
 infinity}. $C$ is homotopic to the meridian
 of $\partial\mx_{-1}$ with respect to $\mx_{-1}$. In general and AdS
spacetime is said standard if its universal covering is standard.
\smallskip

{\bf The convex core.}  There exists a spacelike plane $P$ not
intersecting $\Yy(C)$ (see \cite{M}). In the Klein model we can cut
$\Pm^3$ along the projective plane $\hat P$ containing $P$ and we have
that $\Yy(C)$ is contained in $\mr^3=\Pm^3\setminus\hat P$. Since $C$
is nowhere timelike, then for every point $p\in C$ the plane $P(p)$
tangent to $\partial\mx_{-1}$ at $p$ (that cuts $\mx_{-1}$ at a null
totally geodesic plane) does not separate $C$. It follows that the
convex hull $\Kk(C)$ of $C$ in $\mr^3$ is actually contained in
$\mx_{-1}$.  We realize that $\Kk(C)$ does not depend on the choice of
$\hat P$, and is called the \emph{convex core} of $\Yy(C)$.

Support planes of $\Kk(C)$ are non-timelike and the closure
$\overline{\Yy(C)}$ of $\Yy(C)$ in $\mx_{-1}$ coincides with the set
of dual points of spacelike support planes of $\Kk(C)$ whereas the set
of points dual to null support planes of $\Kk(C)$ coincides with $C$.
$\overline{\Yy(C)}$ is convex and the closure of $\Yy(C)$ in
$\overline\mx_{-1}$ is $\overline{\Yy(C)}\cup C$.  It follows that
$\Kk(C)\subset\overline{\Yy(C)}$.  A point $p\in\partial\Kk(C)$ lies
in $\Yy(C)$ if and only if it is touched only by spacelike support
planes.
 
Being the boundary of a convex set in $\mr^3$, $\partial\Kk(C)\cup C$
is homeomorphic to a sphere. In particular $\partial\Kk(C)$ (that is
the boundary of $\Kk(C)$ in $\mx_{-1}$) is obtained by removing a
circle from a sphere, so it is the union of two disks. These
components will be called {\it the past and the future boundary} of
$\Kk(C)$ (with respect to the time orientation), and denoted
$\partial_-\Kk(C)$ and $\partial_+\Kk(C)$ respectively.  Given any
inextendible timelike ray contained in $\Kk(C)$, its future end-point
lies on the future boundary, and the past end-point lies on the past
boundary.

$\partial_+\Kk(C)\cap\Yy(C)$ is obtained by removing from
$\partial_+\Kk(C)$ the set of points that admits a null support plane.
Now suppose that a null support plane $P$ passes through
$x\in\partial_+\Kk(C)$. Then $P\cap\Kk(C)$ is a triangle with a
vertex at $x(P)$, two ideal edges (that are segments on the leaves
of the double foliation of $\partial\mx_{-1}$) and a complete geodesic
of $\Kk(C)$. It follows that the set $\partial_+\Kk(C)\cap\Yy(C)$ is
obtained by removing from $\partial_+\Kk(C)$ (at most) numerable many
ideal triangles, so it is homeomorphic to a disk.  The only case for
$\partial_+\Kk(C)\cap\Yy(C)$ to be empty is that the curve $C$ is
obtained by joining the end-points of a spacelike geodesic $l$ with
the end-points of its dual geodesic $l^*$; in that case
$\Yy(C)=\Kk(C)$, and we call it the {\it degenerate standard
spacetime}.
So, from now on, we incorporate in the definition of \emph{standard
AdS spacetime} that it is not degenerate.  Moreover, since we will be
mainly interested in $\partial_+\Kk(C)\cap\Yy(C)$, from now on we will
use $\partial_+\Kk(C)$ just to denote that set.

\begin{prop}\label{ads:boundary:prop} 
$\partial_+\Kk(C)$ is locally $\mathrm C^0$-isometric to $\mh^2$. 
\end{prop}
\begin{remark}\label{ads:bend:inj:rem:3}
{\rm If $\partial_+ \Kk$ is complete then it is isometric to $\mh^2$.
In general $\partial_+ \Kk$ is not complete, not even in the special
case when $C$ is the graph of a homeomorphism of $S^1$ onto itself. Moreover,
it can be not complete even when there are no null triangles on the boundary.}
\end{remark}

{\bf The past part of a standard spacetime.}
The \emph{ past part } $\Pp = \Pp(C)$ of a standard AdS spacetime
$\Yy(C)$ is the past in $\Yy(C)$ of the future boundary $\partial_+
\Kk$ of its convex core.  The complement of $\partial_+ \Kk$ in the
frontier of $\Pp(C)$ in $\mx_{-1}$ is called the {\it past boundary}
of $\Yy(C)$, denoted by $\partial_- \Pp$.

\begin{prop}\label{ads:ct:prop}
Let $\Pp$ be the past part of some $\Yy(C)$.  Then $\Pp$ has
cosmological time $\tau$ and this takes values on $(0,\pi/2)$. For
every point $p\in \Pp$ there exist only one point $\rho_-(p)\in
\partial_- \Pp$, and only one point $\rho_+(p)\in \partial_+ \Kk$ such
that
\medskip\par\noindent
\emph{1. }$p$ is on the timelike segment joining $\rho_-(p)$ to $\rho_+(p)$.
\medskip\par\noindent
\emph{2. }$\tau(p)$ is equal to the length of the segment $[\rho_-(p),p]$.
\medskip\par\noindent
\emph{3. }the length of $[\rho_-(p),\rho_+(p)]$ is $\pi/2$. 
\medskip\par\noindent 
\emph{4. }$P(\rho_-(p))$ is a support plane for
$\Pp$ passing through $\rho_+(p)$ and $P(\rho_+(p))$ is a support
plane for $\Pp$ passing through $\rho_-(p)$.
\medskip\par\noindent \emph{5. }The map $p\mapsto \rho_-(p)$ is
continuous. The function $\tau$ is $\mathrm C^1$ and its gradient at
$p$ is the unit timelike tangent vector $\mathrm{grad}\,\tau(p)$ such
that
\[
    \exp_p\left(\tau(p)\mathrm{grad}\,\tau(p)\right)=\rho_-(p)\,.
\]
\end{prop}

Summing up, given the past part $\Pp$ of a standard AdS spacetime
$\Yy(C)$, we can construct:
\smallskip\par\noindent
the cosmological time $\tau:\Pp\rightarrow (0,\pi/2)$;
\smallskip\par\noindent
the future retraction  $\rho_+:\Pp\rightarrow\partial_+\Kk$;
\smallskip\par\noindent
the past retraction $\rho_-:\Pp\rightarrow\partial_-\Pp$.

\begin{cor}\label{ads:ct:cor}
\emph{1.} Given $r$ in the past boundary of $\Yy$, $\rho_-^{-1}(r)$ is
the set of points $p$ such that the ray starting from $r$ towards $p$
meets at time $\pi/2$ the future boundary of $\Kk$.
\medskip\par\noindent \emph{2.} The image of $\rho_-$ is the set of
points of $\partial_-\Pp$ whose dual plane meets $C$ at least in two
points.
\medskip\par\noindent
\emph{3.} The image of $\rho_+$ is the whole $\partial_+\Kk$.
\end{cor}

The image of the past retraction  is called the
\emph{initial singularity} of $\Yy(C)$ .

For every surface of finite type $S$, Stand$_{-1}(S)$ denotes the 
Teichm\"uller-like space of standard AdS spacetimes admitting a Cauchy
surface homeomorphic to $S$. The following is a fundamental step towards
the classification.
\begin{teo}\label{allstandard} {\rm Stand}$_{-1}(S) = \Mm\Gg\Hh_{-1}(S)$.
\end{teo}
Note that a consequence of this theorem is that, similarly to the
flat case, the developing maps of finite type $MGH$ AdS spacetimes
are {\it embedding} onto convex domains.
 
The fact that every spacetime in $\Mm\Gg\Hh_{-1}(S)$ is standard
follows from the following more general result (Section 7
of~\cite{M}).

\begin{prop}\label{memgen:adsstand:prop}
Let $Y$ be an Anti de Sitter {\rm simply connected} spacetime, and
$F\subset Y$ be a complete Cauchy surface. Then the developing map $
Y\rightarrow\mx_{-1}$ is an embedding onto a convex subset of
$\mx_{-1}$.

The closure of $F$ in $\overline{\mx}_0$ is a closed
disk and its boundary $\partial F$ is a nowhere timelike curve of
$\partial\mx_{-1}$.

If $Y$ is the maximal globally hyperbolic Anti de Sitter spacetime
containing $F$ then $Y=\Yy(C)$.  The curve $\partial F$ determines
$Y$, namely $p\in Y$ iff the dual plane $P(p)$ does not meet $\partial
F$.

Conversely $\partial F$ is determined by $Y$, in fact
$\partial F$ is the set of accumulation points of $Y$ on
$\partial\mx_{-1}$.  If $F'$ is another complete spacelike Cauchy
surface of $Y$ then $\partial F'=\partial F$.
\end{prop}

The main step in order to prove the opposite inclusion is the
following proposition (recently achieved also by Barbot~\cite{Ba}(2)
with a different approach with respect to \cite{Be-Bo}), that also
holds for arbitrary standard spacetimes.

\begin{prop}\label{memgen:adscompl:prop}
If $\Pp$ is the past part of $\Yy(C)$ then every level surface
$\Pp(a)$ of the cosmological time is complete.
\end{prop}

\begin{cor}\label{memgen:adscompl:cor}
 Every level surface $\Pp(a)$ of the past part $\Pp$ of a standard AdS
spacetime $\Yy(C)$ is a complete Cauchy surface of $\Yy(C)$ and this
last is the maximal globally hyperbolic AdS spacetime that extends
$\Pp$.
\end{cor}
\begin{remark}{\rm 
$\tau$ extends to the cosmological time of $\Yy(C)$, that takes
values on some interval $(0, a_0(C))$, for some well defined $\pi/2 <
a_0(C) < \pi$. Notice however that $\tau$ is C$^1$ only on the past part.}
\end{remark}

{\bf The map $\mG_{-1}$.}  Let $(F,\lambda)\in \Mm\Ll(S)$, $F\subset
F^\Cc \subset \hat {F} = \mh^2/\Gamma $ $(F^\Cc,\lambda)$, with
universal coverings $\mh^2 \to \hat {F}$, $H \to F^\Cc$, and
$\mathring {H}\to F$ respectively, as usual.  Fix an embedding of
$\mh^2$ in $\mx_{-1}$ as a spacelike plane (for instance as $P(Id)$).
The key ingredient to construct $\mG_{-1}$ is the AdS version of the
{\it bending of $H\subset \mh^2$ along the lamination $\lambda$} (see
below).  This produces a {\it convex embedding}
$\varphi_\lambda:\mathring H\rightarrow\mx_{-1}$. Recall that to
construct the $H$-hull (via the Wick rotation) we used the bending map
$f_\lambda: \mathring H\rightarrow\mh^3$, that is a local convex
embedding, and then we followed the geodesic rays normal to
$f_\lambda(\mathring H)$, in the {\it non-convex} side bounded by
$f_\lambda(\mathring H)$. Eventually the developing map $d_\mh$ has
been obtained by requiring that the integral lines of the cosmological
times would be sent to the integral lines of the normal flow. Also in
the present situation we construct a C$^1$ developing map
$$d^\lambda_{-1}: \Uu^\lambda_0 \to \mx_{-1}$$ by requiring that the
integral lines of the cosmological time of $\Uu^0_\lambda$ are sent to
the integral line of the normal flow. An important difference, with
respect to the hyperbolic case, is that the normal flow is followed now
in the {\it convex} side bounded by $\varphi_\lambda(\mathring H)$
(otherwise singularities would be reached). It turns that the image of
$d^\lambda_{-1}$ is the past part of a standard AdS spacetime, that plays
here the role of a sort of AdS-hull. More precisely we have:
\begin{teo} (1) $d^\lambda_{-1}$ is an embedding onto the past part
$\Pp^\lambda$ of a determined $\Uu^\lambda_{-1}= \Yy(C^\lambda)$,
which is the universal covering of $Y^\lambda_{-1}\in
\Mm\Gg\Hh_{-1}(S)$.  The image of the AdS bending map
$\varphi_\lambda$ coincides with $\partial_+\Kk(C^\lambda)$. The map
$d^\lambda_{-1}$ continuously extends to an isometry between the
respective initial singularity.

(2) $\Uu^\lambda_{-1}$ is produced by the rescaling of
$\Uu^\lambda_0$, directed by the gradient of the cosmological time
$T$, with universal rescaling functions
\[
\begin{array}{ll}   
\alpha =\frac{1}{1+T^2} \ , \qquad &   
\beta=\frac{1}{(1+T^2)^2} \ .
\end{array}
\]

(3) The cosmological time $\tau$ on $\Pp^\lambda$ is given by
 \[
   T(p)=\tan\tau(p)\ .
\]
\end{teo}

In such a way we construct an {\it injective} map
$$\mG_{-1}: \Mm\Ll(S) \to \Mm\Gg\Hh_{-1}(S) \ .$$
The following general proposition (specialized to  Stand$_{-1}(S)$) 
implies that $\mG_{-1}$ is in fact a {\it bijection}.
\begin{prop}\label{ads:flat:teo} For every standard AdS domain $\Yy=\Yy(C)$,
the rescaling of its past part $\Pp$, directed by the gradient of the
cosmological time $\tau$, with universal rescaling functions
\[
   \alpha=\frac{1}{\cos^2\tau} \qquad \beta=\frac{1}{\cos^4\tau}
\]
produces a regular domain, whose cosmological time is given by the formula
\[
    T=\tan\tau\ .
\] 
\end{prop} 
It is not too hard to see, by means of local considerations, that such a
rescaling produces a flat spacetime.  Showing that it is regular
domain is actually more demanding. This is equivalent to show that the
future boundary of the convex core is isometric to a straight convex
set pleated at a measured lamination. The key point is the fact that
level surfaces of $\tau$ are complete (Proposition
\ref{memgen:adscompl:cor}).
\smallskip

{\bf On AdS bending.} We are going to outline more precisely the
construction of $d^\lambda_{-1}$.  The AdS bending runs similarly to
the hyperbolic one, but having some remarkable differences (that are
eventually responsible, for example, that the AdS developing maps are
embeddings, in contrast with the hyperbolic ones). We also stress that
orientations play a subtle role in the AdS bending procedure. The
basic diffence arises from the different behaviour of the ``angles''
between hyperbolic planes (that is spacelike planes) and of
``rotations'' around spacelike geodesics in $\mx_{-1}$, with respect
to $\mh^3$. In fact, given two spacelike planes $P_1,\ P_2$ meeting
each other along a geodesic $l$, the dual points $x_i=x(P_i)$ lie on
the geodesic $l^*$ dual to $l$.  Then we define the {\it angle between
$P_1$ and $P_2$} as the distance between $x_1$ and $x_2$ along $l^*$.
Notice that:
\smallskip

{\it Fix $P_1$, the by varying $P_2$, the angles between them are well
defined numbers that span the whole of the interval $(0,+\infty)$}.
\smallskip 

Define a {\it rotation around a spacelike geodesic
$l$}  simply to be an isometry of $\mx_{-1}$ which
point-wise fixes $l$. We have
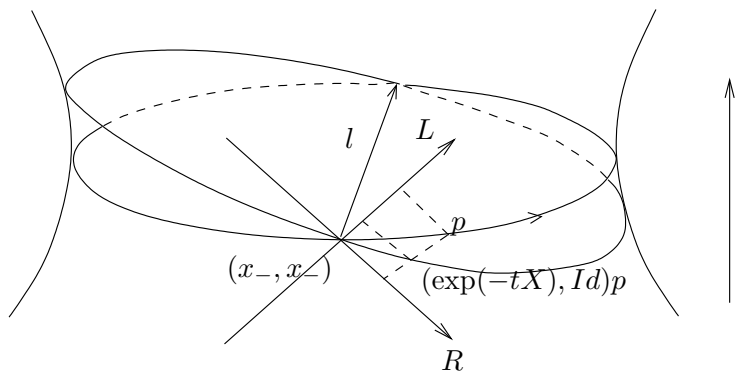
\begin{figure}
\begin{center}
\input{positive-rotation.pstex_t}
\caption{$(\exp(-tX),Id)$ rotates planes around $l$ in the
positive sense.}\label{ADS:rot:fig}
\end{center}
\end{figure}

\begin{lem}
Rotations around a geodesic $l$ act freely and transitively on the
dual geodesic $l^*$. Such action induces an isomorphism between the
set of rotations around $l$ and the set of translations of $l^*$. 

By duality, rotations around $l$ act freely
and transitively on the set of spacelike planes containing $l$.  
Given two spacelike planes $P_1,\ P_2$ such that $l\subset P_i$, then
there exists a unique rotation $T_{1,2}$ around $l$ such that
$T_{1,2}(P_1)=P_2$.
\end{lem}  

\begin{lem} 
An isometry of $\mx_{-1}$ is a rotation around a geodesic if and
only if it is represented by a pair $(x,y)$ such that $x$ and $y$ are
isometries of $\mh^2$ of hyperbolic type with the same translation length. 

Given two spacelike planes $P_1,P_2$ meeting along a geodesic $l$,
let $(x,y)$ be the rotation taking $P_1$ to $P_2$. Then the  
translation length $\tau$ of $x$ coincides with the angle between $P_1$ 
and $P_2$.
\end{lem}

There is a natural definition of positive rotation around an oriented
spacelike geodesic $l$ (depending only on the orientations of $l$ and
$\mx_{-1}$).  Thus, an orientation on the dual line $l^*$ is induced
by requiring that positive rotations act by positive translations on
$l^*$. In particular, if we take an oriented geodesic
$l$ in $P(Id)$, and denote by $X$ the infinitesimal generator
of positive translations along $l$ then it is not difficult to show that the
positive rotations around $l$ are of the form $(\exp(-tX),\exp(tX))$ for
$t>0$. Actually, by looking at the action on the boundary we  deduce
that both the maps $(\exp(-tX),Id)$ and
$(Id,\exp(tX))$ rotate planes through $l$ in the positive direction
(see Fig.~\ref{ADS:rot:fig}).
\smallskip

Given $\lambda$ on $H$ as usual, we construct now an {\it AdS bending
cocycle}

\[
    B^\lambda = (B^\lambda_-,B^\lambda_+):
\mathring H \times\mathring H \rightarrow  PSL(2,\R)\times PSL(2,\R)
\]
 
which formally satisfies the similar properties like the quake
cocycles of Section \ref{bend-quake}, or the above hyperbolic bending cocycle.
In fact $B^\lambda_-$ and $B^\lambda_+$ are exactly the Epstein-Marden
cocycles (like the quake cocycles), corresponding to the {\it
real-valued} measured laminations $-\lambda$ and $\lambda$.  Here
$-\lambda = (L,-\mu)$, that is we take the {\it negative-valued}
measure $-\mu$. Although this is no longer a measured lamination in
the ordinary sense , the construction of \cite{Ep-M} does apply.
Besides the usual cocycle properties, $B^\lambda$ also verifies that
if $x,y$ lie in different strata then $B^\lambda_+(x,y)$
(resp. $B^\lambda_-(x,y)$) is a non-trivial hyperbolic transformation
whose axis separates the stratum through $x$ and the stratum through
$y$. Moreover the translation length is bigger than the total mass of
$[x,y]$.

All this is very simple on the usual local model for finite laminations.
In fact, take  take a finite measured geodesic lamination $\lambda$ of
$\mh^2$.  Take a pair of points $x,y \in\mh^2$ and enumerate the
geodesics in $\lambda$ that cut the segment $[x,y]$ in the natural way
$l_1,\ldots,l_n$.  Moreover, we can orient $l_i$ as the boundary of
the half-plane containing $x$.  With a little abuse, denote by $l_i$ 
also the geodesic in $P(Id)$
corresponding to $l_i$, then let $B^\lambda (x,y)$ be the isometry of
$\mx_{-1}$ obtained by composition of positive rotations around $l_i$
of angle $a_i$ equal to the weight of $l_i$.  In particular, if $X_i$
denotes the unit positive generator of the hyperbolic transformations with
axis equal to $l_i$, then we have
\[
\begin{array}{l}
B_\lambda(x,y)=
(B^\lambda_-(x,y),B^\lambda_+(x,y))\in PSL(2,\R)\times PSL(2,\R)\qquad\textrm
{where}\\
B^\lambda_-(x,y)=\exp(-a_1 X_1/2)\circ\exp(-a_2
X_2/2)\circ\ldots\circ \exp(-a_n X_n/2)\\ 
B^\lambda_+(x,y)=\exp(a_1
X_1/2)\circ\exp(a_2 X_2/2) \circ\ldots\circ\exp(a_n
X_n/2)
\end{array}
\] 
with the following possible modifications: $a_1$
is replaced by $a_1/2$ when $x$ lies on $l_1$ and $a_n$ is replaced by
$a_n/2$ when $y$ lies on $l_n$ The factor $1/2$ in the definition of
$\beta_\pm$ arises because the  translation length of $\exp tX$ is
$2t$. 
\smallskip

By means of the bending cocycle we construct a {\it AdS bending map}:
take a base point $x_0$ in $\mathring H$ and set
\[\varphi_\lambda: \mathring H\ni x\mapsto B^\lambda(x_0,x)x.\]
 
\begin{prop}\label{ads:bend:descr:prop}
The bending map $\varphi_\lambda$ is an isometric $\mathrm C^0$
embedding of $\mathring H$ onto an achronal set $\mx_{-1}$.
\end{prop}

Let $\Uu=\Uu^0_\lambda$ be the flat spacetime 
encoded by $\lambda$.
Just as in the hyperbolic case we want to ``pull-back'' the bending cocycle
$B^\lambda$ to a {\it continuous} bending cocycle
\[   \hat{B}^\lambda:\Uu \times\Uu \rightarrow PSL(2,\R)\times 
PSL(2,\R) \ .\] In fact we get a natural extension such that:

(1) For every $p,q\in \Uu$ such that $N(p)$ and $N(q)$ do not lie
on weighted part of the lamination, then
\[    \hat{B}^\lambda(p,q)=B^\lambda(N(p),N(q))\ . \]
(2) $\hat{B}^\lambda$ on the whole of $\Uu$ is constant along the
     integral geodesics of the gradient of the cosmological time $T$.

(3) It is locally Lipschitz (with respect to the Euclidean distance
on $\Uu$), and the the Lipschitz constant on $K\times K$ ($K$ being
any compact set in $\Uu$) depends only on the image of the Gauss map
$N(K)$, the maximum of the total masses of geodesic paths of $H$
joining points in $N(K)$, and the maximum and the minimum of the
cosmological time $T$ on $K$.
\smallskip

Finally we can define our developing map $d^\lambda_{-1}:
\Uu^\lambda_0 \to \mx_{-1}$.  For every $p\in\Uu^0_\lambda$, we define
$x_-(p)$ as the dual point of the plane
$\hat{B}^\lambda(p_0,p)(P(Id))$ (that is $\hat B^\lambda(p_0,p)(Id)$,
and $x_+(p)=\hat{B}^\lambda(p_0,p)(N(p))$.  Take representatives $\hat
x_-(p)$ and $\hat x_+(p)$ in $\SL{2}{R}$ such that the geodesic
segment between $\hat x_-(p)$ and $\hat x_+(p)$, is future directed.
Finally set
\[
    d^\lambda_{-1}(p)=[\cos\tau(p)\hat x_-(p)+\sin\tau(p)\hat x_+(p)]
\]
where $\tau(p)= \arctan T(p)$ .
\smallskip

As usual, we end with a few explicit computations for our favourite
local model, that is when  when $\lambda$ is a single weighted
geodesic. 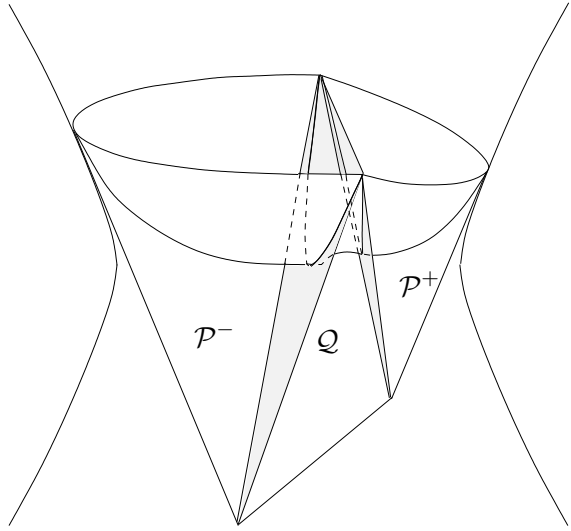
\begin{figure}
\begin{center}
\input{1geo-Ads.pstex_t}
\caption{The domain $\Pp$ with its decomposition. Also the surface
    $\Pp(a)$ is shown.}
\end{center}
\end{figure}

Let us set $\lambda_0=(l_0,a_0)$ and choose a base point
$p_0\in\mh^2-l_0$.  The surface $P=\varphi_\lambda(\mh^2)$ is simply
the union of two half-planes $P_-$and $P_+$ meeting each other along a
geodesic (that, with a little abuse of notation, is denoted by $l_0$).
We can suppose that $p_0$ is in $P_-$, and $l_0$ is oriented as the
boundary of $P_-$.  If $v_\pm$ denote the dual points of the planes
containing $P_\pm$ we have $v_-=Id$ and $v_+=\exp - a_0X_0$, $X_0$
being the standard generator of translations along $l_0$. The vector
$X_0$ is tangent to $P(id)$ along $l_0$, orthogonal to it, and points
towards $p_0$.

By definition , the image, say $\Pp$, of $\Delta_0=d^{\lambda_0}_{-1}$
is the union of three pieces: the cone with vertex at $v_-$ and basis
$P_-$, say $\Pp_-$, the cone with vertex at $v_+$ and basis $\Pp_+$,
and the join of the geodesic $l_0$ and the segment $[v_-,v_+]$, say
$\Qq$.

Fix a point in $l_0$, say $p_0$, and denote by $v_0$ the unit tangent
vector of $l_0$ at $p_0$ (that we will identify with a matrix in
$M(2,\mr)$).  Consider the coordinates on $\Uu_0$, say $(T,u,\zeta)$
introduced in Section~\ref{theWR}. With respect to these coordinates
we have
\[
 \Delta_0(T,u,\zeta)=\left\{\begin{array}{ll} \sin\tau\big
  (\ch\zeta(\ch u\ \hat p_0 + \sh u\ v_0) - \sh\zeta\ X_0 \big)\ +\
  \cos\tau\ \hat v_- & \textrm{ if } \zeta<0\\ \sin\tau( \ch u\ \hat
  p_0+\sh u\ v_0)\ +\ \cos\tau\exp(-\zeta \tan\tau\ X_0) & \textrm{ if
  } \zeta\in[0, a_0/T]\\ \sin\tau \big(\ch \zeta'(\ch u\ \hat p_0 +
  \sh u\ v_0) - \sh\zeta' X_0 \big)\ +\ \cos\tau\ \hat v_+& \textrm{
  otherwise}
\end{array}\right.
\]
where $\zeta'=\zeta-a_0/T$, $\tau=\arctan T$ and $\hat p_0, \hat v_+,
\hat v_-\in SL(2,\mr)$ are chosen as above.

Clearly $\Delta_0$ is $\mathrm C^{\infty}$ for $\zeta\neq 0, a_0/T$. A
direct computation shows that the derivatives along the coordinate
fields glue on $\zeta=0$ and $\zeta= a_0/T$ and this proves that
$\Delta_0$ is $\mathrm C^1$.

By a direct computation we have
\[
 \Delta_0^*(\eta)=\left\{\begin{array}{ll}
  -\d\tau^2+\sin^2\tau(\d\zeta^2+\ch^2\zeta\d u^2)
  %-\frac{1}{(1+T^2)^2}\d T^2+ \frac{T^2}{1+T^2}(\d\zeta^2+\ch^2\zeta\d u^2)
 & {\rm if}\ \zeta<0\\
  -\d\tau^2+\sin^2\tau(\d\zeta^2+\d u^2) & {\rm if}\ \zeta\in[0, a_0/T]\\
 -\d\tau^2+\sin^2\tau(\d\zeta^2+\ch^2\zeta'\d u^2)
 & {\rm otherwise.}
\end{array}\right.
\]

Since $\d\tau^2=\frac{1}{(1+T^2)^2}$ and $\sin^2
\tau=\frac{T^2}{1+T^2}$, we finally see that $\Delta_0$ is obtained by
a rescaling directed by the gradient of $T$ with the right rescaling
functions.
\smallskip

{\bf Compatible holonomy.}  The holonomy representation of
$Y^\lambda_{-1}\in \Mm\Gg\Ll(S)$, $h^\lambda_{-1}: \pi_1(S)\to
PSL(2,\R)\times PSL(2,\R)$, compatible with $d^\lambda_{-1}$ is
obtained as follows. If $x_0\in\mathring H$ is the usual fixed base
point of the construction, then for every $\gamma\in \pi_1(S)=\Gamma$
\[
    h^\lambda_{-1}(\gamma)= B^\lambda(x_0,\gamma x_0)\circ(\gamma,\gamma)\,.
\]

\begin{remark}\label{AdS_hyp}
{\rm It follows from the previous discussion, that the spacetimes in
$\Mm\Gg\Hh_{-1}(S)$ have a few analogies with the hyperbolic
3-manifolds arising as $H$-hulls of {\it quasi-Fuchsian} projective
surfaces belonging to $\Pp(S)$. For instance the curves at infinity
$C\subset \partial \mx_{-1}$ of $\Yy(C)$ play a similar r\^ole of the
Jordan curves that bound the universal coverings embedded in $\partial
\mh^3S^2_\infty$ of quasi-Fuchsian surfaces. However, there are
important difference that make the AdS behaviour much more ``tame''.
For example such Jordan curves are in general rather wild, while the
curves $C$ are Lipschitz. Moreover, taking for example
$S$ compact, for every $(F,\lambda)\in \Mm\Ll(S)$, along the ray
$(F,t\lambda)$ there is a critical value $t_0>0$ such that
$S^{t\lambda}\in \Pp(S)$ is quasi-Fuchsian only for $t<t_0$. On the
other hand, the description of $Y^{t\lambda}$ is qualitatively the
same for {\it every} $t>0$; in particular all AdS developing maps are
embeddings.
}
\end{remark}

%%% Local Variables: 
%%% mode: latex
%%% TeX-master: "HAND"
%%% End: 

%% file: 1geo-flat.pstex_t
\begin{picture}(0,0)%
\includegraphics{1geo-flat.pstex}%
\end{picture}%
\setlength{\unitlength}{2901sp}%
\begingroup\makeatletter\ifx\SetFigFont\undefined%
\gdef\SetFigFont#1#2#3#4#5{%
  \reset@font\fontsize{#1}{#2pt}%
  \fontfamily{#3}\fontseries{#4}\fontshape{#5}%
  \selectfont}%
\fi\endgroup%
\begin{picture}(6996,4522)(-7,-4990)
\put(3112,-1451){\makebox(0,0)[lb]{\smash{\SetFigFont{11}{12.0}{\rmdefault}{\mddefault}{\updefault}{$\Vv_0$}%
}}}
\put(751,-1501){\makebox(0,0)[lb]{\smash{\SetFigFont{11}{12.0}{\rmdefault}{\mddefault}{\updefault}{
$\Uu_0^-$}%
}}}
\put(5251,-1501){\makebox(0,0)[lb]{\smash{\SetFigFont{11}{12.0}{\rmdefault}{\mddefault}{\updefault}{$\Uu_0^+$}%
}}}
\put(2336,-4941){\makebox(0,0)[lb]{\smash{\SetFigFont{11}{12.0}{\rmdefault}{\mddefault}{\updefault}{$v_0$}%
}}}
\end{picture}

%% file: 1geo-hyp.pstex_t
\begin{picture}(0,0)%
\includegraphics{1geo-hyp.pstex}%
\end{picture}%
\setlength{\unitlength}{2072sp}%
\begingroup\makeatletter\ifx\SetFigFont\undefined%
\gdef\SetFigFont#1#2#3#4#5{%
  \reset@font\fontsize{#1}{#2pt}%
  \fontfamily{#3}\fontseries{#4}\fontshape{#5}%
  \selectfont}%
\fi\endgroup%
\begin{picture}(6321,6319)(1788,-6838)
\put(3121,-5926){\makebox(0,0)[lb]{\smash{\SetFigFont{11}{12.0}{\rmdefault}{\mddefault}{\updefault}{$\Ee^-_0$}%
}}}
\put(5416,-6421){\makebox(0,0)[lb]{\smash{\SetFigFont{11}{12.0}{\rmdefault}{\mddefault}{\updefault}{$\Pi$}%
}}}
\put(6676,-5026){\makebox(0,0)[lb]{\smash{\SetFigFont{11}{12.0}{\rmdefault}{\mddefault}{\updefault}{$\Ee_0^+$}%
}}}
\end{picture}

%% file: 1geo-ds2.pstex_t
\begin{picture}(0,0)%
\includegraphics{1geo-ds2.pstex}%
\end{picture}%
\setlength{\unitlength}{1989sp}%
\begingroup\makeatletter\ifx\SetFigFont\undefined%
\gdef\SetFigFont#1#2#3#4#5{%
  \reset@font\fontsize{#1}{#2pt}%
  \fontfamily{#3}\fontseries{#4}\fontshape{#5}%
  \selectfont}%
\fi\endgroup%
\begin{picture}(4548,6895)(887,-7843)
\put(2641,-6061){\makebox(0,0)[lb]{\smash{\SetFigFont{11}{12.0}{\rmdefault}{\mddefault}{\updefault}{$\Uu_1^-$}%
}}}
\put(4036,-6361){\makebox(0,0)[lb]{\smash{\SetFigFont{11}{12.0}{\rmdefault}{\mddefault}{\updefault}{$\Vv_1$}%
}}}
\put(4726,-5821){\makebox(0,0)[lb]{\smash{\SetFigFont{11}{12.0}{\rmdefault}{\mddefault}{\updefault}{$\Uu_1^+$}%
}}}
\end{picture}

%% file: positive-rotation.pstex_t
\begin{picture}(0,0)%
\includegraphics{positive-rotation.pstex}%
\end{picture}%
\setlength{\unitlength}{3729sp}%
\begingroup\makeatletter\ifx\SetFigFont\undefined%
\gdef\SetFigFont#1#2#3#4#5{%
  \reset@font\fontsize{#1}{#2pt}%
  \fontfamily{#3}\fontseries{#4}\fontshape{#5}%
  \selectfont}%
\fi\endgroup%
\begin{picture}(4844,2462)(1369,-5091)
\put(2826,-4420){\makebox(0,0)[lb]{\smash{\SetFigFont{11}{12.0}{\rmdefault}{\mddefault}{\updefault}{$(x_-,x_-)$}%
}}}
\put(4319,-4105){\makebox(0,0)[lb]{\smash{\SetFigFont{11}{12.0}{\rmdefault}{\mddefault}{\updefault}{$p$}%
}}}
\put(4079,-3504){\makebox(0,0)[lb]{\smash{\SetFigFont{11}{12.0}{\rmdefault}{\mddefault}{\updefault}{$L$}%
}}}
\put(4251,-5042){\makebox(0,0)[lb]{\smash{\SetFigFont{11}{12.0}{\rmdefault}{\mddefault}{\updefault}{$R$}%
}}}
\put(4116,-4502){\makebox(0,0)[lb]{\smash{\SetFigFont{11}{12.0}{\rmdefault}{\mddefault}{\updefault}{$(\exp(-tX),Id)p$}%
}}}
\put(3615,-3560){\makebox(0,0)[lb]{\smash{\SetFigFont{11}{12.0}{\rmdefault}{\mddefault}{\updefault}{$l$}%
}}}
\end{picture}

%% file: 1geo-Ads.pstex_t
\begin{picture}(0,0)%
\includegraphics{1geo-Ads.pstex}%
\end{picture}%
\setlength{\unitlength}{2569sp}%
\begingroup\makeatletter\ifx\SetFigFont\undefined%
\gdef\SetFigFont#1#2#3#4#5{%
  \reset@font\fontsize{#1}{#2pt}%
  \fontfamily{#3}\fontseries{#4}\fontshape{#5}%
  \selectfont}%
\fi\endgroup%
\begin{picture}(5424,5072)(3139,-7093)
\put(4936,-5356){\makebox(0,0)[lb]{\smash{\SetFigFont{11}{12.0}{\rmdefault}{\mddefault}{\updefault}{$\Pp^-$}%
}}}
\put(6106,-5386){\makebox(0,0)[lb]{\smash{\SetFigFont{11}{12.0}{\rmdefault}{\mddefault}{\updefault}{$\Qq$}%
}}}
\put(6916,-4861){\makebox(0,0)[lb]{\smash{\SetFigFont{11}{12.0}{\rmdefault}{\mddefault}{\updefault}{$\Pp^+$}%
}}}
\end{picture}

%% file: HANDADS.tex
\section{Causal AdS spacetimes, earthquakes and black holes}\label{moreAdS}
Beyond the classification achieved in the previous Section, the AdS case
displays a rich phenomenology that we are going to point out.

\subsection{On holonomy pregnancy} 
Let us recall first the following results of \cite{M}, in the case of
compact $S$.
\begin{teo}\label{AdS_compact} 
If $S$ is compact, and $Y\in \Mm\Gg\Hh_{-1}(S)$,  then:
\smallskip

(a) The holonomy $h=(h_L,h_R)$ of $Y$ is made by a couple of Fuchsian
representations of $\pi_1(S)$, and every such a couple arises in this
way (by varying $Y$).
\smallskip

(b) $Y$ is completely determined by its holonomy $h=(h_L,h_R)$. In
fact $Y=\Yy(C)$, where $C$ is the graph in $S^1_\infty\times
S^1_\infty=\partial\mx_{-1}$ of the unique orientation preserving
homeomorphism that conjugates the action of $h_L$ on $S^1_\infty=
\partial \mh^2$ with the one of $h_R$. This curve $C$ is the unique
$h$-invariant curve on $\partial \mx_{-1}$.
\end{teo}

In~\cite{Ba}(2, 3) we can find the following generalization of point (a).  
Here we use the notations of Section \ref{ML(S)}.
\begin{prop}\label{discrete-rep} 
Let  $Y\in \Mm\Gg\Hh_{-1}(S)$ be of finite type, with
holonomy representation $$h=(h_L,h_R):\pi_1(Y)\rightarrow
PSL(2,\mr)\times PSL(2,\mr) \ .$$ Then both $h_L$ and $h_R$ are holonomy
representations of hyperbolic structures belonging to $\Tt(S)$.
Conversely given a pair of representations $h=(h_L,h_R)$ corresponding
to elements of $\Tt(S)$, then there exists a spacetime $Y\in
\Mm\Gg\Hh_{-1}(S)$ whose holonomy is $h$.
\end{prop}
Concerning point (b), the following partial generalization holds.
Here we adopt the notations of Corollary \ref{preserve-type}.

\begin{prop}$Y\in \Mm\Gg\Hh_{-1}(S)$  is
completely determined by its holonomy providing that $Y=Y^\lambda_{-1}$,
for some $(F,\lambda)\in \Vv_\cG(F)\cap \Mm\Ll_\cG(F)^0$.
\end{prop}
This is essentially a consequence of the proof of the Earthquake
Theorem considered below. On the other hand, non-equivalent spacetimes
in $\Mm\Gg\Hh_{-1}(S)$ can actually share the same holonomy. This is
the theme of the following construction.

\subsection{Canonical causal AdS spacetimes with prescribed holonomy }
\label{CAUSAL}

We mostly refer to~\cite{Ba}(2, 3).  Let us fix a representation
$h=(h_L,h_R)$ of $\pi_1(S,p_0)$ as in Theorem \ref{discrete-rep}. We
stress that the representation is fixed, not only its conjugation
class; for this reason we have also fixed a base point $p_0\in S$.  We
consider the domain
$$\tilde\Omega(h)$$
of points $x\in \mx_{-1}$ such that, for every
$\gamma \in \pi_1(S,p_0)$, $x$ and $h(\gamma)(x)$ are not causally
related. We have:
\begin{prop}  
$\tilde\Omega(h)$ is simply connected and $h$-invariant; the action of
$\pi_1(S,p_0)$ on it is free and properly discontinuous. The quotient,
say $\Omega(h)$, is a causal AdS spacetime homeomorphic to
$S\times\mr$.
\end{prop}

Let us consider now $$\Mm\Gg\Ll(h)$$ the set of all AdS $MGH$ spacetimes
$Y$ homeomorphic to $S\times\mr$, determined by a compatible couple of
developing map and holonomy representation $(d_Y,h_Y)$ such that
$h_Y=h$. Equivalently, we are considering spacetimes homeomorphic to
$S\times\mr$ of the form $Y=\Yy(C)/h$ such that the nowhere time-like
curve at infinity $C$ is $h$-invariant. Note again that each such a
maximal globally hyperbolic spacetime is fixed and {\it not}
considered up to Teichm\"uller-like equivalence. However, it is not
hard to see that:
\begin{lem} The natural map $\Mm\Gg\Ll(h)\to \Mm\Gg\Hh_{-1}(S)$ is
injective.
\end{lem}
\begin{remark}{\rm If $h'= ghg^{-1}$ is conjugate to $h$, then
$\tilde\Omega(h')= g\tilde\Omega(h)$, as well $\Mm\Gg\Ll(h')=
g\Mm\Gg\Ll(h)$, so that they have the same image in
$\Mm\Gg\Hh_{-1}(S)$. By fixing a representative $h$ in any conjugation
class, we get in this way a partition of  $\Mm\Gg\Hh_{-1}(S)$.
}
\end{remark} 
Since every $Y \in \Mm\Gg\Ll(h)$ is causal, we have
\begin{lem}\label{extend} 
There is a natural $h$-invariant embedding of $\Yy(C)$ in
$\tilde\Omega(h)$, hence of $Y$ in $\Omega(h)$.
\end{lem}
In fact we can prove 
\begin{prop}  $\Omega(h)$ is the union of the $Y\in \Mm\Gg\Ll(h)$,
as well as  $\tilde\Omega(h)$ is the union of the $h$-invariant $\Yy(C)$'s.
\end{prop}
See below for a description of the curves $C$ arising in this way.
\smallskip

{\bf $\partial_\infty\tilde\Omega(h)$ and the limit set.}
The inclusion of Lemma \ref{extend} is in general {\it strict}, in
particular $\Omega(h)$ can be {\it not} globally hyperbolic.  In a
sense, it is just the {\it maximal causal extension} of every globally
hyperbolic $Y\in \Mm\Gg\Ll(h)$. The reason is that the adherence,
$$\partial_\infty\tilde\Omega(h)$$ of $\tilde\Omega(h)$ on the
boundary of $\mx_{-1}$ can have non-empty interior part. Such an
adherence can be explicitly described by means of the holonomy of the
peripheral loops.  Let $\gamma\in\pi_1(S,p_0)$ freely homotopic to a
loop surrounding a point in $V$. If $h_L(\gamma)$ (resp. $h_R(\gamma)$)
is of hyperbolic type, we set $I_L(\gamma)$ (resp.  $I_R(\gamma)$) to
be the interval of $S^1$ whose end-points are the fixed points of
$h_L(\gamma)$ (resp. $h_R(\gamma)$) and that does not meet the limit
set of $h_L(\gamma)$ (resp. $h_R(\gamma)$). If $h_L(\gamma)$ is
parabolic, then let $I_L(\gamma)$ be the fixed point of
$h_L(\gamma)$. Similarly for $I_R(\gamma)$.  Then the ``rectangle''
$R(\gamma)=I_L(\gamma)\times I_R(\gamma)$ is contained in the
adherence of $\tilde\Omega(h)$, and in fact
$\partial_\infty\tilde\Omega(h)$ is the closure $\bigcup_\gamma
R(\gamma)$. The closure of the complement of the union of these
rectangles $R(\gamma)$ in $\partial_\infty\tilde\Omega(h)$ can be
regarded as a {\it limit set} 
$$\Lambda = \Lambda(h)$$ in the sense that it is contained in the
closure of the orbit of any point $x\in\overline{\mx}_{-1}$.  

A rectangle is {\it non-degenerate} if both $h_L(\gamma),h_R(\gamma)$
are of hyperbolic type.

\begin{lem} 
If some rectangle $R(\gamma)$ is non-degenerate, then the interior of
$\partial_\infty\tilde\Omega(h)$ is not empty.  $\Omega(h)$ is not
globally hyperbolic if and only if there is some non-degenerate
rectangles .
\end{lem}
For it is possible to find points $p,q$ close to $R(\gamma)$ such that
$\fut(p)\cap\pass(q)$ is not pre-compact in $\mx_{-1}$ and this
contradicts the global hyperbolicity.
\smallskip

{\bf The asymptotic regions.}  Notice that a non degenerate rectangle
$R(\gamma)$ has exactly two vertices that are the end-points of a
spacelike geodesic $l_\gamma$ that is invariant for
$(h_L(\gamma),h_R(\gamma))$. The boundary lines of $R(\gamma)$
together with $l_\gamma$ span a surface $\tilde H(\gamma)$ embedded in
$\tilde\Omega(h)$ made by two null triangles intersecting at
$l_\gamma$. This surface divides $\tilde\Omega(h)$ in two
components. The component whose adherence in
$\partial_\infty\tilde\Omega(h)$ is $R(\gamma)$ is called an {\it
asymptotic region} of $\tilde\Omega(h)$, denoted by $\tilde A(\gamma)$,
$\tilde H(\gamma)$ is its {\it horizon}. $(A(\gamma),H(\gamma))$ is
invariant for $(h_L(\gamma),h_R(\gamma))$ and the quotient embeds in
$\Omega(h)$ giving us an asymptotic region $A(\gamma)$ with horizon
$H(\gamma)$. This last is the union of two null annuli along a
spacelike closed geodesic. The length of this spacelike geodesic is
called the \emph{size} of the horizon, whereas the \emph{momentum} is
the twist factor for the parallel transport along it.  If $l_L,l_R$
are the translation lengths of $h_L(\gamma)$ and $h_R(\gamma)$, the
size is simply $$s=(l_L+l_R)/2$$ whereas the momentum is
$$m=(l_L-l_R)/2 \ . $$ $\Omega(h)$ has exactly $k$ asymptotic regions,
where {\it $k$ is the number of points $p\in V$ such that the surrounding
circle is of hyperbolic type for both $h_L$ and $h_R$}.
\smallskip

{\bf More about $\Yy(C) \subset \tilde\Omega(h)$.}  Clearly the
$h$-invariant curve at infinity $C$ is contained in
$\partial_\infty\tilde\Omega(h)$.  On the other hand, every nowhere
time-like meridian of $\partial \mx_{-1}$ contained in
$\partial_\infty\tilde\Omega(h)$ is determined by drawing in each
non-degenerate rectangle $R(\gamma)$ an arc $l_\gamma$ joining the
vertices that are the end-points of the spacelike geodesic of the
corresponding horizon. In the degenerate case the segment $l_\gamma$
coincide with $R(\gamma)$. The closure $C$ of the union of these
$l_\gamma$'s is a nowhere timelike meridian; moreover, if the segments
are chosen in $h$-invariant way (that is
$l_{\alpha\gamma\alpha^{-1}}=h(\alpha) l_\gamma$), then $C$ is is the
curve at infinity of the universal covering of some $Y \subset \Omega(h)$.

\subsection{AdS bending and Earthquake Theorems}\label{ge:quake}
By extending  the arguments given in \cite{M} in the
case of compact $S$, we have for a general $S$ of finite type:
\begin{prop}\label{earth:class:prop}
Let $Y \in \Mm\Gg\Hh_{-1}(S)$ be encoded by $(F,\lambda)\in
\Mm\Ll(S)$, and $h=(h_L,h_R)$ be its holonomy.  Then $h_L$ ($h_R$) is
the holonomy of the surface $F_L = \beta^L_\lambda(F)$ ($F_R=
\beta^R_\lambda(F)$), that is the surface in $\tilde {\Tt}(S)$
obtained by the left (right) earthquake on $F$ along $\lambda$.
\end{prop}   

We stress that Proposition~\ref{earth:class:prop}, together with
Theorem \ref{AdS_compact}, actually give an AdS proof of the
Earthquake Theorem \ref{quake-teo} when $S$ is compact. For, given
$F^0,F^1$ two hyperbolic structures on a compact surface $S$, there
exists a unique spacetime $Y^\lambda_{-1} \in \Mm\Gg\Hh_{-1}(S)$ whose
holonomy is $h=(h^0,h^1)$, where $h^j$ is the hyperbolic holonomy of
$F^j$. Then the left earthquake along $2\lambda$ transforms $F^0$ into
$F^1$.
\smallskip

We consider the subset
$$\Mm\Gg\Hh_\cG(h)$$ of $\Mm\Gg\Hh(h)$ consisting of the spacetimes
$Y$ that satisfy the further condition of being encoded by couples
$(F,\lambda)\in\Mm\Ll_\cG(S)$.  In order to get such an AdS proof of
the full Earthquake Theorem \ref{quake-teo}, we need to characterize
the spacetimes $Y=\Yy(C)/h \in \Mm\Gg\Hh_\cG(h)$ in terms of the curve
at infinity $C$.  Consider again the general description of an
$h$-invariant meridian $C$ given above.  A case of particular interest
is when the segments $l_\gamma$ are chosen on the boundary of
$R(\gamma)$. Meridians $C$ obtained in this way are called {\it
extremal}.  Notice that for each asymptotic region there are only two
ways to chose such an arc: an upper extremal arc and a lower extremal
arc.  Thus, there are exactly $2^k$ $h$-invariant extremal arcs where
$k$ is defined as above. This holds also when $k=0$; in such a case 
$\Omega(h)=\Yy(C)/h$ is globally hyperbolic, and $C$ is its extremal
meridian. Finally we the following nice geometric
characterization (see \cite{BSK}).
\begin{prop}\label{geom-char} $\Yy(C)$ is the universal covering of some 
$Y\in \Mm\Gg\Hh_\cG(h)$ if and only if $C$ is an $h$-invariant extremal
meridians.
\end{prop}

\begin{cor}\label{omega-mgh}
$\Omega(h)$ is globally hyperbolic if and only if
it belongs to $\Mm\Gg\Hh_\cG(h)$ and is encoded by $(F,\lambda)$ 
such that $F\in \Tt_{g,r}$ and the lamination does not enter the cusps
\end{cor}

We are ready to prove Theorem \ref{quake-teo}. Let $F^0$, $F^1$ be the
the interior of the convex cores of $\mh^2/h^0$, $\mh^2/h^1$
respectively, that are both homeomorphic to $S$. Set $h=(h^0,h^1)$ and
take $\Omega(h)$.  Let us apply Proposition \ref{earth:class:prop} to
every $Y\in \Mm\Gg\Hh_\cG(h)$, encoded by some $(F,\lambda)\in
\Mm\Ll_\cG(S)$. As the convex cores are uniquely determined by the
holonomy, and $\Tt_\cG(S)$ is closed under earthquakes, it follows
that $F^0 = \beta^L_\lambda(F)$, $F^1= \beta^R_\lambda(F)$, so that
$F^1 = \beta^L_{2\lambda}(F^0)$. The determined lack of uniqueness in
Theorem \ref{quake-teo}, the ``enhanced'' version \ref{quake-teo-bis},
as well as Corollary \ref{nec-comp} are now rather easy consequences of
Proposition \ref{geom-char}, Lemma \ref{ray-quake}, and the definition of
the (enhanced) quake-flow.

\subsection{Convex core of $\Omega(h)$ and black hole}
\label{BH_convexcore}
We denote by $$\Omega_\cG(h)\subset \Omega(h)$$ the union of the
spacetimes belonging to $\Mm\Gg\Hh_\cG(h)$. We do similarly for
$\tilde\Omega_\cG(h)\subset \tilde\Omega(h)$. It follows from
Proposition \ref{geom-char} that the connected components of
$\Omega(h)\setminus \Omega_\cG(h)$ coincide with the asymptotic
regions defined above. Similarly for $\tilde\Omega(h)\setminus
\tilde\Omega_\cG(h)$.

The limit set $\Lambda$ is contained in the adherence of
$\tilde\Omega_\cG(h)$ which is the union of a finite number of
globally hyperbolic spacetimes. Hence there is a spacelike plane $P$
that does not intersect $\tilde\Omega_\cG(h)$, so that we can take the
convex hull
$$\tilde\Kk(h)$$ of $\Lambda$ in $\mr^3=\Pm^3\setminus\hat P$, where
$\hat P$ is the projective plane containing $P$. It turns that
$\tilde\Kk(h)$ is contained in the closure of $\tilde\Omega_\cG(h)$, it is
$h$-invariant and does not depend on the choice of $P$. $\tilde\Kk(h)$
is called the {\it convex core} of $\tilde\Omega(h)$, as well as its
quotient $$\Kk(h)$$ is the convex core of $\Omega(h)$. We can see that
$\tilde\Omega(h)$ coincides with the set of points in $\mx_{-1}$ whose
dual plane does not intersect $\tilde\Kk(h)$, and that every plane
dual to some point of $\tilde\Kk(h)$ does not intersect
$\tilde\Omega(h)$.  The boundary of $\tilde\Kk(h)$ contains the
spacelike geodesics of the horizons of $\tilde\Omega_\cG(h)$. Such
geodesics disconnects $\partial \tilde\Kk(h)$ into two $h$-invariant
pleated surfaces whose quotients are homeomorphic to $S$. One, say 
 $\partial_+ \tilde\Kk(h)$, is in the future of the other one, say
$\partial_- \tilde\Kk(h)$, an they are called the {\it future} and the 
{\it past boundary} of  $\tilde\Kk(h)$ respectively.  
It turns that $\partial_+\tilde\Kk(h)$ is obtained via
the AdS bending of $(F_+,\lambda_+)\in \Mm\Ll_\cG(S)$ (according to Section
\ref{WR}) so that it is  the future boundary of the past part of a
specific $\Yy(C_+)/h \in \Mm\Gg\Hh_\cG(h)$; in fact:
\smallskip

{\it The extremal $h$-invariant meridian $C_+$ is obtained
by taking the {\rm lower} extremal arc in each rectangle.}
\smallskip

Similarly $\partial_-\tilde\Kk(h)$ is the {\it past} boundary of the
{\it future} part of a specific spacetime $\Yy(C_-)/h \in
\Mm\Gg\Hh_\cG(h)$ ( whose future boundary of the past part is obtained
by bending a certain $(F_-,\lambda_-)\in \Mm\Ll_\cG(S)$).  The
corresponding extremal $h$-invariant meridian $C_-$ is obtained by
taking the {\it upper} extremal arc in each rectangle. This makes
sense also when $\Omega(h)$ is globally hyperbolic; in such a case
$C_- = C_+$.

Assume now that $\Omega(h)$ is not globally hyperbolic.  For each
boundary component $c_i$ of $F_+$, we have that $l_{c_i}$ is the size
of the corresponding horizon, whereas ${\rm I}_{c_i}(\lambda_+)$ is
the corresponding momentum. It follows that $\lambda_+$ belongs to the
closure of $\Vv_\cG(F_+)$ (recall Corollary \ref{preserve-type}).  In fact
this property uniquely characterizes $Y_+$ within
$\Mm\Gg\Hh_\cG(h)$. In particular this selects a {\it privileged} one
among the earthquakes of Theorem \ref{quake-teo}. For
$(F_-,\lambda_-)$ we have the somehow {\it opposite} behaviour, that
is for every boundary component ${\rm I}_{c_i}(\lambda_-)>l_{c_i}$.

Set
$$\tilde B(h)= \Yy(C_-),\ \ B(h)=\tilde B(h)/h, \ \ \tilde W(h)= \Yy(C_+),
\ \ W(h)= \tilde W(h)/h \ . $$ Denote by $KB(h)$, $KW(h)$ the respective
convex cores as MGH spacetimes.
We have
\begin{prop}\label{B-W}  (1) $\Kk(h) = KB(h) \cap KW(h)$.
\smallskip

(2) $\tilde \Omega_\cG(h)= \tilde B(h) \cup \tilde W(h)$,
$\Omega_\cG(h) = B(h) \cup W(h)$.
\end{prop}

\begin{figure}[ht]
\begin{center}
\includegraphics[width=6cm]{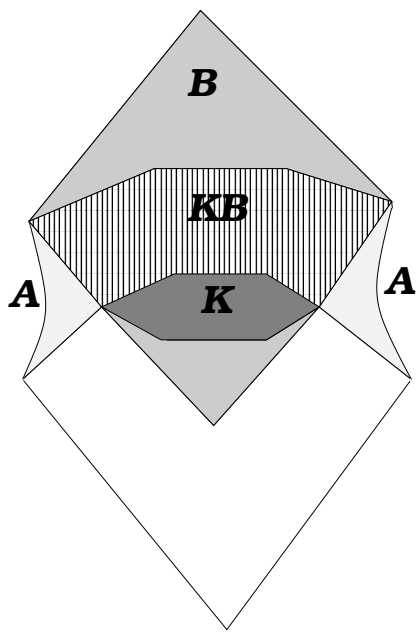}
\caption{\label{BH-fig} A schematic picture of $\Omega$ with its
convex core $\Kk$ and asymptotic regions $\Aa$. It contains the black
hole $\Bb$ with is convex core $\Kk\Bb$ containing $\Kk$. }
\end{center}
\end{figure}

In Physics literature the special globally hyperbolic spacetime $B(h)$
($W(h)$) is known as the {\it multi black hole} ( {\it multi white
hole}) contained in the causal spacetime $\Omega(h)$. The attribute
``multi'' mostly refer to the fact that it has a ``multi''
horizon. $B(h)$ looks like a honest black hole in the sense that every
future inextensible causal curve emanating from any event in $B(h)$
never leave $B(h)$ and eventually reach the final singularity
$\Sigma_+$ of $B(h)$ in finite time. In particular, lightlike rays
emanating from $B(h)$ do not reach $\partial_\infty\Omega(h)$.  So the
final singularity $\Sigma_+$ is an actual singularity for the
spacetime $\Omega(h)$ itself, as it reflects its future timelike
geodesic incompleteness ( at the initial singularity of $B(h)$ that is
contained in its interior, $\Omega (h)$ is perfectly non singular).
The initial singularity $\Sigma_-$ of the white hole $W(h)$ plays a
similar r\^ole with respect to the past. However, $\Sigma_+$ is
``censured'' by the multihorizon of $B(h)$, while $\Sigma_-$ is a
``naked'' singularity. In Figure \ref{BH-fig} we see a schematic picture of
$\Omega(h)$ with its convex core and its black hole.

{\bf Asymptotic regions and BTZ black holes.}
Every asymptotic region $A=A(\gamma)$ of $\Omega(h)$ has by itself a
natural extension to a maximal {\it causal} AdS spacetime
$\Bb=\Bb(\gamma)$, homeomorphic to $(S^1\times \R)\times \R)$. $\Bb$
contains a maximal globally hyperbolic spacetime $\Bb\Hh$, with a
complete Cauchy surface homeomorphic to the annulus $S^1\times
\R$, which is known as the {\it BTZ black hole} contained in $\Bb$ (see
\cite{BTZ}, \cite{Ca}).  $\Bb$ has been particularly studied because
it supports {\it Kerr-like metrics} with several qualitative analogies
with the classical rotating black hole solutions of (3+1) gravity. Let
us briefly recall this matter. It is convenient to lift $\mx_{-1}=
PSL(2,\R)$ in $\hat{\mx}_{-1}=SL(2,\R)$ so that it is given by the
matrices of the form
$$X=\left(\begin{array}{cc}
T_1+X_1& T_2+X_2\\
-T_2+X_2& T_1-X_1
\end{array}\right)$$ 
such that ${\rm det}(X)=1$, $0<T_1^2 -X_1^2 <1$, $X_1, T_1$ have a
definite sign. We fix also a suitable $SL(2,\R)$-lifting of the isometry
$(h_L(\gamma),h_R(\gamma))$ corresponding as above to the given asymptotic
region. Let us assume for simplicity that it is of the form
$$\left(\left(\begin{array}{cc}
\exp(r_+ - r_-) & 0\\
0 & \exp(r_- - r_+)
\end{array}\right),\left(\begin{array}{cc}
\exp((r_+ + r_-)) &  0\\
0 & \exp(-(r_+ + r_-))
\end{array}\right)\right)$$
\smallskip
and that $r_+> r_- \geq 0$. This isometry generates 
a group $\Gg$ that acts on the whole of $\hat{\mx}_{-1}$, with a constant
vector field $\xi$ as infinitesimal generator, and we have $q(\xi)
=(T_2^2-X_2^2)r_+ + (T_1^2 - X^2_1)r_-$ ($q$ is defined in Section
\ref{grav}). Roughly speaking $\tilde\Bb$ is the {\it maximal} region of
$\mx_{-1}$ such that:
\smallskip

(1) $q(\xi)>0$ on $\tilde\Bb$, so that we can take the function
    $r=q(\xi)^{1/2}>0$;

(2) $\{r_+>r>r_- \subset \tilde\Bb$;
\smallskip

(3) $\tilde\Bb$ is $\Gg$-invariant, the group acts
nicely and the quotient $\Bb$ is a causal spacetime homeomorphic to
$(S^1\times \R)\times \R$. 
\medskip

$\tilde \Bb$ admits a $\Gg$-invariant ``tiling'' by regions of three
types I, II, III contained in $\{r>r_+\}$, $\{r_+>r>r_-\}$,
$\{r_->r\}$ respectively. Each region is bounded by suitable null
horizons at which $r=r_\pm$. We can see that our asymptotic
regions $\tilde A$ are of type III.

By ``joining'' (the lifting of) the spacelike line $l_\gamma$ with the
two liftings of the dual line $l_\gamma^*$ respectively, we get two
``tetrahedra'' say $\tilde \Bb\Hh$ and $\tilde \Ww\Hh$ embedded in
$\tilde \Bb$ intersecting at $l_\gamma$. These are the two regions of
type II that form the whole of $\{r_+>r>r_-\}$. One projects onto the
BTZ black hole $\Bb\Hh$, the other one covers the white hole embedded
in $\Bb$, say $\Ww\Hh$.  Note that both $\tilde \Bb\Hh$ and $\Ww\Hh$
are instances of ``degenerate'' globally hyperbolic spacetimes in the
sense of Section \ref{AdS}.
\smallskip

For suitable coordinates $(v,r,\phi)$ on $\Bb$, where $(r,\phi)$
look like polar coordinates on the $v$-level surfaces, the Kerr-like metric
is of the form
$$ds^2 = (M-r^2)dv^2 + f^{-1}dr^2 + r^2 d\phi^2 -Jdvd\phi$$ where
$$M=r^2_+ + r^2_-, \ \ J=2r_+r_- ,\ \ M\geq J, \ \ f= -M+r^2 +
\frac{J^2}{4r^2}$$ and they are related to the previously defined
``size'' and ``momentum'' by
$$M+J= s^2,\ \ M-J=m^2 \ . $$ Each region of $\Bb$ support this
metric, the null horizons of the regions being just ``coordinate
singularities''.  
\smallskip

BTZ black holes naturally arise in the framework of Wick
rotation-rescaling theory for the {\it elementary} surfaces of finite
type, that is having {\it Abelian} fundamental group: $S=S^1\times \R$ and
$S= S^1\times S^1$. This displays an interesting r\^ole of {\it
quadratic differentials} instead of geodesic laminations. See Chapter
7 of \cite{Be-Bo} for more details.

\subsection{(Broken) $T$-symmetry}\label{T-symm}
Let $Y\in \Mm\Gg\Hh_{-1}(S)$. By reversing the time orientation we get
onother spacetime $Y^* \in \Mm\Gg\Hh_{-1}$.  This involution is called
$T$-{\it symmetry} as well the involution induced on $\Mm\Ll(S)$ via
the map $\mG_{-1}$. If $\Yy(C)$ is the universal covering of $Y$, then
the universal covering of $Y^*$ is $\Yy^* = \Yy(C^*)$, where $C^*$ is
the image of the curve $C$ under the involution of
$\partial\mx_{-1}=S^1_\infty\times S^1_\infty$
$$ (x,y)\mapsto(y,x)\ .$$ Moreover, the holonomy $h^*$ of $Y^*$ is
obtained by exchanging the components of the holonomy of $h$ of $Y$
$$ h=(h_-, h_+) \leftrightarrow h^*=(h_+, h_-) \ .$$ If $B(h)$ is the
black hole of $\Omega(h)$, then $B(h)^* = W(h^*)$.  The opposite
behaviour ``${\rm I}_c(\lambda_+)\leq l_c$ vs ${\rm
I}_c(\lambda_-)>l_c$'' at the future boundary of the respective past
parts (see above), can be considered as the basic feature of ``broken
$T$-symmetry''.  A particular instance is when $B(h)$ is encoded by
$(F,\lambda)$ such that $F\in \Tt_{g,r} = \Tt(S)\cap \Tt_\cG(S)$ (the
smallest stratum of $\Tt_\cG(S)$), and $\lambda$ enters the cusps (in
the Figure we show an example of $B(h)$ where $F$ has $g=0$, $r=3$,
and the lamination is like in Example \ref{exa-via-twist-shear} with
respect to a standard ideal triangulation of $F$ by two triangles).
In this case $F^*$ belongs indeed to a higher dimensional cell of
$\Tt_\cG(S)$ and the white hole $W(h^*)$ has the property that
${\rm I}_{\lambda^*}(c) = l(c)$ at every boundary curve, and every
asymptotic region has null momentum.

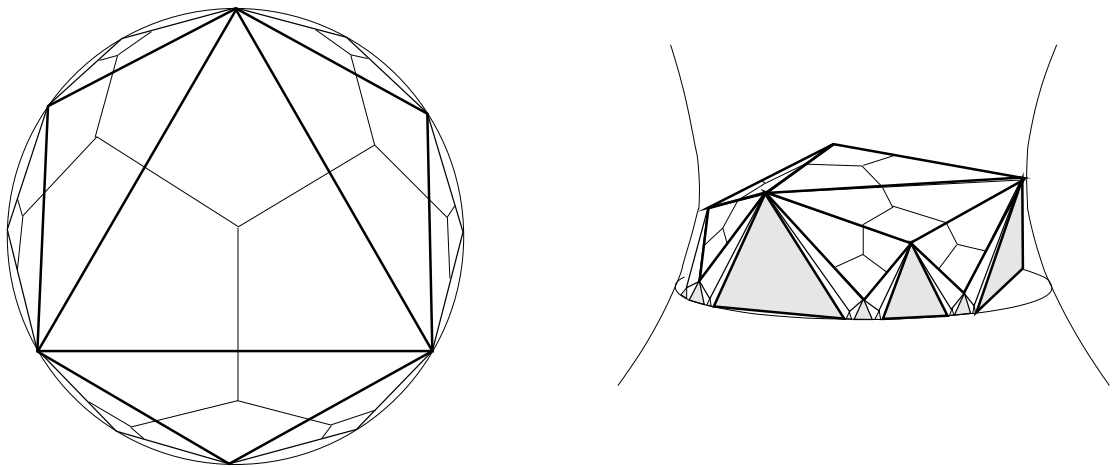
\begin{figure}
\begin{center}
\input{3Cusp-bend.pstex_t}
\caption{The convex core of a black hole $B(h)$. On the left the
    lamination with its dual spine. On the right the bending of
    $\mh^2$ along $\lambda$ in $\mx_{-1}$. Grey regions are null
    components of the past boundary of $\Kk_\lambda$.}
\end{center}
\end{figure}

%% file: 3Cusp-bend.pstex_t
\begin{picture}(0,0)%
\includegraphics{3Cusp-bend.pstex}%
\end{picture}%
\setlength{\unitlength}{2072sp}%
\begingroup\makeatletter\ifx\SetFigFont\undefined%
\gdef\SetFigFont#1#2#3#4#5{%
  \reset@font\fontsize{#1}{#2pt}%
  \fontfamily{#3}\fontseries{#4}\fontshape{#5}%
  \selectfont}%
\fi\endgroup%
\begin{picture}(13219,5518)(1689,-6581)
\end{picture}

%% file: HANDPART.tex
\section{Including particles}\label{part}
In 3-dimensional gravity massive point particles can be modeled as
cone singularities along timelike lines. In particular, the rest mass
$m$ of a particle is related to the curvature $k$ concentrated along
its timelike geodesic ``world line'' by
$$k=2\pi m, \ \ k=2\pi - \alpha$$ where $\alpha$ is the cone angle. If
we require that the mass is positive, then it is bounded by $0\leq m
\leq 1$, while $2\pi \geq \alpha \geq 0$. However, there are no
real geometric reasons to exclude cone angles bigger than $2\pi$.
\smallskip

It is a natural question whether Wick Rotation-rescaling theory does
apply also on cone spacetimes. In such a perspective,
it is quite natural to extend the space $\Tt_\cG(S)$ defined in
Section \ref{ML(S)}, by extending the notion of ``type'' to $\theta =
V_\Hh \cup V_C$, and allowing hyperbolic structures $F$ on $S$ whose
completion $F^\Cc$ is compact, as usual has geodesic boundary
components corresponding to the points of $V_\Hh$, and possibly has
{\it conical singularities} at the points of $V_C$. In particular we
allow also holonomy of {\it elliptic type} at the circles surrounding
these points. The parabolic holonomies correspond now to cone angles
equal to $0$, hence to particles of extremal mass. In order to preserve
the conical structure, we can consider measured geodesic laminations
on such cone surfaces $F$ that have {\it compact} support $L$ in $F$,
that is whose closure in $F^\Cc$ does not intersect the
singularities.
\smallskip

For the sake of simplicity (and following \cite{BS} to which we will
refer for most results stated in this Section), from now on we will
consider the particular case such that $V_\Hh = \emptyset$. Not even
in this simplest case a complete answer to the above question is
known.  Only a few partial results are known, mostly concerning the
case of ``small'' cone angles ($<\pi$), or equivalently the
case of particles with ``big'' masses.
\smallskip

Let $\hat S$, $S = \hat S \setminus V$, $V =\{p_1,\ldots , p_r\}$ be
as usual.  Let $g$ be the genus of $\hat S$. Fix an $r$-tuple of angles
$\Aa=(\alpha_1,\ldots,\alpha_r)$, such that the ``Gauss-Bonnet
inequality''
$$\sum_j (1-\frac{\alpha_j}{2\pi})> 2-2g $$ holds; notice that we are not
requiring here that the cone angles are smaller than $2\pi$.

We denote by
\[
     \Tt_\cG(S,\Aa)
\]
the Theichm\"uller space of hyperbolic structures $F$ on $S$ whose
completion $F^\Cc$ has conical singularities at $p_1,\ldots,p_r$, of
cone angles $\alpha_1,\ldots,\alpha_r$.

By a general result of Troyanov \cite{Tr}, we have
\begin{prop}\label{troyanov} The natural map 
$$\Tt_\cG (S,\Aa)\rightarrow\Tt_{g,r}$$ that
associates to every $F\in \Tt_\cG (S,\Aa)$ the unique complete hyperbolic
structure of finite area  $\tilde F$ on $S$ in the same conformal
class of $F$, is a bijection.
\end{prop}
This means in particular that $\Tt_\cG (S,\Aa)$ is not empty. If $\alpha_j
= 0$ for every $j$, then $\tilde F = F$ and $\Tt_\cG (S,0)$ just coincides
with $\Tt_{g,r}$.
\smallskip

For every $F\in \Tt_\cG(S,\Aa)$, we denote by $\Mm\Ll_\cG(F,\Aa)$ the space of
measured geodesic laminations on $F$ with compact support. The space
of all such $(F,\lambda)$'s is denoted by $\Mm\Ll_\cG (S,\Aa)$.  When $\Aa
= 0$, then $\Mm\Ll_\cG (F,0)$ just coincides with $\Mm\Ll_\cG (\tilde F)^0$
(defined in Section \ref{ML(S)}). More generally we have:

\begin{prop}\label{zzz:prop} Assume that for every $j$, $\alpha_j <\pi$. 
Then there is a natural identification between $\Mm\Ll_\cG(F,\Aa)$ and 
$\Mm\Ll_\cG(\tilde F)^0$.
\end{prop}
We give a brief sketch of the proof of this proposition, assuming some
familiarity with ``train-tracks''. Since $\alpha_i<\pi$, for small
$\eps >0$, the complement, say $\Sigma_\eps$, of a regular
neighborhood of $V$ in $F^\Cc$ of ray $\eps$ is convex.  As a
consequence , any non-peripheral loop on $F$ admits a geodesic
representative whose distance from $V$ is at least $\eps$.

Given $\tilde \lambda \in \Mm\Ll_\cG(\tilde F)^0$, its support is
contained in $\Sigma_\eps$, for $\eps$ sufficiently small. Since
$\Sigma_\eps$ is compact, it follows that the leaves of $\tilde
\lambda$ are quasi-geodesic in $F$. Since $\Sigma_\eps$ is convex,
they can be stretched to become geodesic with respect to $F$. The union
of all these leaves makes a geodesic lamination $\lambda$ on $F$. A
train-track carrying $\tilde \lambda$ carries also $\lambda$ so this
last can be equipped with a transverse measure corresponding to the
measure on $\tilde \lambda$.
 
\begin{remark} {\rm 
If some cone angle $\alpha_i$ is bigger than $\pi$, then
Proposition~\ref{zzz:prop} fails. In fact it is not difficult to
construct a surface with cone angles bigger than $\pi$ and a loop $c$
whose geodesic representative passes through the singular point.}
\end{remark}

\subsection {Maximal globally hyperbolic spacetimes with particles}
Since the causal structure of a spacetime with  timelike geodesic world lines
of conical singularities extends also on the singular locus, we can extend
as well the notion of {\it Cauchy surfaces}. These turn to be spacelike
with conical singularities. Such a cone spacetime is said {\it
globally hyperbolic} if it contains a Cauchy surface. Similarly to the
smooth case, we can restrict our study  to {\it maximal globally
hyperbolic} ones. More precisely we set
$$\Mm\Gg\Hh_\kappa(S,\Aa)$$ the Teichm\"uller-like space of 
cone spacetime structures $h$ on $\hat S\times\mr$, of constant curvature
$\kappa\in\{-1,0,1\}$, such that:
\smallskip

- $h$ is non singular on $S\times \mr$; 
\smallskip

- $h$ has a timelike geodesic line of conical singularity of angle $\alpha_i$
at each $\{p_i\}\times\mr$;
\smallskip

- $h$ is a maximal globally hyperbolic and has a Cauchy surface orthogonal
to the singular set;

- these structures are considered up to isotopies of $\hat S\times\mr$
preserving $V\times\mr$.
\medskip

By easily adapting the constructions of Section~\ref{WR} we get, for
every $\kappa=0, \pm 1$, the map
$$\mG_\kappa: \Mm\Ll_\cG (S,\Aa) \to \Mm\Gg\Hh_\kappa(S,\Aa) \ .$$ The
main differences are that even for $\kappa = 0,-1$ the developing maps
are no longer embeddings; moreover, also the asymptotic complex
projective structures produced by the Wick rotations have conical
singularities.  As in the smooth case, the maps $\mG_\kappa$ are
injective. For the so obtained spacetimes have cosmological time, and
one can recover the corresponding data $(F,\lambda)$ by looking at the
level surfaces of the cosmological time.  Moreover, by construction,
canonical Wick rotations and rescalings, directed by the gradient of
the cosmological times, with the usual universal rescaling functions,
apply to the spacetimes belonging to the images of the maps
$\mG_\kappa$.
\smallskip

On the other hand, the question of having an {\it intrinsic
characterization of the images} Im$(\mG_\kappa)$ is largely open. In
particular one asks to determine, for every $\kappa$, the angle
assignments $\Aa$ such that $\mG_\kappa$ gives a parametrization of
the whole of $\Mm\Gg\Hh_\kappa(S,\Aa)$ (possibly inverting the time
orientation). We have (see~\cite{BS})
\begin{prop}\label{smallangle} If all the cone angles are less than $\pi$, 
then the spacetimes belonging to {\rm Im}$(\mG_{-1})$ are precisely
those admitting a {\rm convex} Cauchy surface orthogonal to the
singular locus.
\end{prop}
In fact, under such a ``big masses'' hypothesis, being in the image of
$\mG_{-1}$ turns to be equivalent to admit a \emph{convex core}, that
is a minimal convex subset. The convex core is homeomorphic to
$\hat S\times\mr$ and its boundary is the union of two $\mathrm
C^{0,1}$-spacelike, intrinsically hyperbolic bent cone surfaces
$\partial_+\Kk(Y)$ and $\partial_-\Kk(Y)$, orthogonal to the singular
locus. Just like the non-singular case, $Y$ is encoded by
$(F,\lambda)$ iff the future boundary of its convex core is
obtained by bending $F$ along $\lambda$; similarly for the past
boundary, via $T$-symmetry.

One would expect that for big masses, the map $\mG_{-1}$ actually is a
bijection, that is a convex Cauchy surface should always exist.
 
If some cone angle is bigger or equal than $\pi$, it is known that in
general the maps $\mG_\kappa$ are not onto, even if all masses are positive. 
For example
in~\cite{BG}(2), by applying a so called ``patchwork'' construction,
one produces flat MGH cone spacetimes with positive masses and with
some cone angles equal to $\pi$, that do not belong to the image of
$\mG_0$. In fact it is remarkable that these spacetimes have
nevertheless cosmological time whose level surfaces are orthogonal to
the singular locus, and are {\it flat} instead of hyperbolic at the
singular points of cone angle $\pi$. The canonical rescalings apply to
them so that we finally get also spacetimes that do not belong to the
images of $\mG_{\pm 1}$.

\subsection{Earthquakes on hyperbolic cone surfaces}
 As every lamination $\lambda\in \Mm\Ll_\cG(F,\Aa)$ avoids the conical points, 
the notion of earthquake along such a lamination is defined as well.
Similarly to the non-singular case, we have (see \cite{BS})
\begin{teo}\label{zzz:teo} If  $Y$  has big masses, belongs to 
${\rm Im}(\mG_{-1})$, and is encoded by $(F,\lambda)\in
\Mm\Ll_\cG(S,\Aa)$, then the left (resp. right) earthquake on $F$ along
$\lambda$ produces surfaces $\beta_L(F,\lambda)$
(resp. $\beta_R(F,\lambda)$) $\in\Tt_\cG (S,\Aa)$ whose holonomy coincides
with the right (resp. left) holonomy of $Y$.
\end{teo}

Under the big masses hypothesis, let us consider the map
$$ \mu: {\rm Im}(\mG_{-1})\rightarrow \Tt_\cG (S,\Aa)\times\Tt_\cG (S,\Aa)$$
that associates to every $Y$ the points obtained by left and 
right earthquake on $(F,\lambda)$ respectively, as above. We have (\cite{BS}) 
\begin{teo}\label{wwww:teo}
The following equivalent facts hold
\begin{enumerate}
\item Given $F,F'\in\Tt_\cG (S,\Aa)$ there exists a unique
  $\lambda\in\Mm\Ll_\cG(F,\Aa)$ such that $\beta_L(F,\lambda)=F'$.
\item The map $\mu$ is bijective.
\end{enumerate}
\end{teo}
Notice that the first statement is in purely hyperbolic terms.  The
equivalence between the two statements follows from Theorem
\ref{zzz:teo}. This equivalence between the hyperbolic and Lorentzian
formulations plays an subtle r\^ole in the proof of Theorem
\ref{wwww:teo}.  In fact by means of the hyperbolic formulation the
map $\mu$ is proved to be locally injective, whereas Lorentzian
geometry is used to prove that it is a proper map.
\smallskip

Finally we mention that in Chapter 7 of \cite{Be-Bo}, we have
described a quite different family of spacetimes with cone angles $\geq \pi$
({\it i.e.} possibly with negative masses) that are governed by
quadratic differentials rather than by measured geodesic laminations,
and such that Wick rotation-rescaling machinery does apply to them.

%% file: HANDBIB.tex
%%% Local Variables: 
%%% mode: latex
%%% TeX-master: "MEM2"
%%% End: 